\documentclass[12pt]{amsart}
\usepackage{amsmath,amsthm,amssymb, amsrefs}
\usepackage{a4wide}
\usepackage{fancybox}
\usepackage{maple2e}
\usepackage[backref]{hyperref}%

\theoremstyle{plain} 
\newtheorem*{main theorem}{Main Theorem}
\newtheorem*{theorem*}{Theorem}
\newtheorem*{emptytheorem*}{}
\newtheorem{theorem}{Theorem}
\newtheorem{lemma}{Lemma}[section]
\newtheorem{sublemma}{Sublemma}[lemma]
\newtheorem{proposition}{Proposition}[section]
\newtheorem*{proposition*}{Proposition}

\newtheorem*{conjecture*}{Conjecture}
\theoremstyle{definition}

\newtheorem*{definition*}{Definition}
\theoremstyle{remark}

\newtheorem*{remark*}{Remark}

\newcommand{\nc}{\newcommand}
 \newcommand{\cal}[2]{\ensuremath{\mathcal{#1}^{(#2)}}}
\nc{\tcal}[2]{\ensuremath{\mathcal{\tilde{#1}}^{(#2)}}}
\nc{\hcal}[2]{\ensuremath{\mathcal{\hat{#1}}^{(#2)}}}

\begin{document}
\author{Stefano Luzzatto}
\address{Department of Mathematics, 
Imperial College, London, UK}
\email{stefano.luzzatto@imperial.ac.uk}
\urladdr{http://www.ma.ic.ac.uk/\textasciitilde luzzatto}

\author{Hiroki Takahasi}
\address{Department of Mathematics,  
Kyoto University, Kyoto, Japan}
\email{takahasi@kusm.kyoto-u.ac.jp}
\date{\today} %

\title[Computable conditions for non-uniform hyperbolicity ]
{Computable conditions for the occurrence of  
non-uniform hyperbolicity in families
of one-dimensional maps}

\begin{abstract}
    We formulate and prove a Jakobson-Benedicks-Carleson type theorem on the
    occurrence of nonuniform hyperbolicity (stochastic dynamics) in families of
    one-dimensional maps, based on \emph{computable starting
    conditions} and providing  \emph{explicit, computable,} lower bounds for
    the measure of the set of selected parameters. 
    As a first application of our results we show
    that the set of parameters corresponding to maps 
    in the quadratic family \( f_{a}(x) =
    x^{2}-a \)  which have an absolutely continuous invariant
    probability measure is at least \( 10^{-5000} \).
 \end{abstract}

\maketitle

\section{Introduction}

In this paper we consider families 
\(
\{f_{a}\}_{a\in\Omega}
\)
of \( C^{2} \) interval maps with a single quadratic critical point
and the parameter \( a \) belonging some some interval \( \Omega \).
Families of interval maps can exhibit a wide variety of dynamical 
behaviour, ranging from the existence of attracting periodic 
orbits to the existence of absolutely continuous (with respect to
Lebesgue) invariant measures with strong mixing
properties, as well as all kinds of intermediate and pathological phenomena. 
These dynamical phenomena can depend very sensitively, and very
discontinuously, on the parameter even in very smooth parametrized families,
see \cite{MelStr93} for a comprehensive survey. 
    We say that \( a\in \Omega \) is a \emph{regular} parameter if \( 
    f_{a} \) has an attracting periodic orbit; we say that \( a\in
    \Omega \) is a \emph{stochastic} parameter if \( f_{a} \) admits
    an ergodic invariant probability measure \( \mu \) which is absolutely
    continuous with respect to Lebesgue and has a positive Lyapunov
    exponent 
    \( \lambda(\mu)=\int \log|f'| d\mu \). 
    Then  we define 
\[ 
\Omega^{-}=\{a: a\text{ is regular }\}
\quad
\text{and}\quad
\Omega^{+}=\{a: a\text{ is stochastic}\}.
\]
We remark that the positivity of the Lyapunov exponent is sometimes a 
non-trivial but automatic consequence of the existence of an
absolutely continuous invariant probability measure. This is the case 
for example for maps with a single critical point such as those
considered here \cite{Kel90} and is probably true for more general 
maps with multiple critical points (it is of course false for maps
without critical points such a rigid circle rotations). 
We include it here in the definition since it is a feature which plays a key
role in the giving rise to stochastic behaviour, see
\cites{KelNow92, You92, BruLuzStr03, AlvLuzPindim1, Luz05} for results
and references to results concerning the precise ``random-like'' properties of
such maps.

For families of maps with a single critical point it is also known
that
the sets \( \Omega^{-} \) and \( \Omega^{+} \) are 
disjoint, see \cite{MelStr93}. Moreover, 
for generic families with a quadratic critical point both \(
\Omega^{-} \) and \( \Omega^{+} \) have positive measure and their 
union has full measure \cites{Lyu02, AviMor03, AviMor03b, AviLyuMel03}. 
These are therefore the only two ``typical''
phenomena. The topological structure of the two sets is however very
different:  \(
\Omega^{-} \) is \emph{open and dense} in \( \Omega \) 
\cites{GraSwi97, Lyu97, Koz03, KozSheStr04} and thus \( \Omega^{+} \),
which is contained in the complement of \( \Omega^{-} \), 
is \emph{nowhere dense}.  The fact that it has positive measure is
therefore non-trivial. This was first shown in the
ground-breaking work of Jakobson
\cite{Jak81} and was generalized, over the
years, in several papers; we mention in particular \cites{ColEck83, BenCar85, Ryc88,
MelStr93, ThiTreYou94, Tsu93, Jak01, Yoc01, Jak04} for smooth maps with
non-degenerate critical points, \cite{Thu99} for maps with a degenerate
(flat) critical point, and \cites{LuzTuc99, LuzVia00} for maps with
both critical points and singularities with unbounded derivatives.

The fact that stochastic behaviour occurs for a
nowhere dense set of parameters means that, notwithstanding the fact
that it also occurs with positive probability, it is a difficult set to
actually ``pinpoint'' in practice. In this paper we are concerned
with the problem of estimating \emph{explicitly }the  measure 
of the set \( \Omega^{+} \).
Notwithstanding the extensive amount of research in the area, none of 
the existing results provide any explicit quantitative bounds on the relative
measure of \( \Omega^{+} \) in \( \Omega \). The arguments are
constructive to some extent:  they ``construct'' a set of
stochastic parameters in a given parameter interval \( \Omega \).
However they all apply only under the assumptions that \( \Omega \) is
a \emph{sufficiently small}
neighbourhood of some \emph{sufficiently good} parameter value \( a^{*} \).
Both of these assumptions are problematic in different ways which we 
discuss briefly in the following paragraphs. 

The size of the neighbourhood \( \Omega \) 
of \( a^{*} \) to which the arguments
apply, as well as the proportion of \( \Omega^{+} \) in \(
\Omega \), is not computed in any of the existing arguments.
This is to some extent more of a technical issue than a
conceptual one: explicit estimates can probably be obtained from
the existing papers by more careful control of the interdependence of 
the constants involved. It should be noted nevertheless that this
interdependence is quite subtle; the issue might be technical but this
does not make it non-trivial. 

A more delicate issue is the assumption that the parameter interval \( 
\Omega \) contains some good parameter \( a^{*} \). This gives rise to
the problem of verifying the presence of such a parameter in \( \Omega \)
and of computing the required quantitative information related to 
the ``goodness'' of \( a^{*} \). It turns out that this is essentially
impossible in general.  Several conditions of various kinds have been 
identified which imply that a given parameter value is stochastic 
\cites{Jak78,
Mis81, ColEck83, NowStr91, BruLuzStr03, BruSheStr03, AlvLuzPindim1} but,
apart from some exceptional cases, all these conditions require
information about an \emph{infinite number of iterations}  and thus
are, to all effects and purposes, \emph{uncheckable}.  
It can also be shown that in a formal theoretical sense the set \(
\Omega^{+} \) is \emph{undecidable}, see \cite{ArbMat04}. 

The main objective of this
paper is to overcome these difficulties. We present a 
\emph{quantitative parameter exclusion argument} which gives explicit 
lower bounds on the proportion of \( \Omega^{+} \) in some parameter
intervals \( \Omega \) and, perhaps most importantly, base this
argument on \emph{explicitly computable} starting
conditions which can be computed in finite time and with finite
precision. 
In Section~\ref{explicit} we give an
application  to the quadratic family and 
obtain a first ever explicit 
lower bound for measure of \( \Omega^{+} \). 
In the remaining sections we prove the main theorem.

\section{Statement of results}
\label{s:results}
To simplify the calculations we shall
consider one-parameter families of \( C^{2} \) interval maps of the
form
\[ 
f_{a}(x) = f(x) - a
\]
for some map \( f: I \to I \) with a single quadratic critical point \( 
c \) and the parameter \( a \)
belonging to some interval \( \Omega \). In particular the
critical point \( c \in I \) does not change with the
parameter; this is not an essential condition but simplifies some of 
the already very technical estimates. For the same reason we shall
suppose without loss of generality that the interval \( I \) contains 
strictly the interval \( [-1,1] \), that the critical point \( c=0 \),
and that \( f_{a}(c)>1 \) for all parameters \( a\in \Omega \).  These
conditions always hold up to a linear rescaling which will not affect 
the argument.

We remark that we are thinking of \( \Omega \) as being
``small'' in a sense which will become clear. This is a natural
and general setting: the estimates can  easily be
applied to a ``large'' parameter interval by subdividing it into
 smaller subintervals and considering each such subinterval
independently.  For future reference we let
 \( \varepsilon = |\Omega|  \).  

\subsection{Computable starting conditions}
 
We start by formulating several conditions in terms of  six
constants
\[ N, \delta, \iota, C_{1}, \lambda, \alpha_{0}, \lambda_{0}. \]
First of all we define  critical neighbourhoods
 \[
 \Delta :=(-\delta, \delta) \subset (-\delta^{\iota}, 
 \delta^{\iota})=:\Delta^{+}. 
 \] 
 We shall assume without loss of generality that both \( \log \delta \) and \(
 \iota\log\delta \) are integers.
 
 \textbf{(A1) Uniform expansivity outside a critical
neighbourhood:} 
for all  \( a\in\Omega,  x\in I,  n\geq 1  \)
 such that   \( x, f_a(x),\ldots, f_a^{n-1}(x) \notin \Delta  \)
    we have:
\begin{equation*}
    |(f^{n}_{a})'(x)| \geq 
    \begin{cases} 
C_{1}  e^{\lambda n} &\text{ if }
    f_a^n(x)\in\Delta^+ 
\\
e^{\lambda n} &\text{ if }
x \in f_{a}(\Delta^{+}) \text{ and/or  } f_a^n(x)\in\Delta. 
\end{cases}
\end{equation*}

\textbf{(A2) Random distribution of critical orbits:}
There exists  \( \tilde N\geq N \) such that \( |f_{a}^{n}(c)| \geq \delta^{\iota} \)
for all \( n\leq \tilde N \) and 
\begin{equation}\label{def1}
    |\Omega_{\tilde N}|:=|\{f^{\tilde N}_{a}(c): a \in \Omega\}|\geq \delta^{\iota}. 
\end{equation}

\textbf{(A3) Bounded recurrence of the critical orbit:} 
for all \( a\in \Omega \) and all \( N\geq n \geq 1 \) we have
\[ 
|f_{a}^{n}(c)| \geq  e^{-\alpha_{0}n}.
\]

\textbf{(A4) Non-Resonance:} There exists an integer 
\( \tilde N \geq 1 \)  such
that 
\[ 1+\sum_{i=1}^{k}\frac{1}{(f^{i}_{a})'(c_{0})} \neq 0 
\ \forall \  k\in [1, \tilde N], \text{ and }
 1-\left|\sum_{i=1}^{\tilde N }\frac{1}{(f^{i})'(c_{0})}\right|-
\frac{e^{-\lambda_{0}(\tilde N +1)}}{1-e^{-\lambda_{0}}} > 0. \]

We emphasize that all of these conditions are  \emph{computable}
 in the sense that they can be verified in a finite number of
steps and using explicit numerical calculations relying only on finite
precision and depending only on a finite number of iterations. 
This does not mean however that their verification in practice is
trivial or even easy; we shall discuss below some of the computational issues 
which arise when applying our results to specific situations.

Condition (A1) says that some uniform expansivity estimates hold 
outside the critical neighbourhoods, uniformly for all parameter values. 
This is in some sense the most important condition of all, the general
principle being that \emph{if we have uniform expansivity on a
sufficiently large region of the phase space for all
parameter values, then we have nonuniform expansivity on
all the phase space for a large region of the parameter space.}
For a single parameter value this expansivity is known to hold under
extremely weak conditions \cite{Man85}, although in general the
constants \( C_{1} \) and \( \lambda \) will depend on \( \delta \). 
By continuity it then holds for nearby maps though, again, the size of
an allowed perturbation in general will depend at least on \( \delta \) and \( 
\lambda \). Condition (A1) as stated therefore is not a strong
assumption \emph{per se}, but becomes strong if we want it to hold for
large \( \lambda \), small \( \delta \) and/or a large interval of
parameter values. 
On the other hand it is important to
have a ``large enough'' parameter interval relative to the choice of \( 
\delta, \iota \) for otherwise (A1) could be satisfied even though
all parameters in \( \Omega \) have an attracting periodic orbit (this
may happen, for instance, if the attracting periodic orbit always has 
at least one point in \( \Delta \), then its attracting nature could be
invisible to derivative estimates outside \( \Delta \)). 

The
appropriate condition on the size of \( \Omega \) is given 
in (A2) which says that the size of the interval given by the images
of the critical points for all parameter values at some time is large 
enough. Technically this gives a sufficiently ``random'' distribution 
of the images to provide the first  step of the probabilistic 
induction argument, showing that
 critical orbits have small probability of returning close to the
 critical point. This is, conceptually, the core of the overall
argument and this condition  is essentially implicit in all
arguments of this kind, in the language of Benedicks-Carleson
\cites{BenCar85, BenCar91} and other papers which follows similar
strategies, such as \cites{LuzTuc99, LuzVia00}. Intervals of parameter
values satisfying these conditions are called \emph{escaping
components}. 

Condition (A3) has been used by Benedicks and Carleson and other
people and has proved extremely useful in arguments related to
shadowing the critical point. Notice however that this condition only 
refers to an initial finite number of iterates, as part of the
inductive construction we will guarantee that it continues to hold for
all good parameter values for all time, see Section~\ref{indassu}. We 
do not feel however that it is as essential a condition as (A1) and
(A2), but rather more of a technical simplifying assumption, albeit
one that is not easy to remove in this context. Binding
period arguments analogous to the ones developed here using this
condition have been generalized in \cite{BruLuzStr03} in order to
study the statistical properties of a large class of maps with
several critical points. 

Condition (A4) is easily checked for parameter 
neighbourhoods of particularly good parameter values, such as when the
critical orbit is pre-periodic or non-recurrent.  A similar condition 
played an important role in the work of Tsujii \cite{Tsu93} in generalizing
parameter exclusion arguments to neighbourhood of quite general maps. 
Again however we do feel that this is more of a technical simplifying 
assumption rather than a deep condition. In principle it should be
possible to weaken or possibly even dispense with both (A3) and (A4)
although this would certainly require some non-trivial technical
improvements in the argument.

 \subsection{Conditions on the constants}   \label{constcond}
 Our assumptions on the family \( \{f_{a}\} \) will be that conditions
 (A1)-(A4) hold for a set of constants 
 \( N, \delta, \iota, C_{1}, \lambda, \alpha_{0}, \lambda_{0} \)
 satisfying certain non-trivial formal relationships (C1)-(C4) 
 which we proceed to
 formulate. 
 The first condition is formulated purely in terms of the constants introduced above: 
 \begin{equation*}\label{C1}\tag{\textbf{C1}}
\lambda > \lambda_{0} > \alpha_{0} \geq  \log
    \delta^{-1/N} \quad \text{ and } \quad \log\delta^{-\iota} \geq
    e^{-\frac{1+\lambda_{0}}{2}}. 
    \end{equation*}
 This gives a first constraint on the relative values of the constants.
We remark that the second expression is essentially trivial since generally
speaking \( \log\delta^{-\iota} \) is large and \( \lambda_{0} \) is
small; however it does not appear to be a formal consequence of any of
the other assumptions and we do use it below, therefore we add it here
as a formal assumption. 
 Additional constraints are imposed indirectly via the definition of
   a  set of auxiliary constants in the order given in the following
   list:
\[ 
M_{1}, M_{2}, L_{1}, L_{2}, \alpha_{1}, N_{1}, \mathcal D_{1},\mathcal D_{2}, \mathcal D_{3}, 
\gamma_{0}, \gamma_{1}, 
 \gamma_{2}, \gamma, \hat{D}, \hat{\hat{D}},  \mathcal D, \Gamma_{1}, k_{0},
 \tau_{1}, \tau_{0}, C_{3}, \tilde C_{3},  \tau, \alpha, \tilde \eta, \eta.
\]
As these constants are introduced we shall define conditions 
 which implicitly impose conditions on the original
 constants \( N, \delta, \iota, C_{1}, \lambda, \alpha_{0}, \lambda_{0} \).
The entire procedure is aimed at obtaining a value for the last
constant \( 
\eta \), which appears in the statement of the
Theorem below.  
 
 We remark first that the constants 
\( M_{1}, M_{2}, L_{1}, L_{2}, N_{1},\mathcal D_{2}, \mathcal D_{3} \)
 require some amount of computation concerning the geometry and/or
 the dynamics  of the family \( \{f_{a}\} \). 
In particular they are not defined exactly: they are ``estimates''
which are required to satisfy some
lower or upper bounds.  We shall use the notation \( 
:\geq \) to denote the fact that the constant on the left is 
required to be an upper bound for the expression on the right. 
Similarly for \(
: \leq  \). 
The four constants \( \alpha_{1},\gamma_{1}, \gamma_{2}, \alpha \) 
are chosen with some freedom
within certain ranges depending on the previously chosen constants.
The remaining constants are defined directly in
terms of those defined previously.

We start with fixing constants \( M_{1}, M_{2}, L_{1},L_{2} \) so that
\[
M_{1}:\geq \max\{|f_a'(x)|\}  \quad \text{and} 
\quad M_{2}:\geq \max\{|f_a''(x)|\}, 
\]
where the maximum is taken over all \( (x,a)\in I\times\Omega\), 
and \( L_{1} \) and \( L_{2} \) are chosen such that 
\begin{equation}\label{quadratic0} 
    \begin{cases}
L_1^{-1}|x-c|^{2} \geq |f_a(x)-f_a(c)| \geq L_1 |x-c|^{2}
\\
L_1^{-1}|x-c| \geq
|f_a'(x)| \geq L_1 |x-c| 
\end{cases}\forall x\in \Delta, a\in\Omega 
\end{equation}
and \( L_{2}^{-1}|x-c| \geq 
|f_a'(x)|\geq L_{2} |x-c| \ \forall x\in I, a\in\Omega. \) 
Notice that \( L_{1} \) is used in bounds in \( \Delta \) whereas \(
L_{2} \) applies to the entire interval \( I \). 
We then choose some  
\begin{equation}\label{alpha1}
    \lambda_{0} \geq : \alpha_{1}:\geq \alpha_{0}
 \end{equation}
 and define the following constants. First of all let 
\begin{equation}\label{N1}
    N_{1}:\leq \max\{i: d(f^{i}(\Delta_{0}), c) \geq 1\},
 \end{equation}
where \( d(f^{i}(\Delta_{0}), c) \) denotes the distance 
of iterates of 
\(\Delta_{0}=f(\Delta) \)  from the critical point \( c \), minimized over all
parameters in \( \Omega \). We can assume without loss of generality
that \( N_{1}\geq 1 \), recall the discussion at the beginning of
Section~\ref{s:results}. 
\begin{equation}\label{D1}
\mathcal D_{1}:= \exp\left(
\frac{M_{2}}{L_{2}}\left(
\frac{e^{-\alpha_1}}{1-e^{-\alpha_1}}+
\frac{e^{-(\alpha_1-\alpha_0)(N_1+1)}}{(1-e^{-(\alpha_{1} 
- \alpha_{0})(N_1+1)})(1-e^{-(\alpha_1-\alpha_0)})}\right)\right), 
\end{equation}
\( \mathcal D_{1} \)
 bounds the distortion during
binding periods, see Section~\ref{binding} and Sublemma~\ref{dist1}.
\begin{equation}\label{D2}
   \mathcal D_{2} :\geq 
   \max_{a\in\Omega}\max\left\{
       \max_{1\leq k\leq \tilde N}
       \left\{\left|1+\sum_{i=1}^{k}\frac{1}{(f^{i}_{a})'(c_{0})}\right|
\right\},
   \left|1+
   \sum_{i=1}^{\tilde N}\frac{1}{(f^{i}_{a})'(c_{0})}\right| + 
   \frac{e^{-\lambda_{0}(\tilde N+1)}}{(1-e^{-\lambda_{0}})} 
\right\},
   \end{equation}
   \begin{equation}
   \label{F3}
   \mathcal D_{3}^{-1} :\leq 
   \min_{a\in\Omega} \min\left\{ 
       \min_{1\leq k\leq
   \tilde N}\left\{\left|1+\sum_{i=1}^{k}\frac{1}{(f^{i}_{a})'(c_{0})}\right|
   \right\}, \  
  1-\left|\sum_{i=1}^{\tilde N}\frac{1}{(f^{i})'(c_{0})}\right|-
   \frac{e^{-\lambda_{0}(\tilde N+1)}}{1-e^{-\lambda_{0}}}
   \right\},
   \end{equation}
   recall that \( \tilde N \) is given (with some freedom) by
   condition (A4). The constants \( \mathcal D_{2}, \mathcal D_{3} \) appear in 
the context of the estimates which compare derivatives with respect to
the parameter and derivatives with respect to the space variable, see 
Lemma~\ref{pardist}. Condition (A4) is designed specifically to ensure that \( \mathcal
D_{3}>0 \). The values of \( \mathcal
D_{2} \) and \( \mathcal D_{3} \), or 
the complexity of the calculation, can be optimized
by  choosing different values for \( \tilde N \). 
Then we define
  
       \begin{equation}\label{gamma0}
	 \gamma_{0} :=
	   \frac{2+\log 2 + 5 \log \log\delta^{-1}}{\log\delta^{-1}},
       \end{equation}
The definition of \( \gamma_{0} \) comes from 
the combinatorial estimates in Section~\ref{combest}.       
We choose 
\begin{equation}\label{min}
\min\left\{1- \frac{\log C_{1}^{-1} +
2\log\log\delta^{-\iota})}{\log\delta^{-\iota}}, 1- \gamma_{0}\right\}
>: \gamma_{1} : > 0.
\end{equation}
The first condition here guarantees the convergence of the infinite
sum in Sublemma~\ref{sumr}. The second condition is appears in the
statement of the key Proposition~\ref{main}. 
it also allows us to
choose some 
\[
1-(\gamma_{0}+\gamma_{1})>:\gamma_{2} : > 0 
\quad\text{ 
and let }
\quad
\gamma := \gamma_{0}+\gamma_{1}+\gamma_{2}\in (0,1),
 \]
 \( \gamma_{2} \) and \( \gamma \) both appear
 in the statement of Proposition~\ref{average}. 
 Now we define
 \begin{equation}\label{dhat}
 \hat  D := 2+\frac{2C_{1}^{-1}\mathcal{D}_2\mathcal{D}_3
  e^{-\lambda}}{1-e^{-\lambda}}+\frac{2\mathcal{D}_1
     \mathcal{D}_2\mathcal{D}_3L_1^{-2}}
     {1-e^{-(\alpha_1-\alpha_0)}}, 
\end{equation}
\begin{equation}\label{dhathat}
  \hat{\hat D} = 
  \left[\left(2+e\frac{(\log\delta^{-\iota})^{2}}
  {(\log\delta^{-\iota}-1)^{2}}\right) 
  \frac{(\log\delta^{-\iota})^2} {(\log\delta^{-\iota})^2 - 
  C_{1}^{-1}\delta^{\iota (1-\gamma_1) }}
  \right] \frac{1}{\log\delta^{-\iota}-1},
\end{equation}  
\( \hat D \) and \( \hat{\hat D} \) are simply shortcuts to long
expressions. We then let
\begin{equation}\label{D}
\mathcal{D} := \mathcal D_{2}\mathcal D_{3}
\exp\left( \frac{M_{2}}{L_{2}} \left(\hat D \hat{\hat D} + 
\frac{C_{1}^{-1}\mathcal{D}_2\mathcal{D}_3 e^{-\lambda}}
{(1-e^{-\lambda})}\right)\right),
\end{equation}
\( \mathcal D \) is the global distortion bound, see  Lemma~\ref{distortion}.
      \begin{equation}\label{Gamma1}
\Gamma_{1}:=\mathcal D\mathcal D_{1}\mathcal D_{2}\mathcal
		  D_{3}e^{1+\lambda_{0}}/L_{1}C_{1},
	  \end{equation} 
  the constant \( \Gamma_{1} \) is a large ``dummy'' 
  constant used to formulate the estimates on the derivative growth at
  the end of binding periods, see \eqref{unibindder}. 
  It allows us to cancel out the small constant arising in the
  estimate \eqref{metric0}. The ``payback'' 
  for having such large \( \Gamma_{1} \) will be in terms of a lower
  bound on \( \gamma_{1} \), see below. 	  
  \begin{equation}\label{k0}
	      k_{0} := \max\left\{\frac{\log (\mathcal 
  D_{1}/L_{1})+\lambda_{0}+\alpha_{1}}{\log\delta^{-\iota}}, 0\right\}, 
  \quad\quad	 
  \tau_{1}:=\frac{2}{\lambda_{0}+\alpha_{1}} 
 \leq
 \frac{2+k_{0}}{\lambda_{0}+\alpha_{1}}=:\tau_{0},     
\end{equation}	  
\( k_{0} \) is of no particular significance, it 
is just used in the definition of \( \tau_{0} \); its definition
come from the proof of Sublemma~\ref{binding1}, see \eqref{pbound}.
 \( \tau_{0}, \tau_{1} \) appear in the binding
period estimates, Sublemmas~\ref{binding1} and~\ref{bindder}. 
We can now formulate the next condition
\begin{equation*}\label{C2}\tag{\textbf{C2}}
       \tau_{0}\alpha_{0}<1.
\end{equation*} 
The product \( \alpha_{0}\tau_{0} \) is related to the length of the
binding periods, see \eqref{up}, Lemma~\ref{uniformbind}. 
Then we define
\begin{equation}\label{c3}
    C_{3}:=
       \mathcal D_{1}^{-\frac{\lambda_{0}+2\alpha_{1}}{\lambda_{0}+\alpha_{1}}} 
       L_{1}^{2+\frac{\alpha_{1}}{\lambda_{0}+\alpha_{1}}}
   \quad\text{ and }\quad
     \tilde C_{3} := \frac{2^{\alpha_{1}\tau_{1}-1}L_{1}^{2} }
     {2\mathcal D \mathcal D_{1}^{2}\mathcal D_{2}\mathcal D_{3}}
     C_{3},
      \end{equation}
 \( C_{3} \) and \( \tilde C_{3} \)  have no
      particular significance, they are just a shorthand way of
      writing some complex expressions arising in  the paper.
      They       are first used in 
      \eqref{binddereq} and in \eqref{bexp} respectively.
      We require
 \begin{equation}\label{C3}\tag{\textbf{C3}}
  \max\left\{
  \begin{aligned}
  &\iota+\frac{
\log (\mathcal D \mathcal D_{2} \mathcal 
  D_{3} C_{1}^{-1})+ 2\log\log\delta^{-\iota} }{\log\delta^{-1}}
  \\
  &\alpha_{1}\tau_{1}+ 
  \frac{\log(\Gamma_{1}\mathcal D_{2}\mathcal 
   D_{3}\tilde C_{3}^{-1} e^{\alpha_{1}\tau_{1}-1})+2\log 
   \log \delta^{-\iota}}{\iota\log\delta^{-1}}
  \end{aligned}
  \right\} \leq \gamma_{1}.
  \end{equation}
  These conditions 
  arise in \eqref{metric4} and \eqref{beexp4} respectively.  
\begin{equation}\label{tau}
     \tau :=
     \frac{\tau_{0}}{1-\gamma_{1}} 
     \left(1+\frac{\log 
     |I|+ 2\log\log\delta^{-1}  - \log \Gamma_{1} }{\log\delta^{-1}}\right) > 
     0. 
 \end{equation}
 The constant \( \tau \) gives the total proportion of iterates
 which belong to binding periods, see Lemma~\ref{inessential}. 
 We remark that it is not clear directly from the definition that \(
 \tau >0 \). However this follows \emph{a fortiori} from estimates
 \eqref{eq:tau} and \eqref{eq:p} at the end of the proof of Sublemma~\ref{inessential}. 
 \begin{equation}\label{alpha}
     \alpha:= \min\left\{\alpha_{0}, \frac{\lambda-\lambda_{0}}{
(\tau(\lambda-\frac{1-\gamma_{1}}{\tau_{0}})+1)}\right\} > 0.
 \end{equation}
 The requirement \eqref{alpha} on \( \alpha \)
 comes from the growth estimates in the proof 
 of Lemma~\ref{posexp}. Notice that the fact that \( \alpha>0 \) comes from 
\(
\lambda-\frac{1-\gamma_{1}}{\tau_{0}}\geq 
\lambda-(1-\gamma_{1})(\frac{\lambda_{0}+\alpha_{1}}{2}) 
\geq \lambda - (1-\gamma_{1}) \lambda_{0} \geq \lambda-\lambda_{0}>0. 
\)
The first inequality comes from the definition of \( \tau_{0}
\), the second from the fact that \( \alpha_{1}<\lambda_{0} \), the 
the third from the fact that \( \gamma_{1} <1 \), the fourth from (C1).
 \[ 
\tilde\eta:= e^{-\gamma_{2}
 \alpha} \left(1+ \sum_{R\geq \log \delta^{-1}} 
 e^{ - (1-\gamma) R}\right) = 
 e^{-\gamma_{2} \alpha} 
 \left(1+\frac{\delta^{1-\gamma}}{1-e^{-(1-\gamma)}}\right),
\]
\( \tilde\eta \) appears in Subsection~\ref{deviations}. 
 We can then formulate our  final condition which includes the
 definition of the main constant which appears in the statement of the
 Main Theorem. 
\begin{equation}\label{C4}\tag{\textbf{C4}}
\eta :=
\sum_{j=N}^{\infty} \tilde\eta^{j}  = \frac{\tilde\eta^{N}}{1-\tilde\eta}< 1.
\end{equation}

\subsection{The main result}
We are now ready to state our main result. First we
 define the set
 \[
 \Omega^*:=\{a\in\Omega:|(f_a^n)'(c_0)|\geq e^{\lambda_0 n}\ 
 \forall n\geq 0\}.
 \]
By classical results, 
 for any \( a\in \Omega^{*} \), \( f_{a} \) admits an ergodic
 absolutely continuous invariant probability measure and thus \(
 \Omega^{*}\subseteq \Omega^{+} \). 
 
 \begin{main theorem}
 \label{main theorem} 
 Let \( f_{a}: I \to I \) be a family of \( C^{2} \) unimodal maps with
 a quadratic critical point, of
 the form \( f_{a}(x) = f(x)-a \), for \( a \) belonging to some
 interval \( \Omega \) of parameter values. 
 Suppose that there exist constants 
 \(
 N, \delta, \iota, C_{1},  \lambda, \alpha_{0}, \lambda_{0} 
 \)
 so that conditions (A1)-(A4) and (C1)-(C4) hold. 
Then
 \[
 | \Omega^*| \geq (1-\eta) |\Omega|. 
 \]
 \end{main theorem}
 In Section~\ref{explicit} we give an application to the
 quadratic family and  in Sections~\ref{indassu} to \ref{s:posmeas} 
 we give a complete proof. 
 First of all we make some general remarks; 
 
 \subsubsection{Computer-assisted proofs}
The spirit of this result is to have a ready-made ``formula'' for
proving rigorously the existence, 
and for obtaining a lower bound for the probability, 
of stochastic dynamics in a given family of maps. 
In some sense we think of this result as a 
\emph{computer-assisted} 
theorem in (non-uniformly) hyperbolic dynamics. 
The combination of rigorous numerical/computational methods 
with deep geometric/analytic/probabilistic methods can be
extremely powerful with each approach contributing to overcome the
limitations of the other. Indeed there exist extremely powerful
techniques in Dynamical Systems and Ergodic Theory which yield highly 
sophisticated results about the fine structure and long term behaviour
of certain systems, but these methods often rely on certain geometric 
assumptions, usually hyperbolicity assumptions, 
which are highly non trivial to verify in practice. On the
other hand numerical/computational methods are much more flexible and 
in principle can be applied to essentially any system but suffer from 
the fundamental limitation of having finite precision and being able
to deal with only finite time properties. These two approaches are
therefore naturally complementary and the research presented in this
paper is precisely a step in the direction of implementing this point 
of view. 

There are of course several other  results in Dynamics, some of which 
have proved ground-breaking, which have
benefited from such a combination of methods, although, as far as we
know, not in the area of non-uniform hyperbolicity. We mention here
just some of these: the pioneering work, also in the context of 
families of interval maps, of Collet, Eckmann and Lanford on the
hyperbolicity of the renormalization operator, \cite{ColEckLan80,
Lan82}; the development of computational implementations of ideas from 
Conley Index theory, see \cite{ZglNis01, Mis02} and in particular the
book \cite{KacMisMro04} for a comprehensive review and references. 
The history of the study of the Lorenz attractor \cite{Lor63} is particularly
interesting from this point of view: over more than 20 years
sophisticated geometrical/analytical methods have been applied to a
``geometric model'' of the attractor, essentially assuming some
hyperbolicity conditions of the equations. Recently these conditions
have been verified by Tucker \cite{Tuc99, Tuc02} using a combination
of numerical and analytic methods and this the abstract results
developed for the geometric model immediately apply to the actual
Lorenz equations. 

Finally, we mention also the extremely interesting
area of so-called \emph{constructive KAM theory} which aims to give 
 a more concrete and quantitative understanding of the classical
 abstract results of KAM theory (see \cite{Lla04} for a comprehensive
survey and references). Although KAM theory deals in many ways with
situations which are the complete opposite of the theory of
hyperbolicity, there are some striking structural similarities such as
the occurrence of certain dynamical phenomena for 
nowhere dense sets  of positive Lebesgue measure. Analogously to the
situation in KAM theory, 
the first twenty years of work in non-uniform hyperbolicity has been
relatively abstract and the present work can perhaps be seen as an initial
step in the direction of an analogous \emph{constructive theory of
non-uniform hyperbolicity} providing a concrete and quantitative 
understanding of the abstract results.

\subsubsection{Computational issues}
Though we have stressed the computability of our assumptions, and
indeed give a concrete application of our result,  
we emphasize that a systematic and extensive 
 application of our results involves several non-trivial theoretical
 and computational issues. 
 
 The first, relatively obvious, issue is the development of
 rigorous computational algorithms to verify conditions (A1)-(A4). 
 Of these conditions, (A1) is probably the most technically and
 computationally demanding requiring some non-trivial algorithm. 
 However, even for condition (A1) and also for the other three
 conditions, the problem essentially boils down to having sufficient
 precision in the calculations and thus does not seem to present
 particular conceptual hurdles. Rather, it is likely that issues 
 regarding the `computational cost'' and ``efficiency'' 
 of the algorithms will come into play, since in general 
 the number of iterations which need to
 be computed will be quite large. 
 
 A second, less obvious but probably ultimately more
 difficult, issue is that of finding a systematic way of adjusting the 
 constants \( N, \delta, \iota, C_{1}, \lambda, \alpha_{0}, \lambda_{0} \)
 as well as the length of the subinterval \( \Omega \) of parameters
 and the auxiliary constants, so that all the conditions are satisfied. 
 A natural situation would be that conditions (A1)-(A4) have been
 verified for a certain set of constants 
 \( N, \delta, \iota, C_{1}, \lambda, \alpha_{0}, \lambda_{0} \) but
 that these constants fail one of the conditions (C1)-(C4). One then
 has many options such as changing any one or 
 more of the constants and trying 
 to verify (A1)-(A4) and (C1)-(C4) again. The interdependence of the
 various conditions is however quite subtle, non-monotone, 
 and trial and error does
 not appear to be a satisfactory strategy in general. 
 It would therefore be essential to find a systematic way of optimizing
 the choices in order to converge to a set of constants in which all
 conditions are satisfied.

\subsubsection{Notation}
We introduce some notation which will be used extensively below. 
We let \(c_{0}=c_{0}(a)= f_{a}(0) \) denote the 
\emph{critical value} 
of \( f_{a} \) and for \( i\geq 0 \), 
\( c_{i}=c_{i}(a) = f^{i} (c_{0})
\). 
A key feature of the argument involves considering a family of
maps from parameter space to dynamical space which tracks the orbits
of the critical points for different parameter values:
for \( n\geq 0 \) and \(
\omega\subseteq\Omega \) let 
\[
\omega_{n}=\{c_{n}(a); a\in \omega\}
\subseteq I.
\]

\subsection{Acknowledgements}
We would like to thank Hiroshi Kokubu and Konstantin Mischaikow for
the fundamental role they have played in the the origin and 
development of this project. The idea of a computational approach to
the parameter exclusion argument has been around probably as long as
the argument itself; Benedicks and Jakobson at least 
have been conscious of such a potential and certainly other 
people would have thought about it. As
far as we are concerned however, the stimulus for this project came
from Mischaikow's suggestion and his enthusiasm and conviction
regarding its feasibility. Mischaikow and Kokubu went on to organize
reading groups on the topic in Atlanta and Kyoto which provided 
important practical motivation. In particular, 
HT would like to express special gratitude to his supervisor,
Hiroshi Kokubu, for his encouragement and support. 

We also thank Artur Avila and Marcelo Viana for pointing out an error 
in an earlier version of the paper, Katie Bloor, Simon Cedervall,
Karla Diaz-Ordaz, Oliver Butterley for their valuable help and
comments during the final  proof-reading of the paper, and an
anonymous referee for their support, encouragement, and general
enthusiasm for the project.

\section{Explicit estimates for the quadratic family}
\label{explicit}

A special, but
already very interesting, example of a family of interval maps
such as those considered above, is the well known
quadratic family
\[ 
f_{a}(x) = x^{2}-a.
\]
A relatively straightforward application of our Main Theorem gives us 
the following
\begin{theorem}\label{quadratictheorem}
    Let \( f_{a}(x) = x^{2}-a \) be the quadratic
    family for 
    \(
  a\in \Omega := [2-10^{-4990}, 2].  
    \) 
 Then 
   \[ 
  |\Omega^{+}| \geq 0.97 |\Omega| > 10^{-5 000}.
   \]
\end{theorem}
The proof actually yields a lower bound for the set 
\[
\Omega^{*}=\{a\in \Omega: 
|(f_a^n)'(c_0)|\geq e^{0.61 n}\ \forall n\geq
0\}\subseteq\Omega^{+}.
\]  
By classical results all these maps are stochastic. 
We emphasize that even though the interval of
parameter values under consideration here is relatively small, this is
the first known rigorous explicit lower bound for the set \( \Omega^{+} \)
in any context. The only other related work we are aware of is by
Sim\'o and Tatjer \cites{SimTat91, SimTat05}, who carried out careful 
numerical estimates of the length of some of the \emph{periodic
windows}, i.e. the open connected components of \( 
\hat\Omega^{-} \) in \( \hat\Omega= [\hat a, 2]  \), where 
\( \hat a \) is the Feigenbaum parameter 
at the limit of the first period-doubling cascade
They are able to compute the
length of some relatively ``large'' windows (up to those of length of 
the order of \( 10
^{-30} \) or so, though not all the windows of up to this size) and 
their calculations suggest a
lower bound for the overall proportion of \( \hat\Omega^{-} \) in \(
\hat\Omega \) of about 10\%. One expects the contribution of the
smaller windows to be negligible but this is probably a hard statement
to prove. 

A combination of the numerical
methods of Sim\'o and Tatjer with the application of the estimates
given here might form a good strategy for obtaining some global 
bounds for the relative 
measures of \( \hat\Omega^{-} \) and \( \hat\Omega^{+} \) 
in the entire parameter interval \( \hat\Omega \). Indeed, our methods
appear naturally suited to the analysis of small parameter intervals
and thus,in some sense, take over precisely where the purely numerical
estimates no longer work. 
Another possible strategy might involve taking advantage of the 
hyperbolicity of the renormalization operator in combination, although
even if such an
argument were feasible, this strategy would probably give less
concrete information about the actual location of parameters of \(
\hat\Omega^{-} \) and \( \hat\Omega^{+} \).  

We prove Theorem~\ref{quadratictheorem} in the next three sections. 
In Section~\ref{choosing} we discuss the choice of the main
constants and check conditions (A1) to (A4); in Section 
\ref{conditionsC} we discuss the choice of the auxiliary constants 
 and check conditions (C1) to (C4); 
 in Section~\ref{maple} we enclose a copy of the 
 Maple worksheet used to carry out the explicit 
 numerical calculations. Once the calculations are set up it is very
 easy to obtain the estimates for other intervals
 of parameter intervals. In particular, choosing smaller intervals
 generally yields smaller values for \( \eta \) and thus a larger
 proportion of stochastic parameters.

\subsection{Two preliminary lemmas}
In this section we prove two preliminary lemmas which motivate our
choice of the constants in the following subsection. 
 
 \begin{lemma}\label{quadraticlambda}
  Let  \( 1 \geq \delta, \iota > 0 \) and \( 2 \geq  a^{*}>1 \). Then, for any 
 \(  a\in [a^{*}, 2] \) the expansivity condition (A1) 
 holds with constants 
 \[ 
 C_{1}= \sqrt{\frac{4-\delta^{2\iota}}{4-\delta^{2}}}
 \quad\text{and}\quad
 \lambda=
 \log \left(2\delta
 \sqrt{\frac{4-\delta^{2}}{4-(\delta^{2}-a^{*})^{2}}}\right).
 \]
 \end{lemma}

 \begin{proof}
 Let \( h: (0,1) \to (-2, 2) \) be given by 
 \( x = h(\theta) = 2 \cos \pi \theta.  \)
 For each \(  a\in [1,2] \) let  \( \zeta_{a}= \pi^{-1}\cos^{-1}(\sqrt{(2+a)/4}) \),
 and consider the family of intervals 
 \(
 J_{a}=  (\zeta_{a},1-\zeta_{a}) \subset [0,1]  
\) and their images \( I_{a}= h(J_{a}) = (-\sqrt{2+a}, \sqrt{2+a}) 
\).
 For \( a=2 \) we have \( f_{2}(I_{2}) = I_{2} \), and for
 \( a\in [0,2)  \) we have \( f_{a}(I_{a}) \subset I_{a} \), since \( f_{a}(0) =
 -a > -\sqrt{2+a}\). Thus we shall always consider \( f_{a}:I_{a}\to
 I_{a} \) to be the restriction of \( f_{a} \) to \( I_{a} \), and 
implicitly consider \( h \) restricted to \( J_{a} \). We then define a family of maps 
 \( g_{a}: J_{a} \to J_{a}\) 
  by 
 \begin{equation}\label{g}
  g_{a}(\theta) = h^{-1}\circ f_{a} \circ h (\theta) = 
 \frac{1}{\pi}\cos^{-1}\left( 2\cos^{2}\pi\theta - \frac{a}{2}
 \right).
 \end{equation}
 Thus \( h:J_{a}\to I_{a} \) defines a smooth conjugacy between \( f_{a} \)
 and \( g_{a} \) and, for any \( n\geq 1 \), we have 
 \(
 f_{a} ^{n}= h \circ g_{a}^{n} \circ h^{-1}.
 \)
 Therefore, for \( x= h(\theta) \) we have 
 \begin{equation}\label{Df}
     \begin{split}
 D_{x}f_{a}^{n}(x) &= Dh (g^{n}_{a}(h^{-1}(x))) \cdot
 D g_{a}^{n}(h^{-1}(x)) \cdot D h^{-1}(x)
 \\ &= 
 \frac{Dh (g^{n}_{a}(\theta)) }{- Dh(\theta)} 
 D g_{a}^{n}(\theta).
 \\ &= \frac{\sin \pi (g^{n}_{a}(\theta)) }{\sin \pi (\theta)} 
 D g_{a}^{n}(\theta) 
 \quad\text{ since } D_{\theta}h(\theta) = -2\pi \sin \pi \theta.
 \end{split}
 \end{equation}
Differentiating  \eqref{g} 
 and using the fact that \( \theta = \pi^{-1}\cos^{-1}(x/2) \), 
 as well as the  standard identity \(
 \sin(\cos^{-1}(x)) = \sqrt{1-x^{2}} \),  
 to get,  for \( a^{*}\in (1, 2], a\in [a^{*}, 2] \)  and \( x\notin\Delta
 \),  
 \begin{equation*}
 \left| D_{\theta}g_{a}(\theta)\right| = 
 \left| \frac{4\cos \pi\theta \sin\pi\theta }{  
 \sqrt{1-\left(
 2\cos^{2}\pi\theta - \frac{a}{2}
 \right)^{2}}} \right| 
 = \left| 2x
 \sqrt{\frac{4-x^{2}}{4-(x^{2}-a)^{2}}} \right| 
 \geq 2\delta
 \sqrt{\frac{4-\delta^{2}}{4-(\delta^{2}-a^{*})^{2}}}= \lambda.
 \end{equation*}
The last inequality uses the fact that  \( |D_{\theta}g_{a}(\theta)| \)
 is monotone increasing in \( |x| \) and \( a \). 
Moreover we have 
\[ 
 \frac{\sin \pi (g^{n}_{a}(\theta)) }{\sin \pi (\theta)} 
 = 
 \frac{\sin \pi (h^{-1}(f^{n}(x))) }{\sin \pi (h^{-1}(x))} 
 = 
 \frac{\sin (\cos^{-1}(f^{n}(x)/2)) }{\sin (\cos^{-1}(x/2))} 
 = 
 \sqrt{\frac{4-(f^{n}(x))^{2}}{4-(x)^{2}}}.
 \]
 If \( |x|\geq \delta \) and \(
  |f^{n}(x)| \leq \delta^{\iota} \), substituting into 
  \eqref{Df}, we get the first part of (A1). Similarly 
 if  \(  |f^{n}_{a}(x)| \leq |x| \), e.g. if \( x\in
 f_{a}(\Delta^{+}) \),
 and/or  \( f_{a}^{n}(x) \in \Delta \) we have the second part of 
(A1). 
\end{proof}

 Notice that for \( a=2 \) we have
 \(
 |D_{\theta}g_{2}(\theta)| 
 \equiv 2 
 \) for all \( \theta \neq 1/2 \) or \( x\neq 0 \) and thus in
 particular  \( |D_{x}f_{2}^{n}(x)| \geq 2^{n} \text{ whenever } 
 |f^{n}_{2}(x)| \leq |x|
 \), 
 which is in itself a quite remarkable and non-trivial estimate for \( f_{2} \).
It can also  be checked directly from \eqref{g} that \( g_{a} \) is the standard
 ``top'' tent map, 
 thus \( h \) defines a smooth conjugacy between the ``top'' quadratic 
 map and the ``top'' tent map.  Indeed this was the main idea used by 
 Ulam and Von Neumann in their paper \cite{UlaNeu47} to prove that \( 
 f_{2} \) admits an ergodic absolutely continuously invariant
 probability measure. Similar calculations to the ones carried out
 above have been used in several papers to prove similar estimates
 concerning the expansion outside critical neighbourhoods for maps
 close to \( f_{2} \). 
 For \( a<2 \), 
 \( g_{a} \) is no longer piecewise linear and has a critical 
 point at \( \theta = 1/2 \) (corresponding to \( x =
 h^{-1}(\theta) = 0 \)).  However, it has the advantage over the maps 
 in the quadratic family that its fold near \( \theta=1/2 \) is
 extremely ``sharp'' and, outside a small neighbourhood of \( 1/2 \)
 the slope remain essentially close to \( \log 2 \), and even tends to
 infinity close to the boundaries of the intervals of definition \(
 J_{a} \). 

Notice also that the Lemma, as stated, does not specify any relations 
between \( \delta \) and \( a^{*} \). Indeed, for a fixed \( a^{*}<2
\), sufficiently small values of \( \delta \) give negative values of \( 
\lambda \). This corresponds to the fact that such a parameter may
have an attracting periodic orbit (recall that an open and dense set
of parameters have attracting periodic orbits). If at least one point 
of this orbit lies inside \( \Delta \) then the attracting nature of
the orbit is invisible to derivative estimates outside \( \Delta \),
in other words the existence of such an orbit is not incompatible to
expansion estimates as in condition (A1) with \( \lambda >0 \).
However, if \( \delta \) is so small that all points of the attracting
orbit lie outside \( \Delta \) then clearly (A1) can only be satisfied
from some \( \lambda < 0 \).

    \begin{lemma} \label{quadratick}
	For all \( a\in \Omega \)  we have
     \[
     f^{n}_{a}(c_{0})
      > 1.5 > \delta^{\iota} \quad \text{for all }\quad 1\leq n \leq 
     N  := \frac{1}{\log 4}\log\frac{1}{
(2-(a^{*})^{2}+a^{*})}.
     \]
     \end{lemma}
     \begin{proof}
For each parameter value \( a \),  \( f_{a} \) has two fixed points, one
of which is given by \( p_{a}=(1+\sqrt{1+4a})/2 \). For \( a \) close 
to \( 2 \) the second image of the critical point lies close to this
fixed point. Indeed \( c_{0}(a^{*}) = f_{a^{*}}(c) = -a^{*} \) 
and  \(c_{1}(a)= 
f_{a }(c_{0}(a)) = (c_{0}(a))^{2}-a = a^{2}-a \). For \( a=2 \)
we have \( c_{1}(2) = 2 = p_{2} \) and the critical point actually
lands on the fixed point, for \( a<2 \) and close to 2 we have 
\( a^{2}-a <  (1+\sqrt{1+4a})/2 < 2\) and thus \( c_{1}(a) \) lies to the
\emph{left} of the fixed point \( p_{a} \), which itself lies to the
left of the fixed point \( p_{2}=2 \) for \( f_{2} \). 

Let \( I_{1}=[c_{1}(a), p_{a}] \) and \( I_{n} = f_{a}^{n-1}(I_{1}) = 
[c_{n}(a), p_{a}]\).
Then a simple application of the Mean Value Theorem, using the fact
that the derivative is bounded above by \( 4 \), gives 
\[ 
|I_{n}|\leq 4^{n-1}|I_{1}| \leq 4^{n-1} (2-a^{2}+a)
\]
for every \( n\geq 1 \). We have used here the fact that the length of
\( I_{1} \) is the distance between \( c_{1}(a) \) and \( p_{a} \)
which is less than the distance between \( c_{1}(a) \) and \( 2 \),
which is precisely given by \( 2-c_{1}(a) = 2-(a^{2}-2) \). Since 
\( p_{a} \) is very close to \( 2 \) for \( a\in \Omega \) it is
sufficient to show that 
\begin{equation}\label{kN}
4^{n-1} (2-a^{2}+a) \leq 1/4.
\end{equation}
In particular that \( f^{n-1} \) is a bijection from \(
I_{1} \) to \( I_{n} \) as long as \eqref{kN} holds, and thus this
implies in particular \( c_{n}(a)\geq 1/2 \) as required. 
Taking logs and solving for \( n \) gives exactly the choice of \( N \)
above.

\end{proof}

\subsection{Condition (A1) to (A4)}
\label{choosing}

In this section we discuss the choice of constants in our specific
setting and show that all the necessary conditions are satisfied. 
As per the statement of the Theorem, we have first of all
\[ 
a^{*} = 2-10^{-4990}. 
\]
We remark that the exponent 4990 is  chosen in order to have a
round figure for the lower bound for the measure of stochastic
parameters after the exclusions have been carried out. 
We now fix 
the values of \( \delta, \iota \) as 
\[ 
\delta = 10^{-1000} \quad\text{and}\quad \iota = 0.8
\]
For these choices we have the following 
\begin{proposition}\label{conditions}
    Conditions (A1) to (A4) are satisfied for constants 
\[ 
    C_{1}:= \sqrt{\frac{4-\delta^{2\iota}}{4-\delta^{2}}}
    \simeq 1, \quad 
    \lambda :=
    \log \left(2\delta
    \sqrt{\frac{4-\delta^{2}}{4-(\delta^{2}-a^{*})^{2}}}\right)
    \simeq 0.693,
    \]
    \[ 
    N:= \left[\frac{1}{
    \log 4(2-(a^{*})^{2}+a^{*})}\right]
    = 8287,
    \]
    \[
    \alpha_{0} := \frac{\log \delta^{-1}}{N}\simeq 0.277, \quad 
    \lambda_{0} := 0.2 \alpha_{0}+0.8\lambda \simeq 0.610.
    \]
\end{proposition}
Proposition~\ref{conditions} says in particular that all constants
can be chosen as formal functions of the three variable \( a^{*},
\delta, \iota \). We shall show below that the same is true also for
the other auxiliary functions.
In the statement we have given the approximate 
values of the constants for our particular choices of \( a^{*}, \delta \) and \( 
\iota \). 

\begin{proof}
Condition (A1) follows immediately from Lemma~\ref{quadraticlambda}. 
For (A2) let \( \tilde N \) 
be the smallest integer for  which 
\[
    |\Omega_{\tilde N}|=|f^{\tilde N}_{a^{*}}(c_{0})-2|\geq  1/4 > \delta^{\iota}.
\]
 For (A3) we just observe  that Lemma~\ref{quadratick} and the
 the choice of \( \alpha_{0}= \log\delta^{-1}/N \) imply 
 \[
 |f^{n}_{a}(c_{0})| > 1 > \delta = e^{\log \delta} \geq 
 e^{-\frac{n}{N}\log \delta^{-1}}   = e^{-\alpha_{0} n} 
 \quad \text{for } n\leq N.
 \]
 It remains only to verify (A4), for which it suffices to choose \(
 N_{1}=1 \). Indeed, then we have 
 \(
{e^{-2\lambda_{0}}}/{(1-e^{-\lambda_{0}})} \approx 0.65
 \)
 and so (A4) is easily seen to hold. 
\end{proof}    
    
\subsection{Conditions (C1) to (C4)}
\label{conditionsC}
We proceed to compute the auxiliary constants and prove 
\begin{proposition}
    Conditions (C1) to (C4) are satisfied.
 \end{proposition}  
 We emphasize that all the computations
are carried out in Maple on a \emph{purely symbolic} level 
and thus there is no risk of approximation
errors being amplified over several calculations due to floating point
arithmetic etc.
Maple essentially calculates the
final constant \( \eta \) as a formal function of the previously
defined constants. 
Only at the final stage, the command \texttt{evalf} asks
for this expression to be evaluated numerically. This command can also
be used at the intermediate stages to obtain approximate 
 values arising from the intermediate calculations, and these have
 been shown above and will be shown below but these
 approximations are not used in the calculation of \( \eta \). 
 
 Condition (C1) follows immediately from the calculations given above:
 \begin{equation}\tag{C1}
     \lambda\simeq 0.693 > \lambda_{0} \simeq  0.610  > \alpha_{0} = \log 
     \delta^{-1/N} \simeq 0.277
    \text{ and }
     \log\delta^{-1} - e^{-\frac{1+\lambda_{0}}{2}} \simeq 2302.
     \end{equation}
 To verify the other conditions we need to compute some the other
 constants. First of all  we have 
 \[ 
 M_{1}=M_{2}=4 \quad\text{ and } \quad L_{1}=L_{2} = 0.5.
 \]
Then we choose 
 \[ \alpha_{1}=0.2 \lambda_{0}  + 0.8 \alpha_{0} \simeq 0.344 \]
 We then set 
\[  N_{1}:= 
     \left[ \frac{1}{\log 4}
     \log\frac{1}{2\delta^{2}}\right] = 3321 < N = 8287. 
     \]
To see that this satisfies \eqref{N1} recall that 
from Lemma~\ref{quadratick} we have   \( |f^{i}_{a}(c_{0})| 
\geq 1.5\) for all \( a\in\Omega \) and all \( 1\leq i
\leq N \). For these iterates, by the Mean Value Theorem we have 
\[
|f^{n}(\Delta_{0})| \leq 4^{n} |\Delta_{0}| = 4^{n}\delta^{2}. 
\]
Thus it is sufficient to solve \( 4^{n}\delta^{2}<1/2 \) in terms of \( 
n \) and this gives exactly the condition on \( N_{1} \) above. 
We can then compute 
\[
\mathcal D_{1} := 
\exp\left(
\frac{1}{
1-e^{-\alpha_{1} 
}}
+
 \frac{e^{-(\alpha_{1}-\alpha_{0})(N_{1}+1)
 }}{(1-e^{-(\alpha_{1}-\alpha_{0}) (N_{1}+1)})^{2}}
 \right)\simeq 31.
 \]
 Notice that the formal expression for \( \mathcal D_{1} \) contains a
 term \( M_{2}/L_{2} \) in the exponent. However it is immediate from 
 the proof of Sublemma~\ref{dist1} that the two terms cancel out in
 the case of the quadratic family. The same is also true for the
 overall distortion constant \( \mathcal D \) computed below. 
To compute \( \mathcal D_{2} \) and \( \mathcal D_{3} \), 
we choose \( \tilde N=N \) and use the fact that \(
(f^{i}_{a})'(c_{0})\geq 3^{i}\) for all \( 1\leq i \leq
N\). This gives 
\[
\mathcal D_{2} := \frac{3}{2}
+ \frac{e^{-\lambda_{0}(N+1)}}{1-e^{-\lambda_{0}}}
\geq     
1+\sum_{i=1}^{N}\frac{1}{3^{i}} + 
    \frac{e^{-\lambda_{0}(N+1)}}{1-e^{-\lambda_{0}}}
\geq 
1+\sum_{i=1}^{N}\frac{1}{(f^{i}_{a^{*}})'(c_{0}(a^{*}))} + 
    \frac{e^{-\lambda_{0}(N+1)}}{1-e^{-\lambda_{0}}}.
\]
and similarly
\[ 
\mathcal D_{3}^{-1}:=
\frac{1}{2}
- \frac{e^{-\lambda_{0}(N+1)}}{1-e^{-\lambda_{0}}}
\leq 1-\sum_{i=1}^{N}\frac{1}{3^{i}} - 
    \frac{e^{-\lambda_{0}(N+1)}}{1-e^{-\lambda_{0}}}
\leq 
1-\sum_{i=1}^{N}\frac{1}{(f^{i}_{a^{*}})'(c_{0}(a^{*}))} - 
    \frac{e^{-\lambda_{0}(N+1)}}{1-e^{-\lambda_{0}}}.
\]
This gives values of 
\[ 
\mathcal D_{2} \simeq 1.5 \quad\text{and}\quad \mathcal D_{3} \simeq
2.
\]
By explicit computation we then get 
\[ 
\gamma_{0} \simeq 0.017,
\]
and choose 
\[ 
\gamma_{1} = 0.835  \min\left\{1- \frac{\log C_{1}^{-1} +
2\log\log\delta^{-\iota})}{\log\delta^{-\iota}}, 1- \gamma_{0}\right\} 
\simeq 0.982,
\]
\[ 
\gamma_{2} = 0.8 (1-\gamma_{0}-\gamma_{1}) \simeq 0.117, 
\]
and so
\[
\gamma:= \gamma_{0}+\gamma_{1}+\gamma_{2} \simeq 0.970.
\]
Then we get 
\[ 
\hat D \simeq 11570, \hat{\hat D} \simeq 0.002,  \mathcal D \simeq
10^{14}, \Gamma_{1}\simeq 10^{17},  k_{0} \simeq 0.003, 
\tau_{0}\simeq 2.098, \tau_{1} \simeq 2.095.
\]
We can therefore verify 
\begin{equation}\tag{C2}
    \tau_{0}\alpha_{0} \simeq 0.58  < 1.
    \end{equation}
By more explicit computations we get 
\[ 
C_{3} \simeq 0.002 \quad  \tilde C_{3}\simeq 10^{-22},
\]
and we can verify (C3): 
   \[ 
   \gamma_{1}\simeq 0.835 > 0.821 \simeq   \max\left\{
  \begin{aligned}
  &\iota+\frac{
\log (\mathcal D \mathcal D_{2} \mathcal 
  D_{3} C_{1}^{-1})+ 2\log\log\delta^{-\iota} }{\log\delta^{-1}}
  \\
  &\alpha_{1}\tau_{1}+ 
  \frac{\log(\Gamma_{1}\mathcal D_{2}\mathcal 
   D_{3}\tilde C_{3}^{-1} e^{\alpha_{1}\tau_{1}-1})+2\log 
   \log \delta^{-\iota}}{\iota\log\delta^{-1}}
  \end{aligned}
  \right\}. 
   \]
We then compute 
\[ 
\tau \simeq 12, \quad 
\alpha  \simeq 0.009,  \quad \tilde\eta \simeq 0.998,  
\]
 and verify 
 \begin{equation}\tag{C4}
 \eta = \frac{\tilde\eta^{N}}{1-\eta} \simeq 0.082 < 1. 
 \end{equation}
This  gives the value of \( \eta \) which appears in the
statement of Theorem~\ref{quadratictheorem}.

\subsection{Maple computations}
\label{maple}

We enclose a copy of the Maple worksheet used to carry out the
explicit calculations. The notation should be self-explanatory. 

\def\emptyline{\vspace{12pt}}
\DefineParaStyle{Maple Heading 4}
\DefineParaStyle{Maple Heading 2}
\DefineParaStyle{Maple Text Output}
\DefineParaStyle{Maple Bullet Item}
\DefineParaStyle{Maple Warning}
\DefineParaStyle{Maple Error}
\DefineParaStyle{Maple Dash Item}
\DefineParaStyle{Maple Heading 3}
\DefineParaStyle{Maple Heading 1}
\DefineParaStyle{Maple Title}
\DefineParaStyle{Maple Normal}
\DefineCharStyle{Maple 2D Input}
\DefineCharStyle{Maple Maple Input}
\DefineCharStyle{Maple 2D Output}
\DefineCharStyle{Maple 2D Math}
\DefineCharStyle{Maple Hyperlink}
\begin{mapleinput}
\mapleinline{active}{1d}{\begin{Maple Normal}\QTR{Maple Maple Input}{}\QTR{Maple Maple Input}{restart;}\end{Maple Normal}}{}
\end{mapleinput}
\begin{mapleinput}
\mapleinline{active}{1d}{\begin{Maple Normal}\QTR{Maple Maple Input}{}\QTR{Maple Maple Input}{delta:=10\symbol{94}(-1000): iota:=.8:epsilon:=10\symbol{94}(-4990):}\QTR{Maple Maple Input}{astar:= 2-epsilon:}\end{Maple Normal}}{}
\end{mapleinput}
\begin{mapleinput}
\mapleinline{active}{1d}{\begin{Maple Normal}\QTR{Maple Maple Input}{}\QTR{Maple Maple Input}{C1:=sqrt((4-delta\symbol{94}(2*iota))/(4-delta\symbol{94}2)):evalf[15](C1);}\end{Maple Normal}}{}
\end{mapleinput}

\mapleresult
\begin{maplelatex}
\QTR{Maple 2D Output}{\mapleinline{inert}{2d}{1.000000000}{%
\[
 1.0
\]
}
}

\end{maplelatex}
\begin{mapleinput}
\mapleinline{active}{1d}{\begin{Maple Normal}\QTR{Maple Maple Input}{}\QTR{Maple Maple Input}{lambda:=ln(2*delta*sqrt((4-delta\symbol{94}2)/}\QTR{Maple Maple Input}{
(4-(delta\symbol{94}2-astar)\symbol{94}2))):evalf[15](lambda);}\end{Maple Normal}}{}
\end{mapleinput}

\mapleresult
\begin{maplelatex}
\QTR{Maple 2D Output}{\mapleinline{inert}{2d}{.693147180559945}{%
\[
 0.693147180559945
\]
}
}

\end{maplelatex}
\begin{mapleinput}
\mapleinline{active}{1d}{\begin{Maple Normal}\QTR{Maple Maple Input}{}\QTR{Maple Maple Input}{N:=floor((ln(1/(2-(astar\symbol{94}2-astar))))/ln(4)):evalf[15](N);}\end{Maple Normal}}{}
\end{mapleinput}

\mapleresult
\begin{maplelatex}
\QTR{Maple 2D Output}{\mapleinline{inert}{2d}{8287.}{%
\[
 8287.0
\]
}
}

\end{maplelatex}
\begin{mapleinput}
\mapleinline{active}{1d}{\begin{Maple Normal}\QTR{Maple Maple Input}{}\QTR{Maple Maple Input}{alpha0:=-ln(delta)/N:evalf[15](alpha0);}\end{Maple Normal}}{}
\end{mapleinput}

\mapleresult
\begin{maplelatex}
\QTR{Maple 2D Output}{\mapleinline{inert}{2d}{.277855085434300}{%
\[
 0.277855085434300
\]
}
}

\end{maplelatex}
\begin{mapleinput}
\mapleinline{active}{1d}{\begin{Maple Normal}\QTR{Maple Maple Input}{}\QTR{Maple Maple Input}{lambda0:=(alpha0*.2+ lambda*.8):evalf[15](lambda0);}\end{Maple Normal}}{}
\end{mapleinput}

\mapleresult
\begin{maplelatex}
\QTR{Maple 2D Output}{\mapleinline{inert}{2d}{.610088761524283}{%
\[
 0.610088761524283
\]
}
}

\end{maplelatex}
\begin{mapleinput}
\mapleinline{active}{1d}{\begin{Maple Normal}\QTR{Maple Maple Input}{}\QTR{Maple Maple Input}{evalf[15](log(delta\symbol{94}(-1/N)));}\end{Maple Normal}}{}
\end{mapleinput}

\mapleresult
\begin{maplelatex}
\QTR{Maple 2D Output}{\mapleinline{inert}{2d}{.277855085434300}{%
\[
 0.277855085434300
\]
}
}

\end{maplelatex}
\begin{mapleinput}
\mapleinline{active}{1d}{\begin{Maple Normal}\QTR{Maple Maple Input}{}\QTR{Maple Maple Input}{NR:=(exp(-2*lambda0))/(1-exp(-lambda0)):evalf[15](NR);}\end{Maple Normal}}{}
\end{mapleinput}

\mapleresult
\begin{maplelatex}
\QTR{Maple 2D Output}{\mapleinline{inert}{2d}{.646331222657726}{%
\[
 0.646331222657726
\]
}
}

\end{maplelatex}
\begin{mapleinput}
\mapleinline{active}{1d}{\begin{Maple Normal}\QTR{Maple Maple Input}{}\QTR{Maple Maple Input}{C1b:=ln(delta\symbol{94}(-1))-exp(-(1+lambda0)/2):evalf[15](C1b);}\end{Maple Normal}}{}
\end{mapleinput}

\mapleresult
\begin{maplelatex}
\QTR{Maple 2D Output}{\mapleinline{inert}{2d}{2302.13802490916}{%
\[
 2302.13802490916
\]
}
}

\end{maplelatex}
\begin{mapleinput}
\mapleinline{active}{1d}{\begin{Maple Normal}\QTR{Maple Maple Input}{}\QTR{Maple Maple Input}{M1:=4:M2:=2:L1:=.5:L2:=.5:J:=4:}\end{Maple Normal}}{}
\end{mapleinput}
\begin{mapleinput}
\mapleinline{active}{1d}{\begin{Maple Normal}\QTR{Maple Maple Input}{}\QTR{Maple Maple Input}{alpha1:=(lambda0*.2+alpha0*.8):evalf[15](alpha1);}\end{Maple Normal}}{}
\end{mapleinput}

\mapleresult
\begin{maplelatex}
\QTR{Maple 2D Output}{\mapleinline{inert}{2d}{.344301820623978}{%
\[
 0.344301820623978
\]
}
}

\end{maplelatex}
\begin{mapleinput}
\mapleinline{active}{1d}{\begin{Maple Normal}\QTR{Maple Maple Input}{}\QTR{Maple Maple Input}{N1:=floor(min(log((1-delta\symbol{94}iota)/delta\symbol{94}2)}\QTR{Maple Maple Input}{
/log(4), N-1 )):evalf[15](N1);}\end{Maple Normal}}{}
\end{mapleinput}

\mapleresult
\begin{maplelatex}
\QTR{Maple 2D Output}{\mapleinline{inert}{2d}{3321.}{%
\[
 3321.0
\]
}
}

\end{maplelatex}
\begin{mapleinput}
\mapleinline{active}{1d}{\begin{Maple Normal}\QTR{Maple Maple Input}{}\QTR{Maple Maple Input}{D1:=exp(((1-exp(-alpha1))\symbol{94}(-1) + }\QTR{Maple Maple Input}{
(exp(-(N1+1)*(alpha1-alpha0))/ }\QTR{Maple Maple Input}{
(1-exp(-(N1+1)*(alpha1-alpha0)))*}\QTR{Maple Maple Input}{
( 1-exp(-(alpha1-alpha0)))))):evalf[15](D1);}\end{Maple Normal}}{}
\end{mapleinput}

\mapleresult
\begin{maplelatex}
\QTR{Maple 2D Output}{\mapleinline{inert}{2d}{30.9713855467660}{%
\[
 30.9713855467660
\]
}
}

\end{maplelatex}
\begin{mapleinput}
\mapleinline{active}{1d}{\begin{Maple Normal}\QTR{Maple Maple Input}{}\QTR{Maple Maple Input}{D2:=1.5 + ((exp(-lambda0*(N+1)))}\QTR{Maple Maple Input}{
/(1-exp(-lambda0))):evalf[15](D2);}\end{Maple Normal}}{}
\end{mapleinput}

\mapleresult
\begin{maplelatex}
\QTR{Maple 2D Output}{\mapleinline{inert}{2d}{1.5}{%
\[
 1.5
\]
}
}

\end{maplelatex}
\begin{mapleinput}
\mapleinline{active}{1d}{\begin{Maple Normal}\QTR{Maple Maple Input}{}\QTR{Maple Maple Input}{D3minus:=0.5-((exp(-lambda0*(N+1)))/(1-exp(-lambda0))):}\QTR{Maple Maple Input}{
D3:=D3minus\symbol{94}(-1):evalf[15](D3);}\end{Maple Normal}}{}
\end{mapleinput}

\mapleresult
\begin{maplelatex}
\QTR{Maple 2D Output}{\mapleinline{inert}{2d}{2.00000000000000}{%
\[
 2.0
\]
}
}

\end{maplelatex}
\begin{mapleinput}
\mapleinline{active}{1d}{\begin{Maple Normal}\QTR{Maple Maple Input}{}\QTR{Maple Maple Input}{gamma0:=(1+ln(2) + 5*ln(-ln(delta)) )/}\QTR{Maple Maple Input}{
(-ln(delta)):evalf[15](gamma0);}\end{Maple Normal}}{}
\end{mapleinput}

\mapleresult
\begin{maplelatex}
\QTR{Maple 2D Output}{\mapleinline{inert}{2d}{0.175464029210645e-1}{%
\[
 0.0175464029210645
\]
}
}

\end{maplelatex}
\begin{mapleinput}
\mapleinline{active}{1d}{\begin{Maple Normal}\QTR{Maple Maple Input}{}\QTR{Maple Maple Input}{gamma1max:=min(1-gamma0, }\QTR{Maple Maple Input}{
1-( (ln((C1\symbol{94}(-1))+2*ln(ln(delta\symbol{94}(-iota))))}\QTR{Maple Maple Input}{
/ln(delta\symbol{94}(-iota))))):}\QTR{Maple Maple Input}{
evalf[15](gamma1max);}\end{Maple Normal}}{}
\end{mapleinput}

\mapleresult
\begin{maplelatex}
\QTR{Maple 2D Output}{\mapleinline{inert}{2d}{.982453597078936}{%
\[
 0.982453597078936
\]
}
}

\end{maplelatex}
\begin{mapleinput}
\mapleinline{active}{1d}{\begin{Maple Normal}\QTR{Maple Maple Input}{}\QTR{Maple Maple Input}{gamma1:=0.85*gamma1max:evalf[15](gamma1);}\end{Maple Normal}}{}
\end{mapleinput}

\mapleresult
\begin{maplelatex}
\QTR{Maple 2D Output}{\mapleinline{inert}{2d}{.835085557517095}{%
\[
 0.835085557517095
\]
}
}

\end{maplelatex}
\begin{mapleinput}
\mapleinline{active}{1d}{\begin{Maple Normal}\QTR{Maple Maple Input}{}\QTR{Maple Maple Input}{gamma2:=.8*(1-(gamma0+gamma1)):evalf[15](gamma2);}\end{Maple Normal}}{}
\end{mapleinput}

\mapleresult
\begin{maplelatex}
\QTR{Maple 2D Output}{\mapleinline{inert}{2d}{.117894431649472}{%
\[
 0.117894431649472
\]
}
}

\end{maplelatex}
\begin{mapleinput}
\mapleinline{active}{1d}{\begin{Maple Normal}\QTR{Maple Maple Input}{}\QTR{Maple Maple Input}{gammatil:=gamma0+gamma1+gamma2:evalf[15](gammatil);}\end{Maple Normal}}{}
\end{mapleinput}

\mapleresult
\begin{maplelatex}
\QTR{Maple 2D Output}{\mapleinline{inert}{2d}{.970526392087632}{%
\[
 0.970526392087632
\]
}
}

\end{maplelatex}
\begin{mapleinput}
\mapleinline{active}{1d}{\begin{Maple Normal}\QTR{Maple Maple Input}{}\QTR{Maple Maple Input}{Dhat:=2+((2*C1\symbol{94}(-1)*D2*D3*exp(-lambda)) / }\QTR{Maple Maple Input}{
(1-exp(-lambda)))+ ((2*D1*D2*D3*(L1\symbol{94}(-2)))/}\QTR{Maple Maple Input}{
(1-exp(-(alpha1-alpha0))) ):}\QTR{Maple Maple Input}{
evalf[15](Dhat);}\end{Maple Normal}}{}
\end{mapleinput}

\mapleresult
\begin{maplelatex}
\QTR{Maple 2D Output}{\mapleinline{inert}{2d}{11570.3753180732}{%
\[
 11570.3753180732
\]
}
}

\end{maplelatex}
\begin{mapleinput}
\mapleinline{active}{1d}{\begin{Maple Normal}\QTR{Maple Maple Input}{}\QTR{Maple Maple Input}{logdiota:=log(delta\symbol{94}(-iota)):logdiota2:=logdiota\symbol{94}2:}\end{Maple Normal}}{}
\end{mapleinput}
\begin{mapleinput}
\mapleinline{active}{1d}{\begin{Maple Normal}\QTR{Maple Maple Input}{}\QTR{Maple Maple Input}{Dhathat:=(  (2+exp(1)*(logdiota2/(logdiota-1)\symbol{94}2) )*}\QTR{Maple Maple Input}{
(logdiota2/(logdiota2-C1\symbol{94}(-1)*delta\symbol{94}(iota*(1-gamma1))))}\QTR{Maple Maple Input}{
*(logdiota-1)\symbol{94}(-1) ):}\QTR{Maple Maple Input}{
evalf[15](Dhathat);}\end{Maple Normal}}{}
\end{mapleinput}

\mapleresult
\begin{maplelatex}
\QTR{Maple 2D Output}{\mapleinline{inert}{2d}{0.256440032673570e-2}{%
\[
 0.00256440032673570
\]
}
}

\end{maplelatex}
\begin{mapleinput}
\mapleinline{active}{1d}{\begin{Maple Normal}\QTR{Maple Maple Input}{}\QTR{Maple Maple Input}{Dist:=D2*D3*exp((Dhat*Dhathat+  }\QTR{Maple Maple Input}{
((C1\symbol{94}(-1)*D2*D3*exp(-lambda))/( 1-exp(-lambda))))):}\end{Maple Normal}}{}
\end{mapleinput}
\begin{mapleinput}
\mapleinline{active}{1d}{\begin{Maple Normal}\QTR{Maple Maple Input}{}\QTR{Maple Maple Input}{evalf[15](Dist);}\end{Maple Normal}}{}
\end{mapleinput}

\mapleresult
\begin{maplelatex}
\QTR{Maple 2D Output}{\mapleinline{inert}{2d}{463434666796857.}{%
\[
 463434666796857.0
\]
}
}

\end{maplelatex}
\begin{mapleinput}
\mapleinline{active}{1d}{\begin{Maple Normal}\QTR{Maple Maple Input}{}\QTR{Maple Maple Input}{Gamma1:=Dist*D1*D2*D3*exp(1+lambda0)/(L1*C1):}\end{Maple Normal}}{}
\end{mapleinput}
\begin{mapleinput}
\mapleinline{active}{1d}{\begin{Maple Normal}\QTR{Maple Maple Input}{}\QTR{Maple Maple Input}{evalf[15](Gamma1);}\end{Maple Normal}}{}
\end{mapleinput}

\mapleresult
\begin{maplelatex}
\QTR{Maple 2D Output}{\mapleinline{inert}{2d}{0.430876756737302e18}{%
\[
 430876756737302000.0
\]
}
}

\end{maplelatex}
\begin{mapleinput}
\mapleinline{active}{1d}{\begin{Maple Normal}\QTR{Maple Maple Input}{}\QTR{Maple Maple Input}{k0:=max((ln(D1/L1)+lambda0+alpha1)}\QTR{Maple Maple Input}{
/ln(delta\symbol{94}(-iota)),0):evalf[15](k0);}\end{Maple Normal}}{}
\end{mapleinput}

\mapleresult
\begin{maplelatex}
\QTR{Maple 2D Output}{\mapleinline{inert}{2d}{0.275809649265672e-2}{%
\[
 0.00275809649265672
\]
}
}

\end{maplelatex}
\begin{mapleinput}
\mapleinline{active}{1d}{\begin{Maple Normal}\QTR{Maple Maple Input}{}\QTR{Maple Maple Input}{tau1:=2/(lambda0+alpha1):evalf[15](tau1);}\end{Maple Normal}}{}
\end{mapleinput}

\mapleresult
\begin{maplelatex}
\QTR{Maple 2D Output}{\mapleinline{inert}{2d}{2.09557809691856}{%
\[
 2.09557809691856
\]
}
}

\end{maplelatex}
\begin{mapleinput}
\mapleinline{active}{1d}{\begin{Maple Normal}\QTR{Maple Maple Input}{}\QTR{Maple Maple Input}{tau0:=(2+k0)/(lambda0+alpha1):evalf[15](tau0);}\end{Maple Normal}}{}
\end{mapleinput}

\mapleresult
\begin{maplelatex}
\QTR{Maple 2D Output}{\mapleinline{inert}{2d}{2.09846800021815}{%
\[
 2.09846800021815
\]
}
}

\end{maplelatex}
\begin{mapleinput}
\mapleinline{active}{1d}{\begin{Maple Normal}\QTR{Maple Maple Input}{}\QTR{Maple Maple Input}{evalf[15](1-tau0*alpha0);}\end{Maple Normal}}{}
\end{mapleinput}

\mapleresult
\begin{maplelatex}
\QTR{Maple 2D Output}{\mapleinline{inert}{2d}{.416929994518240}{%
\[
 0.416929994518240
\]
}
}

\end{maplelatex}
\begin{mapleinput}
\mapleinline{active}{1d}{\begin{Maple Normal}\QTR{Maple Maple Input}{}\QTR{Maple Maple Input}{C3:=D1\symbol{94}(-(lambda0+2*alpha1)/}\QTR{Maple Maple Input}{
(lambda0+alpha1))* L1\symbol{94}(2+(alpha1/(lambda0+alpha1))):}\QTR{Maple Maple Input}{
evalf[15](C3);}\end{Maple Normal}}{}
\end{mapleinput}

\mapleresult
\begin{maplelatex}
\QTR{Maple 2D Output}{\mapleinline{inert}{2d}{0.182183306992661e-2}{%
\[
 0.00182183306992661
\]
}
}

\end{maplelatex}
\begin{mapleinput}
\mapleinline{active}{1d}{\begin{Maple Normal}\QTR{Maple Maple Input}{}\QTR{Maple Maple Input}{C3til:=2\symbol{94}(alpha1*tau1-1)*L1\symbol{94}(2)*C3/}\QTR{Maple Maple Input}{
(2*(D1\symbol{94}2)*Dist*D2*D3):evalf[15](C3til);}\end{Maple Normal}}{}
\end{mapleinput}

\mapleresult
\begin{maplelatex}
\QTR{Maple 2D Output}{\mapleinline{inert}{2d}{0.140784264658478e-21}{%
\[
{ 1.40784264658478\times 10^{-22}}
\]
}
}

\end{maplelatex}
\begin{mapleinput}
\mapleinline{active}{1d}{\begin{Maple Normal}\QTR{Maple Maple Input}{}\QTR{Maple Maple Input}{evalf[10](alpha1*tau1);}\end{Maple Normal}}{}
\end{mapleinput}

\mapleresult
\begin{maplelatex}
\QTR{Maple 2D Output}{\mapleinline{inert}{2d}{.7215113540}{%
\[
 0.7215113540
\]
}
}

\end{maplelatex}
\begin{mapleinput}
\mapleinline{active}{1d}{\begin{Maple Normal}\QTR{Maple Maple Input}{}\QTR{Maple Maple Input}{gamma1min:=max(iota+(((ln(Dist*D2*D3/C1) +  }\QTR{Maple Maple Input}{
2*ln(ln(delta\symbol{94}(-iota) )))/(ln(delta\symbol{94}(-1))))), }\QTR{Maple Maple Input}{
alpha1*tau1+((  }\QTR{Maple Maple Input}{
ln(Gamma1*D2*D3*C3til\symbol{94}(-1)*exp(alpha1*tau1-1))+ }\QTR{Maple Maple Input}{
2*ln(ln(delta\symbol{94}(-iota)))))/}\QTR{Maple Maple Input}{
(iota*ln((delta)\symbol{94}(-1)))):evalf[15](gamma1min);}\end{Maple Normal}}{}
\end{mapleinput}

\mapleresult
\begin{maplelatex}
\QTR{Maple 2D Output}{\mapleinline{inert}{2d}{.821673721125574}{%
\[
 0.821673721125574
\]
}
}

\end{maplelatex}
\begin{mapleinput}
\mapleinline{active}{1d}{\begin{Maple Normal}\QTR{Maple Maple Input}{}\QTR{Maple Maple Input}{evalf[15](gamma1-gamma1min);}\end{Maple Normal}}{}
\end{mapleinput}

\mapleresult
\begin{maplelatex}
\QTR{Maple 2D Output}{\mapleinline{inert}{2d}{0.134118363915209e-1}{%
\[
 0.0134118363915209
\]
}
}

\end{maplelatex}
\begin{mapleinput}
\mapleinline{active}{1d}{\begin{Maple Normal}\QTR{Maple Maple Input}{}\QTR{Maple Maple Input}{tau:=(tau0/(1-gamma1))*}\QTR{Maple Maple Input}{
(1+((ln(J)-ln(Gamma1))/ln(delta\symbol{94}(-1))) +  }\QTR{Maple Maple Input}{
(2*ln(ln(delta\symbol{94}(-iota))))/ln(delta\symbol{94}(-1))):}\QTR{Maple Maple Input}{
evalf[15](tau);}\end{Maple Normal}}{}
\end{mapleinput}

\mapleresult
\begin{maplelatex}
\QTR{Maple 2D Output}{\mapleinline{inert}{2d}{12.5909564129903}{%
\[
 12.5909564129903
\]
}
}

\end{maplelatex}
\begin{mapleinput}
\mapleinline{active}{1d}{\begin{Maple Normal}\QTR{Maple Maple Input}{}\QTR{Maple Maple Input}{alpha:=min(alpha0,(lambda-lambda0)/(tau*}\QTR{Maple Maple Input}{
(lambda-((1-gamma1)/tau0))+1)):evalf[15](alpha);}\end{Maple Normal}}{}
\end{mapleinput}

\mapleresult
\begin{maplelatex}
\QTR{Maple 2D Output}{\mapleinline{inert}{2d}{0.950554901234048e-2}{%
\[
 0.00950554901234048
\]
}
}

\end{maplelatex}
\begin{mapleinput}
\mapleinline{active}{1d}{\begin{Maple Normal}\QTR{Maple Maple Input}{}\QTR{Maple Maple Input}{etatil:=exp(-gamma2*alpha)*(1+((delta\symbol{94}(1-gammatil))/}\QTR{Maple Maple Input}{
(1-exp(-1+gammatil)))):evalf[15](etatil);}\end{Maple Normal}}{}
\end{mapleinput}

\mapleresult
\begin{maplelatex}
\QTR{Maple 2D Output}{\mapleinline{inert}{2d}{.998879976396842}{%
\[
 0.998879976396842
\]
}
}

\end{maplelatex}
\begin{mapleinput}
\mapleinline{active}{1d}{\begin{Maple Normal}\QTR{Maple Maple Input}{}\QTR{Maple Maple Input}{eta:=(etatil\symbol{94}(N))/(1-etatil):evalf[15](eta);}\end{Maple Normal}}{}
\end{mapleinput}

\mapleresult
\begin{maplelatex}
\QTR{Maple 2D Output}{\mapleinline{inert}{2d}{0.827085919317526e-1}{%
\[
 0.0827085919317526
\]
}
}

\end{maplelatex}

\section{Setting up the induction}
\label{indassu} 
The remaining sections of the paper are devoted to proving the main
Theorem. 
In this Section~\ref{indassu} we give the formal inductive construction of the set 
\( \Omega^{*}\). In Section~\ref{binding} we prove the main technical 
lemma concerning the shadowing of the critical orbit. In Section 
\ref{s:posexp} we prove the inductive step in the definition of \( 
\Omega^{*} \) and in Section~\ref{s:posmeas} we obtain the lower bound 
on \( |\Omega^{*}| \). 

\subsection{Inductive assumptions}

Let \( \Omega^{(0)}=\Omega\) and \( 
\mathcal{P}^{(0)}=\{\Omega^{(0)}\} \)
denote
the trivial partition of
\( \Omega \). Given \( n\geq 1 \) suppose that for each \( k \leq n-1
\) there exists a set \( \Omega^{(k)}\subseteq\Omega \) satisfying the
combinatorial structure described in the next paragraph and satisfying
the following four bulleted properties.

There exists a partition \(\mathcal {P}^{(k)}\)
of \( \Omega^{(k)} \) into intervals such that each
\( \omega\in\mathcal P^{(k)} \) has an associated \emph{itinerary}
constituted by the following information 
(for the moment we describe the combinatorial structure 
as abstract data, the geometrical meaning of this data will become
clear in the next section).
To each \( \omega\in\mathcal P^{(k)} \) is associated a 
sequence \(
0=\theta_{0}<\theta_{1}<\dots<\theta_{r}\leq k\),  \( r=r(\omega)\geq 
0 \)
of \emph{escape times}. Escape times are divided into three 
categories, i.e.
{\it substantial}, {\it essential}, and {\it inessential}.
Inessential escapes possess no combinatorial feature and 
are only relevant to 
the analytic bounded distortion argument to be 
developed later.
Substantial and essential escapes play a role in splitting
itineraries into  segments in the following sense.
Let \(
0=\eta_{0}<\eta_{1}<\dots<\eta_{s}\leq k\), \(s=s(\omega)\geq 0 \)
be the maximal sequence of substantial and essential escape times. 
Between any of the two \( \eta_{i-1} \)
and \( \eta_{i} \) (and between \( \eta_{s} \) and \( k \)) there is a
sequence \( \eta_{i-1}<\nu_{1}<\dots<\nu_{t}<\eta_{i}\), \( 
t=t(\omega, i)\geq 0 \) of \emph{essential return times} (or
\emph{essential returns}) and between any two essential returns \(
\nu_{j-1} \) and \( \nu_{j} \) (and between \( \nu_{t} \) and \(
\eta_{i} \)) there is a sequence 
\(\nu_{j-1}<\mu_{1}<\dots<\mu_{u}<\nu_{j}\),
\( u=u(\omega, i, j)\geq
0\) of \emph{inessential return times} (or \emph{inessential
returns}).  Following essential and inessential return (resp. escape)
there is a time interval \( [\nu_{j}+1, \nu_{j}+ p_{j}]
\ \text{ (resp.  \(
[\mu_{j}+1, \mu_{j}+ p_{j}] \) ) } \) with \( p_{j} > 0 \) called the
\emph{ binding period} \label{pagebinding}. A binding period cannot
contain any return and escape times.  
Finally, associated to each
essential and inessential return time (resp. escape) is a positive 
integer \( r \)
called the \emph{return depth} (resp. \emph{escape depth}).  

\subsection*{\textbullet \ Bounded Recurrence}\label{BR}
We define the function 
\(
\mathcal {E}^{(k)}: \Omega^{(k)}\to \mathbb N 
\) 
which associates to each \(
a\in\Omega^{(k)} \) the total sum of all essential return depths of
the element \( \omega\in\mathcal {P}^{(k)} \) containing \( a \)
in its itinerary up to and including time \( k \). Notice that 
\( \mathcal {E}^{(k)} \) is constant on elements of \( 
\mathcal{P}^{(k)} \) by construction. Then, for all
\( a\in \Omega^{(k)} \) \begin{equation*}
\tag*{\( (BR)_{k} \)} \mathcal {E}
^{(k)}(a) \leq \alpha k. 
\end{equation*}

\subsection*{\textbullet \ Slow Recurrence}
For all \( a\in \Omega^{(k)} \) and all \( i\leq k \) we have
\begin{equation*}\label{SR}
\tag*{\( (SR)_{k} \)} 
|c_i(a)| \geq e^{-\alpha_{0} i}.
\end{equation*}

\subsection*{\textbullet \ Hyperbolicity}
For all 
\( a\in \Omega^{(k)} \) \begin{equation}\tag*{\( (EG)_{k} \)}
|(f_{a}^{k+1})'(c_{0})|\geq  e^{\lambda_{0} (k+1)}.  
\end{equation}

\subsection*{\textbullet \ Bounded Distortion}
Critical orbits with the same
combinatorics satisfy uniformly comparable derivative estimates:
For every  \( \omega\in\mathcal  {P}^{(k)} \), every pair of 
parameter values \( a,b\in\omega \) and every \( j\leq \nu+p+1 \) 
where \( \nu \) is
the last return or escape before or equal to time \( k \) and \( p
\) is the associated binding period, we have 
\begin{equation}\tag*{\( (BD)_{k} \)}\label{bd}
\frac{|(f_{a}^{j})'(c_{0})|}{|(f_{b}^{j})'(c_{0})|} \leq \mathcal D
\quad \text{ and } \quad \frac{|c'_{j}(a)|}{|c_{j}'(b)|}\leq \mathcal
D.
\end{equation}
Moreover if $k$ is a 
substantial escape 
a similar distortion estimate holds for all $ j \leq l$ ($l$ is the 
next  chopping time) replacing $\mathcal D$ by $\tilde{\mathcal D}$ and
\( \omega \) by any subinterval \( \omega'\subseteq\omega \) which
satisfies \( \omega'_l \subseteq\Delta^{+} \). 
In particular  for \( j\leq k \), the map \( c_{j}: \omega \to 
\omega_{j} = \{c_{j}(a) : a\in \omega \} \) is a bijection.

\subsection{Definition of \( \Omega^{(n)}\) and  
\( \mathcal P^{(n)} \)}
\label{omega n}  

For \( r\in\mathbb N \), let 
\( I_{r}= [e^{-r}, e^{-r+1}),   I_{-r}=-I_{r}.  \) Then
\[ 
\Delta^{+}=\{0\}\cup \bigcup_{|r|\geq r_{\delta^{+}}+1} I_{r}
\quad\text{ and }\quad \Delta=\{0\}\cup \bigcup_{|r|\geq r_{\delta}+1}
I_{r}. 
\] 
where \( r_{\delta}= \log\delta^{-1}, 
r_{\delta^{+}}= \iota\log\delta^{-1} \). Recall that we can suppose
without loss of generality that \( r_{\delta},
r_{\delta^{+}}\in\mathbb N \).  Subdivide each \( I_{r} \)
into \( r^2\) subintervals of equal length. This defines partitions \(
\mathcal I, \mathcal I^{+} \) of \( \Delta^{+}\) with \( \mathcal I=
\mathcal I^{+}|_{\Delta} \). An interval belonging to either one of
these partitions is of the form \( I_{r, m} \) with \( m\in [1, r^2]
\). Let \( I_{r, m}^{\ell} \) and \( I_{r,m}^{\rho} \) denote the
elements of \( \mathcal I^{+} \) adjacent to \( I_{r,m} \) and let \(
\hat I_{r,m}= I^{\ell}_{r,m}\cup I_{r,m}\cup I_{r,m}^{\rho} \). If \(
I_{r,m} \) happens to be one of the extreme subintervals of \(
\mathcal I^{+} \) then let \( I_{r,m}^{\ell} \) or \( I_{r,m}^{\rho}
\), depending on whether \( I_{r,m} \) is a left or right extreme,
denote the intervals 
\( \left(-\delta^{\iota}-\frac{\delta^{\iota}}{(\log\delta^{-\iota})^{2}}, 
-\delta^{\iota}\right] \) or \(
\left[\delta^{\iota}, \delta^{\iota}+
\frac{\delta^{\iota}}{(\log\delta^{-\iota})^{2}}\right) \) respectively. 
We now use this partition to define a
refinement \(\hat{\mathcal {P}}^{(n)}\) of \( \mathcal {P}^{(n-1)} \). 
Let \( \omega\in\mathcal
P^{(n-1)} \).  We distinguish two different
cases. 

\subsection*{ \textbullet \ Non-chopping times}
    We
say that \( n \) is a non-chopping time for \( \omega\in\mathcal
P^{(n-1)}  \) if one (or
more) of the following situations occur: 
(1) \(
\omega_{n}\cap \Delta^{+}=\emptyset \); 
(2) \( n \) belongs to the
binding period associated to some return or escape time \( \nu<n \) of
\( \omega
\); 
(3)  \( \omega_{n}\cap\Delta^{+}\neq \emptyset \) but \(
\omega_{n} \) does not intersect more than two elements of the
partition \( \mathcal I^{+} \). 
In all three cases we
let \( \omega\in \hat{\mathcal {P}}^{(n)} \). 
In cases (1) and (2) and in case \( (3) \) if \( \omega_{n} \) is not 
completely contained in \( \Delta^{+} \), 
no additional combinatorial information is added 
to the itinerary of \( \omega \). Otherwise, i.e. in case (3) with \( 
\omega_{n}\subset \Delta^{+} \), we consider two distinct
possibilities: 
 if \( \omega_{n}\cap (\Delta\cup I_{\pm
r_{\delta}})\neq \emptyset \) we say that \( n \) is
an inessential return time for 
\( \omega\in \hat{\mathcal {P}}^{(n)}\), 
if 
$\omega_n\subset\Delta^+\setminus 
(\Delta\cup I_{\pm r_{\delta}})$ we say that \( n \) is 
an inessential escape time for 
\( \omega\in \hat{\mathcal {P}}^{(n)}\). 
We define
the corresponding depth by 
\( r= \max\{|r|: \omega_{n}\cap I_{r}\neq \emptyset \}. \)

\subsection*{\textbullet\ Chopping times} In all
remaining cases, i.e. if \( \omega_{n}\cap\Delta^{+}\neq \emptyset \)
and \( \omega_{n} \) intersects at least three elements of \( \mathcal
I^{+} \), we say that \( n \) is a chopping time for \( \omega\in 
\mathcal
{P}^{(n-1)} \).  We define a natural subdivision \begin{equation*} 
\omega =
\omega^{\ell}\cup\bigcup_{(r, m)}\omega^{(r, m)}\cup \omega^{\rho}.
\end{equation*} so that each \( \omega_{n}^{(r, m)} \) fully contains
a unique element of \( \mathcal I_{+} \) (though possibly extending to
intersect adjacent elements) and \( \omega_{n}^{\ell} \) and \(
\omega_{n}^{\rho} \) are components of \( \omega_{n}\setminus
\Delta^{+} \) with \( |\omega_{n}^{\ell}|\geq
\delta^{\iota}/(\log\delta^{-\iota})^{2} \) and 
\( |\omega_{n}^{\rho}| \geq \delta^{\iota}/(\log\delta^{-\iota})^{2} \).
If the connected components of \( \omega_{n}\setminus \Delta^{+} \) fail to satisfy the above condition on their length we
just glue them to the adjacent interval of the form \( w_{n}^{(r, m)}
\).  By definition we let each of the resulting subintervals of \(
\omega \) be elements of \( \hat{\mathcal {P}}^{(n)} \). 
The intervals \(\omega^{\ell}, \omega^{\rho} \) and \( 
\omega^{(r, m)} \) 
with \( |r|<r_{\delta} \) are called \emph{escape components} and are 
said to
have an \emph{substantial escape} and \emph{essential escape} 
respectively
at time \( n \). The corresponding values of \( |r|<r_{\delta} \) are 
the associated \emph{essential escape depths}. All other intervals are 
said to have an \emph{essential return} at time \( n \) and the
corresponding values of \( |r| \) are the associated \emph{essential
return depths.}
We remark that partition
elements \( I_{\pm r_{\delta}} \) do not  belong to \( \Delta \) but
we still say that the associated intervals  \( \omega^{(\pm
r_{\delta}, m)} \) have a return rather than an
escape. 

\medskip
This completes the definition of the partition \(\hat{\mathcal {P}}^{(n)}\) of \(
\Omega^{(n-1)} \) and of  the function \(\mathcal {E}^{(n)} \) on \( 
\Omega^{(n-1)} \). We 
define 
\begin{equation}\label{exc}
\Omega^{(n)}= \{a\in\Omega^{(n-1)}: \mathcal {E}^{(n)}
(a) \leq {\alpha}n \}. 
\end{equation}
Notice  that \( \mathcal E^{(n)} \) is constant on elements of
\( \hat{\mathcal {P}}^{(n)}\). Thus \( \Omega^{(n)} \) is the union of 
elements of \( \hat{\mathcal P}^{(n)} \) and we can define
\[
\mathcal {P}^{(n)}=\hat{\mathcal P}^{(n)}|_{\Omega^{(n)}}.
\]
We shall prove the following
\begin{proposition}\label{space}
Conditions \( (BR)_{n}, (EG)_{n}, (SR)_n, (BD)_{n} \)
all hold for \( \Omega^{(n)} \).    
\end{proposition}
Notice that condition \( (BR)_{n} \) is satisfied by construction,
this is precisely the criteria by which the elements of \( \Omega^{(n)} \)
are chosen, see \eqref{exc}. Thus the key part of the proposition are 
the remaining three conditions. This is the general step of an
induction, it allows us to iterate the construction indefinitely and
thus define the set 
\[ 
\Omega^{*}=\bigcap_{n\geq 0}\Omega^{(n)}.  
\]
Notice that in particular  for every \( a\in\Omega^{*} \) the map \( f_{a} \) has 
an exponentially growing derivative along the critical orbit and thus 
exhibits stochastic dynamics. 
Our main Theorem therefore reduces to the following 

\begin{proposition}\label{parameter}
    \( |\Omega^{*}|\geq (1-\eta)|\Omega| \).
 \end{proposition}   
We shall prove Proposition~\ref{space} in Section~\ref{s:posexp} and Proposition~\ref{parameter} 
 in Section~\ref{s:posmeas}. First of all however, 
 in Section~\ref{binding}, we obtain several estimates 
 which are used extensively in the proofs of both
Propositions. First we make a few remarks concerning the 
conditions introduced in this section, all of which 
have played important roles
in various results in one-dimensional dynamics. 

The exponential growth
condition was first introduced by Collet and Eckmann \cite{ColEck83}
who showed that it implies the existence of an absolutely continuous
invariant measure for unimodal maps. Indeed, this is in some sense
the most important condition for us as it guarantees that parameters
in \( \Omega^{*} \) are stochastic. For the purposes of the induction 
however this is not sufficient. Indeed, as often is the case with
inductive arguments, it is easier to assume stronger conditions even if
that means that we have to prove stronger estimates. The slow
recurrence condition in the slightly different form \( e^{-\sqrt n} \)
was used by Benedicks and Carleson \cite{BenCar85}, and in the form
given here, in conjunction with the exponential growth condition in
\cite{BenCar91}.  The bounded recurrence condition, in the precise way
in which it is stated here, is new, but it is closely related to the
\emph{free period} assumption of \cite{BenCar91}.

\section{The binding period}
 \label{binding}

In this section we make precise the definition of the \emph{binding 
period} which is part of the combinatorial information given above, and
obtain several analytic estimates which play an important role in the 
following. In accordance with our inductive assumptions, we 
assume throughout this section that
 that the sets  \( \Omega^{(k)} \) and \( \mathcal P^{(k)} \)
are defined 
and the conditions \( (BR)_{k}, (SR)_{k},
(EG)_{k} \) and \( (BD)_{k} \) hold for all \( k\leq n-1 \). 

Now let \( k\leq n-1 \) and suppose that 
\( \omega\in\mathcal P^{(k)}\) has an 
essential 
or inessential return or an essential or inessential 
escape at time \( k \) with return depth \( r \). 
For each \( a\in\omega \) we let 
\begin{equation}\label{bindcond}
p(c_{k}(a)) = \min\{i: |c_{k+1+i}(a) - c_{i}(a)| \geq e^{-
\alpha_{1} i}\}.
\end{equation}
This is the time for which the future orbit of \( c_{k}(a) \) can be 
thought of as \emph{shadowing} or \emph{being bound to} the orbit of the 
critical point (that is, in some sense, the number of iterations for 
which the orbit of \( c(a) \) repeats its early history after the \( k 
\)'th iterate). 
Then we define 
the \emph{binding period} of \( \omega_{k} \) as 
\[ 
p(\omega_{k}) = \min_{a\in \omega}\{p(c_{k}(a))\}.
\]
The main result of this section is the following
\begin{lemma} \label{uniformbind} 
    For  every \( a\in\omega \)  we have 
\begin{equation}\label{up}
p \leq \tau_{0} \log |c_{k}(a)|^{-1}\leq 
\min\{\tau_{0}(r+1),  
\tau_{0}\alpha_{0} k\}
    < k
\end{equation}
and
\begin{equation}\label{unibindder}
    |(f^{p+1})'(c_{k}(a))| \geq   \mathcal D_{2}\mathcal D_{3}
    \Gamma_{1}e^{r (1-\gamma_{1})} 
    r^{2} 
    \geq
  \mathcal D_{2}\mathcal
    D_{3}   \Gamma_{1}r^{2} e^{- (1+\lambda_{0})}
    e^{\frac{1-\gamma_{1}}{\tau_{0}}(p+1)}.
 \end{equation}   
If \( k \) is an essential return  or an essential escape,  then 
\begin{equation}\label{bindexp2}
|\omega_{k+p+1}|\geq  \Gamma_{1}e^{-\gamma_{1}r}.
\end{equation}
\end{lemma}
  
We shall prove this result in a sequence of sublemmas. We 
assume the notation and set up of the Lemma throughout.

\subsection{Binding period estimates for individual parameter values}
\label{individual}
We start by obtaining several estimates related to the ``pointwise''
binding period \( p(c_{k}(a)) \). To simplify 
the notation we write \( x=c_{k}(a) \) and omit the dependence on the 
parameter \( a \) where there is no risk of confusion.

\begin{sublemma}[Distortion]\label{dist1}
For all \( a\in \omega \), all \( y_{0}, z_{0}\in [x_{0}, c_{0}] \) 
and all \( 0\leq
j\leq \min\{p, k\} \) we have
\[
\frac{|(f_a^{j})'(z_{0})|}{|(f_a^{j})'(y_{0})|} \leq  
  \mathcal
D_{1}.
\] 
\end{sublemma}

We remark that the distortion bound is formally calculated for 
iterates \( j\leq \min\{p, k\} \) because we need to make use of the 
inductive assumption \( (SR)_{k} \) in the proof.  However, we shall show
in the next lemma that we always have \( p< k \) and therefore the estimates do
indeed hold throughout the duration of the binding period. 

\begin{proof}
By the Mean Value Theorem and the standard inequality \( \log(1+x) < x \)
for \( x>0 \), we have
\begin{align*}
\left|\log\frac{|(f_a^{j})'(z_{0})|}{|(f_a^{i})'(y_{0})|}\right| &=\left|\log
\prod_{i=0}^{j-1} \frac{|f_a'(z_{i})|}{|f_a'(y_{i})|}\right| \\
&\leq\sum_{i=0}^{j-1}
\log\left(1+\frac{|f'(z_{i})-f'(y_{i})|}{|f'(y_{i})|}\right)
\\
&
\leq \sum_{i=0}^{j-1}
\frac{|f'(z_{i})-f'(y_{i})|}{|f'(y_{i})|}
\\
&\leq\frac{M_2}{L_2}\sum_{i=0}^{j-1}\frac{|z_i-y_i|}{|y_i|}.
\end{align*}
By the definition of the binding period we have  
\(|z_i-y_i|\leq e^{-\alpha_1i}\); by the definition of \( N_{1} \) in 
\eqref{N1} we have \( |y_{i}|\geq 1 \) for all \( i\leq N_{1} \); and 
by   \((SR)_k\) and the definition of the binding period in
\eqref{bindcond} we have 
\(|y_i|\geq e^{-\alpha_0i}-e^{-\alpha_1 i}\) for \( i> N_{1} \). 
Therefore we have
\begin{align*}
    \sum_{i=0}^{j-1}\frac{|z_i-y_i|}{|y_i|}
&\leq 
\sum_{i=0}^{N_{1}} {e^{-\alpha_{1} i}}+
\sum_{i=N_{1}+1}^{\infty}
\frac{e^{-\alpha_{1} i}}{e^{-\alpha_0i}-e^{-\alpha_{1}i}}\\
&\leq \frac{e^{-\alpha_1}}{1-e^{-\alpha_1}}+\frac{1}{1-e^{-(\alpha_{1} 
- \alpha_{0})(N_1+1)}}\sum_{i=N_{1}+1}^{\infty}
e^{-(\alpha_{1}-\alpha_0) i}\\
&=\frac{e^{-\alpha_1}}{1-e^{-\alpha_1}}+
\frac{e^{-(\alpha_1-\alpha_0)(N_1+1)}}{(1-e^{-(\alpha_{1} 
- \alpha_{0})(N_1+1)})(1-e^{-(\alpha_1-\alpha_0)})}.
\end{align*}
The result then follows by the definition of \( \mathcal D_{1} \) in
\eqref{D1}. 
\end{proof}

\begin{sublemma}[Duration]\label{binding1}
For all \( a\in\omega \) and \( p=p(c_{k}(a)) \) we have 
\begin{equation}\label{disteq1} 
    p \leq 
 \frac{2\log
    |c_{k}(a)|^{-1}+\log (\mathcal D_{1} L_1^{-1}) +
    \lambda_{0}+\alpha_{1}}
    {\lambda_{0}+\alpha_{1}} \leq 
    \tau_{0}\log
    |c_{k}(a)|^{-1}. 
\end{equation}
 In particular 
 \begin{equation}\label{p}
     p \leq \min\{\tau_{0}(r+1),  
\tau_{0}\alpha_{0} k\}
    < k.
    \end{equation}
\end{sublemma}

\begin{proof} 
    Let \( \hat p = \min\{p, k\} \). We shall show that the above 
    estimates work for \( \hat p \) and obtain as a corollary that 
    \( \hat p < k \) and therefore \( p=\hat p \). 
    For simplicity let \( \gamma_{0}=[x_{0}, c_{0}] \) and \( 
    \gamma_{j}=f^{j}(\gamma_{0}) \). 
    Then the Mean Value Theorem and Sublemma~\ref{dist1} imply
    \( |\gamma_{\hat p-1}|\geq \mathcal D^{-1}_{1}|(f^{\hat p-1})'
(c_{0}) | |\gamma_{0}|  \), condition \( (EG)_{k} \) gives \(
|(f_{a}^{j})'(c_{0})|\geq  e^{\lambda_{0} j} \) for all \( j\leq
\hat p \), the definition of binding gives 
\( |\gamma_{\hat p-1}|\leq e^{-\alpha_{1}
(\hat p-1)} \), and \eqref{quadratic0} 
gives \( |\gamma_{0}|\geq
L_1|x|^{2} \). 
Combining all these statements gives 
\[ 
e^{-\alpha_{1} (\hat p-1)} \geq 
|\gamma_{\hat p-1}|\geq \mathcal D^{-1}_{1}|(f^{\hat p-1})'
(c_{0}) | |\gamma_{0}|  
\geq 
\mathcal D^{-1}_{1}  e^{\lambda_{0} (\hat p-1)} |\gamma_{0}|  
\geq 
\mathcal D^{-1}_{1}  L_1 e^{\lambda_{0} (\hat p-1)} |x|^{2}.
\]
Rearranging gives 
\(
e^{(\lambda_{0}+\alpha_{1}) (\hat p-1)} \leq \mathcal
D_1L_{1}^{-1}|x|^{-2}, 
\) and taking logarithms on both sides we get 
\begin{equation}
\hat p-1\leq \frac{2\log
|x|^{-1}+\log (\mathcal D_{1} L_1^{-1})}
{\lambda_{0}+\alpha_{1}}.
\end{equation}
Now, since \( |x|\leq \delta^{\iota}  \) we can use the definition of \( 
k_{0} \) in \eqref{k0} to get 
\begin{equation}\label{pbound} 
    \begin{aligned}
\hat p &\leq 
\frac{2\log
|x|^{-1}+\log (\mathcal D_{1} L_1^{-1}) + \lambda_{0}+\alpha_{1}}
{\lambda_{0}+\alpha_{1}}
\\ &\leq \frac{1}{\lambda_{0}+\alpha_{1}}\left(2\log |x|^{-1} + 
\frac{\log (\mathcal D_{1}L_{1}^{-1}) +
\lambda_{0}+\alpha_{1}}{\log\delta^{-\iota}}\log\delta^{-\iota}\right)
\\& 
\leq 
\frac{1}{\lambda_{0}+\alpha_{1}}\left(2 + 
\frac{\log (\mathcal D_{1}L_{1}^{-1}) +
\lambda_{0}+\alpha_{1}}{\log\delta^{-\iota}}\right) \log |x|^{-1}
\\ &\leq
\frac{2+k_{0}}{\lambda_{0}+\alpha_{1}} \log |x|^{-1} = 
\tau_{0}\log |x|^{-1}.
\end{aligned}
\end{equation}
By condition \( 
(SR)_{k} \) 
we also have \( |c_{k}(a)|\geq e^{-\alpha_{0}k} \) 
we then get 
\begin{equation}\label{up1} 
\hat p\leq \tau_{0}\alpha_{0} k,  
\end{equation}
which is \( < k \) by \eqref{C2}. In particular, \( p=\hat p < k \). 
Moreover, from \eqref{pbound} and the fact that \( |c_{k}(a)|\geq 
e^{-r}/2 \) we also get
\(
p\leq \tau_{0} (r+\log 2).
\)
\end{proof}

\begin{sublemma}[Expansion]\label{bindder} 
For all \( a\in\omega \) and \( p=p(c_{k}(a)) \) we have
\begin{equation}\label{binddereq}
|(f^{p+1})'(c_{k}(a))| \geq 
C_{3} |c_{k}(a)|^{\alpha_{1}\tau_{1}-1}.
\end{equation}
\end{sublemma}
\noindent Recall that \( \alpha_{1}\tau_{1}:=
2\alpha_{1}/(\lambda_{0}+\alpha_{1}) \), see definition in  \eqref{k0},
is strictly less than 1 because \( \alpha_{1}<\lambda_{0} \) by the
choice of \( \alpha_{1} \) in \eqref{alpha1}. 

\begin{proof} 
    We estimate first of all the average expansion between
\( \gamma_{0} \) and \( \gamma_{p} \). Since \( |\gamma_{0}|\leq
L_{1}^{-1}|x|^{2} \) by \eqref{quadratic0} 
and \( |\gamma_{p}| \geq e^{-\alpha_{1} p} \)
by the definition of \( p \), we have \( |\gamma_{p}|/|\gamma_{0}|
\geq L_{1}e^{-\alpha_{1} p} |x|^{-2}.  \) Applying the Mean Value
Theorem we conclude that there exist some point \(
\xi_{0}\in\gamma_{0} \) for which \( |(f^{p})'(\xi_{0})|\geq
L_{1}e^{-\alpha_{1} p} |x|^{-2} \). Then, using the bounded distortion
Sublemma~\ref{dist1} it follows that 
\( 
|(f^{p})'(x_{0})|\geq \mathcal
D^{-1}_{1} L_{1} e^{-\alpha_{1} p} |x|^{-2} 
\). 
Finally, using that fact
that \( |f'(x)|\geq L_{1}|x| \) we obtain
\begin{equation}\label{finalder} 
    |(f^{p+1})'(x)| = |f'(x)|\
|(f^{p})'(x_{0})| \geq 
\mathcal D^{-1}_{1} L_{1}^{2}e^{-\alpha_{1} p}
|x|^{-1}.  
\end{equation} 
Using the upper bound on \( p \) from Sublemma
\ref{binding1} gives
\[ 
e^{-\alpha_{1} p} 
\geq e^{-\alpha_{1} 
\frac{2\log |x|^{-1}
+ \log (\mathcal D_1L_{1}^{-1})}{\lambda_0+\alpha_{1}}}
=(\mathcal D_1L_{1}^{-1})^{\frac{-\alpha_{1}}{\lambda_0+\alpha_{1}}}
|x|^{\frac{2\alpha_{1}}{\lambda_0+\alpha_{1}}}. 
\]
Substituting this into \eqref{finalder} gives 
the result, recall definition of \( C_{3} \) in \eqref{c3}. 
\end{proof}

\subsection{Parameter dependence}
We are now almost ready to proof our main Lemma~\ref{uniformbind}. 
This involves some estimates regarding the dependence of images of
critical points on the parameter.  We shall therefore need the following 
statement which is of intrinsic interest and of wider scope and will be used 
again later on in the argument.  

\begin{lemma}\label{pardist} 
 For any \( 1\leq k \leq  n-1, \omega\in\mathcal P^{(k)} \) and \( a\in\omega \)
 we have
 \[ 
\mathcal D_{2}\geq \frac{|c_{k+1}'(a)|}{|(f^{k+1}_{a})'(c_{0})|} \geq
\mathcal D_{3}^{-1}.
\]
In particular, for all \( 1\leq i < j\leq k+1 \), there exists \(
\tilde a\in\omega \) such that
\begin{equation}\label{pardisteq}  
\frac{1}{\mathcal D_{2}\mathcal D_{3}}
|(f_{\tilde a}^{j-i})'
(c_{i}(\tilde a)| \leq 
\frac{|\omega_{j}|}{|\omega_{i}|}\leq 
\mathcal D_{2}\mathcal D_{3} |(f_{\tilde a}^{j-i})'((c_{i}(\tilde a))|.
\end{equation}
\end{lemma}
\begin{proof} 
Recall that 
\(
c_{k+1}(a) = f(c_{k}(a)) + a.
\)
Therefore, by the chain rule we have
\begin{equation*}
\begin{aligned}
c_{k+1}'(a) &= 
 1 + f'(c_{k}) c'_{k}(a)
\\ &= 1 + f'(c_{k}) [1+f'(c_{k-1}) c'_{k-1}(a)]
\\&= 
1 + f'(c_{k}) + f'(c_{k}) f'(c_{k-1}) c'_{k-1}(a)
\\&=
1 + f'(c_{k}) + f'(c_{k}) f'(c_{k-1})
[1+f'(c_{k-2})c'_{k-2}(a)]
\\&=
1 + f'(c_{k}) + f'(c_{k}) f'(c_{k-1})
+ f'(c_{k}) f'(c_{k-1}) f'(c_{k-2})c'_{k-2}(a)
\\&=\dots
\\&=
1+
f'(c_{k})+f'(c_{k})f'(c_{k-1})+\ldots 
+f'(c_{k})f'(c_{k-1}) \dots f'(c_{1})f'(c_{0}).
\end{aligned}%
\end{equation*} 
Dividing both sides by 
\(
(f^{k+1})'(c_{0})= f'(c_{k})f'(c_{k-1}) \dots
f'(c_{1})f'(c_{0}) \) gives 
\begin{equation}\label{parder}
\frac{c_{k+1}'(a)}{(f^{k+1}_{a})'(c_{0})} = 
1+ \sum_{i=1}^{k+1}\frac{1}{(f^{i})'(c_{0})}.
\end{equation}
Now let \( \tilde N \) be as in condition
(A4) and in the definition of the constants 
\( \mathcal D_{2}, \mathcal
D_{3} \) in \eqref{D2}-\eqref{F3}. Then, for \( 1\leq k < \tilde N \) we have 
\[
\left|\frac{c_{k+1}'(a)}{(f^{k+1}_{a})'(c_{0})}\right| = 
\left|1+ \sum_{i=1}^{k+1}\frac{1}{(f^{i})'(c_{0})}\right| \leq
\mathcal D_{2}.
\]
Using the fact that by
\( (EG)_{k} \) we have 
\( |(f^{i})'(c_{0})| \geq e^{\lambda_{0}i}
\) for all \( 1\leq i\leq k+1 \), 
for any \(  k  \geq \tilde N \)  we have 
\begin{align*}
\left|\frac{c_{k+1}'(a)}{(f^{k+1}_{a})'(c_{0})} \right| 
&=
\left|1+ \sum_{i=1}^{\tilde N}\frac{1}{(f^{i})'(c_{0})}
+\sum_{i=\tilde N+1}^{k+1}\frac{1}{(f^{i})'(c_{0})}\right|
\\
&\leq 
\left|1+ \sum_{i=1}^{\tilde N}\frac{1}{(f^{i})'(c_{0})}\right|
+\left|\sum_{i=\tilde N+1}^{k+1}\frac{1}{(f^{i})'(c_{0})}\right|
\\ 
&\leq 
\left|1+ \sum_{i=1}^{\tilde N}\frac{1}{(f^{i})'(c_{0})}\right| + 
\sum_{i=\tilde N+1}^{k+1}e^{-\lambda_{0}i}
\\ 
&\leq 
\left|1+ \sum_{i=1}^{\tilde N}\frac{1}{(f^{i})'(c_{0})}\right| + 
\sum_{i=\tilde N+1}^{\infty}e^{-\lambda_{0}i}
\\ 
&\leq 
\left|1+ \sum_{i=1}^{\tilde N}\frac{1}{(f^{i})'(c_{0})}\right| + 
\frac{e^{-\lambda_{0}(\tilde N+1)}}{1-e^{-\lambda_{0}}} \leq \mathcal
D_{2}.
\end{align*}
The lower bound is obtained similarly.  
If \( k <  \tilde N \) we have 
\[ 
\left|\frac{c_{k+1}'(a)}{(f^{k+1}_{a})'(c_{0})}\right|  
= 
\left|1+\sum_{i=1}^{k+1}\frac{1}{(f^{i})'(c_{0})}\right| \geq \mathcal
D_{3}^{-1}, 
\]
and if \( k\geq \tilde N  \) we have 
\begin{align*}
    \left|\frac{c_{k+1}'(a)}{(f^{k+1}_{a})'(c_{0})}\right| 
    &\geq
    \left|1+ \sum_{i=1}^{\tilde N}\frac{1}{(f^{i})'(c_{0})}
    +\sum_{i=\tilde N+1}^{k+1}\frac{1}{(f^{i})'(c_{0})}\right|
 \\& \geq 
 \left|1+ \sum_{i=1}^{\tilde N}\frac{1}{(f^{i})'(c_{0})}\right|
     - \left|\sum_{i=\tilde N+1}^{k+1}\frac{1}{(f^{i})'(c_{0})}\right|
\\&\geq
\left|1+\sum_{i=1}^{\tilde N}\frac{1}{(f^{i})'(c_{0})}\right|-
\sum_{i=\tilde N+1}^{\infty} e^{-\lambda_{0}i} \geq \mathcal
D_{3}^{-1}.
\end{align*}
This completes the proof of the first set of
inequalities. For the second, 
consider the composition \(
c_{j} \circ c_{i}^{-1} :\omega_{i}\to\omega_{j} \). 
By the Mean Value Theorem there exists
\( \tilde a \in\omega \) such that
\[
\frac{|\omega_{j}|}{|\omega_{i}|}
=| (c_{j} \circ c_{i}^{-1})'((c_{i}(\tilde a))|.
\]
Then by the chain
rule and the first set of inequalities we get the desired statement.
\end{proof}

Lemma~\ref{pardist} implies that all parameters in  a
given interval \( \omega \) satisfy comparable derivative estimates 
during the binding period \( p=p(\omega_{k}) \). This allows us to
extend the expansion estimates at the end of a binding period to
parameter for which we do not necessarily have \( p(c_{k}(a)) =
p(\omega_{k}) \), recall that the
binding period of the entire interval is defined as the minimum of
these binding periods. 

\begin{sublemma}
    Let  \( p=p(\omega_{k}) \). Then, for all \( a\in \omega \) we have 
    \begin{equation}\label{bexp} 
   |(f^{p+1})'(c_{k}(a))| \geq \tilde C_{3} |c_{k}(a)|^{\alpha_{1}\tau_{1}-1}.
    \end{equation}
\end{sublemma}
\begin{proof}
For those  parameter values \( \tilde a\in\omega \) such that \(
p(c_{k}(\tilde a))= p(\omega_{k}) \), the result follows immediately from 
Sublemma~\ref{bindder}. Notice that there must exist such a \( \tilde
a\). 
 For a generic \( a \) we argue as 
follows. First of all we consider the iterates coming  after the return.
For any \(
a, \tilde a\in\omega \), we have
\begin{equation}\label{par1} 
   \begin{aligned}
       |(f^{p}_{a})'(c_{k+1}(a))| & \geq \frac{1}{\mathcal D_{1}}
    |(f^{p}_{a})'(c_{0}(a))|  \quad \text{by Sublemma~\ref{dist1}},\\
    &\geq \frac{1}{\mathcal D_{1}\mathcal D_{2}}|c'_{p}(a)|
    \quad \text{by Lemma~\ref{pardist}},
    \\ &\geq \frac{1}{\mathcal D\mathcal D_{1}\mathcal D_{2}}
    |c'_{p}(\tilde a)| \quad \text{by the inductive assumption }
    (BD)_{k},
\\    &\geq \frac{1}{\mathcal D\mathcal D_{1}\mathcal D_{2}\mathcal D_{3}}
    |(f^{p}_{\tilde a})'(c_{0}(\tilde a))|   \quad \text{by
    Lemma~\ref{pardist}},
    \\
    &\geq \frac{1}{\mathcal D\mathcal
    D_{1}^{2}\mathcal D_{2}\mathcal D_{3}} 
    |(f^{p}_{\tilde a})'(c_{k+1}(\tilde a))| \quad \text{by Sublemma~\ref{dist1}}.
\end{aligned}
\end{equation} 
Now we deal with the actual return iterate. 
By construction we have \(
|\omega_{k}|\leq |c_{k}(\tilde a)|/2.  \) 
Therefore, using
\eqref{quadratic0}, we get
\begin{equation}\label{par2} 
 \begin{aligned} 
|f'(c_{k}(a))| &\geq
L_{1}|c_{k}(a)| \geq L_{1}(\vert c_{k}(\tilde a)\vert -
|\omega_{k}|) \geq L_{1} |c_{k}(\tilde a)|/2 \geq L_{1}^{2}
|f'(c_{k}(\tilde a))|/2.  
 \end{aligned} 
\end{equation} 
Since \(
|(f^{p+1})'(c_{k}(a))| = |f'(c_{k}(a))|\ |(f^{p})'(c_{k+1}(a))| \) and
similarly for \( \tilde a \), \eqref{par1} and \eqref{par2} imply 
\begin{equation}\label{par3}
|(f^{p+1})'(c_{k}(a))| \geq
\frac{ L_{1}^{2}}{2 \mathcal D\mathcal D_{1}^{2}\mathcal D_{2}
\mathcal D_{3}} 
    |(f^{p+1})'(c_{k}(\tilde a))| 
\end{equation}
for all \( a, \tilde a \in\omega  \).  Now, choosing \(  \tilde a \) such that 
\( p(\omega_{k}) = p(c_{k}(\tilde a)) \) and applying \eqref{binddereq} 
and the fact that 
  \( |c_{k}(a)|\leq 2 |c_{k}(\tilde a)| \), we get 
  \begin{equation} \label{bexp1}
      |(f^{p+1})'(c_{k}( \tilde a))| 
       \geq  C_{3}|c_{k}(\tilde a)|^{\alpha_{1}\tau_{1}-1} 
   \geq  C_{3} 2^{\alpha_{1}\tau_{1}-1}
      |c_{k}(a)|^{\alpha_{1}\tau_{1}-1}.
\end{equation}
Substituting \eqref{bexp1} into \eqref{par3} and using the definition 
of \( \tilde C_{3} \) in \eqref{c3} gives the result. 
  \end{proof}
  
\begin{proof}[Proof of Lemma~\ref{uniformbind}]
The upper bound \eqref{up} on \( p \) follows 
immediately from \eqref{disteq1}. For the expansion estimates,
\eqref{bexp} and the fact that \( |c_{k}(a)|\leq e^{-(r-1)} \),
imply
\[ 
|(f^{p+1})'(c_{k}( a))| \geq 
\tilde C_{3} |c_{k}(a)|^{\alpha_{1}\tau_{1}-1}
\geq
\tilde C_{3} e^{-(r-1)(\alpha_{1}\tau_{1}-1)}.
\]
To get the first inequality in \eqref{unibindder} it is therefore
sufficient to show that 
\begin{equation}\label{iff1}
\tilde C_{3} e^{-(r-1)(\alpha_{1}\tau_{1}-1)} \geq 
r^{2}\mathcal D_{2}\mathcal D_{3} 
\Gamma_{1}e^{(1-\gamma_{1})r}.
\end{equation}
Rearranging and taking logs we get that \eqref{iff1} holds if and only
if
 \begin{equation}\label{beexp4}
 \gamma_{1}\geq  \frac{\log (\Gamma_{1}\mathcal D_{2}\mathcal 
 D_{3}\tilde C_{3}^{-1} e^{\alpha_{1}\tau_{1}-1})+2\log 
 r}{r}+\alpha_{1}\tau_{1}.
 \end{equation}
 Since the expression on the right hand side is strictly decreasing 
 with with \( r \), it is sufficient for this inequality to be 
 verified for \( r=\iota \log\delta^{-1} \) which follows from
 \eqref{C3}. 
This gives the first inequality in \eqref{unibindder}. 
To get the second inequality in \eqref{unibindder} we use the fact
that \( r> p\tau_{0}^{-1}-1 \) which comes from 
\( p\leq \tau_{0}(r+1) \) in  \eqref{p}, and get 
\begin{align*}
|(f^{p+1})'(c_{k}( a))| &\geq \Gamma_{1}\mathcal D_{2}\mathcal D_{3} 
e^{(1-\gamma_{1})r} r^{2}
\\ &\geq \Gamma_{1}\mathcal D_{2}\mathcal
D_{3} e^{p \frac{1-\gamma_{1}}{\tau_{0}} }e^{-(1-\gamma_{1})}r^{2}
\\& = 
\Gamma_{1}\mathcal D_{2}\mathcal
D_{3} e^{(p+1)
 \frac{1-\gamma_{1}}{\tau_{0}} }
e^{ -\frac{1-\gamma_{1}}{\tau_{0}} }e^{-(1-\gamma_{1})}r^{2}
\\& = \Gamma_{1}\mathcal D_{2}\mathcal
D_{3} r^{2} e^{-\frac{(\tau_{0}+1)(1-\gamma_{1})}{\tau_{0}}} e^{(p+1)
 \frac{1-\gamma_{1}}{\tau_{0}} }.
\end{align*}
We then use the fact that \( \tau_{0}> 2/(\lambda_{0}+\alpha_{1}) \)
and \( \alpha_{1}<\lambda_{0} \) by definition, to get 
\[ 
\frac{(\tau_{0}+1)(1-\gamma_{1})}{\tau_{0}} \leq 
\frac{\tau_{0}+1}{\tau_{0}} = 1+\frac{1}{\tau_{0}} <
1+\frac{\lambda_{0}+\alpha_{1}}{2} < 1+\lambda_{0}. 
\]
This completes the proof of \eqref{unibindder}. 
Now suppose 
 that \( k \) is an essential return or escape, in particular \( 
 |\omega_{k}|\geq e^{-r}/r^{2} \).  Therefore 
 Lemma~\ref{pardist} implies that  there exists some parameter \( 
 a\in\omega \) such that 
 \[ 
 |\omega_{k+p+1}| \geq \frac{1}{\mathcal D_{2}\mathcal D_{3}} 
 |(f^{p+1})'(c_{k}( a))|  \ |\omega_{k}| 
 \geq \frac{e^{-r}}{r^{2}\mathcal D_{2}\mathcal D_{3}} 
 |(f^{p+1})'(c_{k}( a))|,
 \]
and the result follows immediately by \eqref{unibindder}. 
\end{proof}
 
 \subsection{Uniform bounds on the sum of all inessential return
 depths}
 
 The following result is essentially a Corollary of Lemma~\ref{uniformbind}. It shows that the total length of the binding
 periods  associated to a sequence of inessential returns can be
 bounded in terms of the return depth of the immediately preceding
 \emph{essential} return. 

 \begin{lemma}\label{inessential}
     Let \( \omega\in\mathcal P^{(n)} \) and suppose that \( \omega \) 
     has an essential return with return depth \( r_{0} \) 
     at some time \( \mu_{0}\leq n-1 \). Let 
     \( \mu_{1}, \ldots, \mu_{u} \) denote the following inessential 
     returns which occur after time \( \mu_{0} \) and before any 
     subsequent chopping time and before time \( n \). Let \( p_{i} \), \( i=0,\ldots, u \) 
     denote the corresponding binding periods. Then 
     \[ 
     \sum_{i=0}^{u}
 (    p_{i}+1) \leq \tau r_{0}.
     \] 
 \end{lemma}

 \begin{proof}
 Since \( \mu_{1}, \dots, \mu_{u} \) are all returns (to \( \Delta \)),
 condition (A1)  gives  
 \[
 |(f_{a}^{\mu_{i + 1} - (\mu_i + p_i + 1)})'(c_{\mu_i +
 p_i + 1}(a))| \geq 1
 \]
 for \(  0 \leq i\leq u-1\) and 
 \(
 \forall\ a\in\omega \). Notice that we do not have (and will not need)
 such an estimate for \( i=u \) since we do not have information about 
 the location of \( \omega_{n} \) (nor indeed do we even know if \(
 n\geq \mu_{u}+p_{u}+1 \)).  Equation \eqref{unibindder} in Lemma~\ref{binding} 
 together with the fact that \( r_{i}^{2}\geq
 r_{\delta}^{2}= (\log \delta)^{2} \geq e^{-(1+\lambda_{0})} \) by condition (C1), 
 gives 
 \[ 
 |(f^{\mu_{i}+p_{i}+1})'(c_{\mu_{i}(a)})| \geq \mathcal D_{2}\mathcal
 D_{3}\Gamma_{1} 
 e^{\frac{1-\gamma_{1}}{\tau_{0}}(p+1)}
 \]
 for all  \( 0 \leq i\leq u\) and  \( a\in\omega \). Notice that this
 estimate does include the case \( i=u \). 
 Combining these estimates by the chain rule, and applying Lemma~\ref{pardist} we then get 
 \[
 \frac{\vert\omega_{\mu_{i+1}}\vert}{\vert\omega_{\mu_{i}}\vert}
 \geq
 \Gamma_{1}e^{\frac{1-\gamma_{1}}{\tau_{0}}(p_{i}+1)}
 \quad \text{for }  i=0,\ldots, u-1, 
 \ \text{ and }\ 
 \frac{\vert\omega_{\mu_u+p_u+1}\vert}{\vert\omega_{\mu_0}\vert}
 \geq \Gamma_{1} e^{\frac{1-\gamma_{1}}{\tau_{0}}(p_{u}+1)},
 \]
 and  therefore
 \[
 \frac{\vert\omega_{\mu_u+p_u+1}\vert}{\vert\omega_{\mu_0}\vert}=
 \frac{\vert\omega_{\mu_u+p_u+1}\vert}{\vert\omega_{\mu_u}\vert}
 \frac{\vert\omega_{\mu_u}\vert}{\vert\omega_{\mu_{u-1}}\vert}\cdots
 \frac{\vert\omega_{\mu_2}\vert}{\vert\omega_{\mu_1}\vert}
 \frac{\vert\omega_{\mu_1}\vert}{\vert\omega_{\mu_0}\vert}
 \geq 
 \Gamma^{u+1} e^{\frac{1-\gamma_{1}}{\tau_{0}} \sum_{i=0}^{u} (p_{i}+1)}.
 \]
 Now since \( \mu_{0} \) is an essential return we have \(
 |\omega_{\mu_{0}}|\geq e^{-r_{0}}/r_{0}^{2} \) and thus
 \[ 
 |I|
 \geq 
 |\omega_{\mu_{u}+p_{u}+1}| 
 \geq |\omega_{\mu_{0}}|
 \Gamma^{u+1} e^{\frac{1-\gamma_{1}}{\tau_{0}} \sum_{i=0}^{u} (p_{i}+1)}
 \geq 
 \Gamma^{1} e^{\frac{1-\gamma_{1}}{\tau_{0}} \sum_{i=0}^{u}
 (p_{i}+1)}e^{-r_{0}}/r_{0}^{2}. 
 \]
 We have replaced \( \Gamma^{u+1}_{1} \) by \( \Gamma_{1} \) in the
 inequality because we do not know how many inessential returns there
 are. There may be none, in which case we have \( u=0 \). 
 Taking logs and rearranging we get 
 \begin{equation}\label{eq:tau}
     \begin{aligned}
 \sum_{i=0}^{u} (p_{i}+1) &\leq 
 \frac{\tau_{0}}{1-\gamma_{1}} 
 \log \left(\frac{|I| e^{r_{0}}r_{0}^{2}}{\Gamma_{1}}\right) 
 \\ &= r_{0} \frac{\tau_{0}}{1-\gamma_{1}} \left(1+\frac{\log |I| + 2\log
 r_{0} - \log \Gamma_{1}}{r_{0}}\right) 
 \\
 &\leq 
 r_{0} \frac{\tau_{0}}{1-\gamma_{1}} \left(1+\frac{\log |I| + 2\log
 r_{\delta} - \log \Gamma_{1}}{r_{\delta}}\right) 
 \\ &=
 \tau r_{0}.
 \end{aligned}
 \end{equation}
 The last inequality follows from the fact that the fraction in
 parenthesis is decreasing in \( r \), therefore \( r_{0} \) can be
 replaced by \( r_{\delta} \) giving precisely the definition of \(
 \tau \) in \eqref{tau}.  Notice moreover that 
 \begin{equation}\label{eq:p}
 \sum_{i=0}^{u}(p_{i}+1) \geq p_{0}+1 > 0,
 \end{equation}
 and therefore the above inequality implies \emph{a fortiori} that \(
 \tau > 0 \). 
 \end{proof}

\section{Positive exponents in dynamical space} 
\label{s:posexp}

In this section we prove Proposition~\ref{space}. 
Notice that the combinatorics and the recurrence condition \( (BR)_{n}
\) are satisfied for every  $a\in\Omega^{(n)}$
by construction. We therefore just need to verify the slow recurrence, 
exponential growth, and bounded distortion conditions. We shall do
this in three separate subsections.

\subsection{Slow recurrence}

\begin{lemma} \label{srlemma}
    For every \( a\in\Omega^{(n)} \), \( (SR)_{n} \) holds.
\end{lemma}
\begin{proof}
    The statement clearly holds if \( n \) is not an escape or a 
    return or a bound iterate for \( \omega \) containing \( a \). 
    Now suppose that \( n \) is an escape iterate. Then it must be 
    that \( n \geq N \) and then \eqref{def1} implies that 
    \( e^{-\alpha_0 n} \leq 
    e^{-\alpha_0 N} < \delta \). Since \( n \) is an escape we must have 
    \( |c_{k}(a)| \geq \delta \geq e^{-\alpha_0 N} \geq e^{-\alpha_0 n} \).
    If \( n \) is an essential return the result follows immediately 
    by the bounded recurrence condition \( (BR)_{k} \).  
    if \( n \) is an inessential return, it follows immediately from 
    the binding period estimates Lemma~\ref{inessential}
    that its return depth must be less 
    than the return depth of the preceding inessential return depth 
    and thus in particular satisfy the required estimate. If \( n \) 
    belongs to a binding period the same reasoning gives the result.

\end{proof}

\subsection{Exponential derivative growth}

\begin{lemma}\label{posexp}
For every \( a\in \Omega^{(n)} \), \( (EG)_{n} \) holds. 
\end{lemma}
\begin{proof}
   Let \( \omega\in\mathcal P^{(n)} \) 
be the element containing \( a \). 
If \( \omega \) has no returns before time \( n \) 
then this implies that \( c_{k}(a) \notin\Delta \) for all \( 
k\leq n \) and therefore 
we have \( |(f^{n+1})'(c_{0}(a))|\geq e^{\lambda (n+1)} >   
e^{\lambda_{0} (n+1)}\) 
by condition (A1) and the fact that \( 
\lambda>\lambda_{0} \). 
If  \( \omega \) has a non-empty sequence of returns before time \( n
\), let 
\(\nu_1<\nu_2<\cdots<\nu_q\leq n\) be all the free (essential and
inessential) returns up to \(n\)
and \(p_i\) the corresponding binding period.
Lemma~\ref{inessential} 
and the bounded recurrence condition \( (BR)_{\nu_{q}} \) 
imply
\begin{equation}\label{pn}
\sum_{i=1}^q (p_{i}+1) \leq \tau  \!\!\!\! \!\!\!\!  
\sum_{\text{essential returns}} \!\!\!\! \!\!\!\!   r_{i}(\omega) 
=\tau \mathcal 
E^{(\nu_q)}(\omega)  \leq 
\tau\alpha \nu_q.
\end{equation}
Splitting the orbit of \( c_{0} \) into free and 
bound iterates and using condition (A1),  the binding period 
estimate \eqref{unibindder},  and (\ref{pn}),  we get  
\begin{equation}\label{exp1}
\begin{aligned}
|(f^{\nu_q})'(c_{0})|   
&\geq 
e^{\lambda
(\nu_q-\sum_{i=1}^q(p_{i}+1))}\Gamma_{1}^{q}
e^{\frac{1-\gamma_{1}}{\tau_{0}}\sum_{i=1}^{q}(p_{i}+1)}\\
&\geq \Gamma_{1}^{q} e^{\lambda
\nu_q-(\lambda-
\frac{1-\gamma_{1}}{\tau_{0}})\sum_{i=1}^{q}(p_{i}+1)}\\
&\geq \Gamma_{1}^{q} e^{\lambda
\nu_q-\tau\alpha\nu_{q}(\lambda-
\frac{1-\gamma_{1}}{\tau_{0}})}\\
& \geq \Gamma_{1}^{q} e^{(\lambda-\tau\alpha(\lambda-
\frac{1-\gamma_{1}}{\tau_{0}}))\nu_{q}}.
\end{aligned}
\end{equation}
Then condition \( (BR)_{\nu_{q}} \) implies in particular \(
|c_{k}|\geq e^{-(r_{q}+1)} \) with \( r_{q}\leq \alpha\nu_{q} \)
which, together with \eqref{quadratic0}, implies 
\begin{equation}\label{exp1a}
|f'(c_{\nu_{q}})|\geq L_{1} e^{-\alpha \nu_{q}-1}.
\end{equation}
If \( n=\nu_{q} \), equations \eqref{exp1} and \eqref{exp1a} give 
\begin{equation}\label{nu=n}
|(f^{n+1})'(c_{0})|   
\geq L_{1}  \Gamma_{1}^{q} e^{-1} 
e^{-(\lambda-\tau\alpha(\lambda-
\frac{1-\gamma_{1}}{\tau_{0}})-\alpha )}
e^{(\lambda-\tau\alpha(\lambda-
\frac{1-\gamma_{1}}{\tau_{0}})-\alpha )(n+1)},
\end{equation}
notice that 
\begin{equation}\label{alphadef}
    \lambda-\tau\alpha(\lambda-
\frac{1-\gamma_{1}}{\tau_{0}})-\alpha \geq \lambda_{0}
\Longleftrightarrow 
\alpha\leq
(\lambda-\lambda_{0})
(\tau(\lambda-\frac{1-\gamma_{1}}{\tau_{0}})+1).
\end{equation}
This condition is therefore satisfied by the definition of \( \alpha \)
in \eqref{alpha}.  In particular, from \eqref{nu=n} we have 
\[ 
|(f^{n+1})'(c_{0})|   
\geq L_{1}  \Gamma_{1}^{q} e^{-1-\lambda_{0}} 
e^{\lambda_{0}(n+1)}\geq e^{\lambda_{0}(n+1)},
\]
since \( \Gamma_{1}^{q} \geq \Gamma_{1}\geq L_{1}^{-1}e^{1+\lambda} \)
by \eqref{Gamma1}. This proves the statement in the case \( n=\nu_{q} \).

If \( n> \nu_{q} \), 
it remains to estimate the derivative along the remaining iterates. We
claim that 
\begin{equation}\label{exp1b}
|f^{(n-\nu_{q})})'(c_{\nu_{q}+1})| \geq \mathcal D_{1}^{-1}
e^{\lambda_{0} (n-\nu_{q})}.
\end{equation}
Indeed, if \( n\leq\nu_q+p_q \), i.e. \( n \) belongs to the binding
period following the return at time \( \nu_{q} \), then this follows
immediately from the inductive assumption \( (EG)_{n-\nu_{q}-1} \) and
the bounded distortion during binding periods in Lemma~\ref{dist1}. 
Otherwise, i.e. if \( n> \nu_q+p_q  \), we consider two possibilities.
If the orbit of \( c_{k} \) never enters \( \Delta \) between time \( 
\nu_{q}+1 \) and time \( n \) we just ignore the binding period and
use (A1) to obtain \(  |f^{(n-\nu_{q}-1)})'(c_{\nu_{q}+1})| 
\geq 
e^{\lambda (n-\nu_{q}-1)}\) which clearly implies \eqref{exp1b}. 
Alternatively, if the orbit of \( c_{k} \) does enter \( \Delta \)
between time \(  \nu_{q}+1 \) and time \( n \) we just apply the
inductive assumption and the bounded distortion to get an estimate
like \eqref{exp1b} up to the last time before \( n \) that \(
c_{k}\in \Delta \) and then apply (A1) to get exponential growth at a 
rate \( \lambda \) for the remaining iterates; also in this case  
we obtain \eqref{exp1b}.

By the chain rule and equations \eqref{exp1}-\eqref{exp1b} we then get
\begin{align*}
|(f^{n+1})'(c_{0})|   &\geq \Gamma_{1}^{q} e^{(\lambda-\tau\alpha(\lambda-
\frac{1-\gamma_{1}}{\tau_{0}}))\nu_{q}} 
L_{1} e^{-\alpha \nu_{q}-1} \mathcal D_{1}^{-1}
e^{\lambda_{0} (n-\nu_{q})}
\\ & = \Gamma_{1}^{q}L_{1}\mathcal D_{1}^{-1} e^{-1-\lambda_{0}} 
e^{(\lambda-\tau\alpha(\lambda-
\frac{1-\gamma_{1}}{\tau_{0}})-\alpha-\lambda_{0})\nu_{q}} 
e^{\lambda_{0}(n+1)} \geq e^{\lambda_{0}(n+1)}.
\end{align*}
The last inequality follows, again, by \eqref{alphadef} and
\eqref{Gamma1}.  This completes the proof of Lemma~\ref{posexp}. 
\end{proof}

\subsection{Bounded distortion}

\begin{lemma}
\label{distortion}
Let \( \omega\in\mathcal{P}^{(n)} \). 
Then 
\[
\frac{|c_{k}'a|}{|c_{k}'b|}\leq
\mathcal{D} 
\] 
for all \( a,b
\in\omega \) and all \( k\leq \nu_{q}+p_{q}+1 \) where \( \nu_{q} \leq n\)
is the last essential or inessential return, or the last essential or 
inessential escape of \( \omega \) and \( p_{q} \) is the corresponding
binding period.  If \( n>\nu_{q}+p_{q}+1 \) 
then the same statement
holds for all \( k\leq n \) restricted to any subinterval \(
\bar\omega \subseteq\omega \) such that \(
|\bar\omega_{k}|\subseteq\Delta^{+} \). 
\end{lemma}

By Lemma~\ref{pardist} it is sufficient to show that
\begin{equation}\label{itinerary}
\frac{|(f_{a}^{k})'(c_{0})|}{|(f_{b}^{k})'(c_{0})|}\leq \frac{\mathcal
D}{\mathcal{D}_2\mathcal{D}_3}.
\end{equation}
Notice that, strictly speaking, Lemma~\ref{pardist} is stated for \(
k\leq n-1 \) but this is not an issue here since the actual assumptions
used and therefore the conclusions of the lemma 
clearly apply up to time \( k \).  
Equation~\eqref{itinerary} essentially says that critical orbits with the same
combinatorics satisfy comparable derivative estimates. 
By standard
arguments (see also the proof of Sublemma~\ref{dist1}),
we have 
\[ 
\log\frac{|(f_{a}^{k})'(c_{0})|}{|(f_{b}^{k})'(c_{0})|}
\leq  \frac{M_{2}}{L_{2}} \sum_{j=0}^{k-1}
D_{j}
\quad \text{where } \quad D_{j}=
\frac{|\omega_{j}|}{d(\omega_{j})}
\]
and   \( d(\omega_{j}) = \inf_{a\in\omega}|c_{j}(a)| \).
Let \( 0<\nu_{1}<\dots<\nu_{q}\leq n \) be  
the sequence of essential and inessential returns and essential and inessential
escapes of \( \omega \).  By construction there is a unique
element \( I_{\rho_{i}, m_{i}} \) in \( \mathcal I^{+} \) associated
to each \( \nu_{i} \).  Let \( p_{i} \) be the length of the binding
period associated to \( \nu_{i} \).  For notational convenience define
\( \nu_{0} \) and \( p_{0} \) so that \( \nu_{0}+p_{0}+1=0 \).  We
suppose first that \( k\leq \nu_{q}+p_{q}+1 \). Then write 
\begin{equation}\label{distsum}
    \sum_{j=0}^{k-1}D_{j} \leq 
    \sum_{j=0}^{\nu_{q}+p_{q}}D_{j}=
\sum_{i=0}^{q-1}
\sum_{\nu_{i}+p_{i}+1}^{\nu_{i+1}+p_{i+1}}D_{j}.
\end{equation}
Notice that we have divided the itinerary of \( \omega \) into a
finite number of blocks corresponding to pieces of itinerary starting 
immediately after a binding period and going through to the end of the 
following binding period.   
In the next sublemma we obtain a bound for the sum over each
individual block. 

\begin{sublemma}\label{lemma1} 
    For each \( i=0, \dots, q-1 \) we
have 
\begin{equation}\label{distsum1}
    \sum_{\nu_{i}+p_{i}+1}^{\nu_{i+1}+p_{i+1}}D_{j} \leq
\hat D \ | \omega_{\nu_{i+1}}|
e^{\rho_{i+1}}. 
\end{equation}
\end{sublemma}
\begin{proof}
 We start first of all by subdividing the sum further into iterates
 corresponding to: \( i \))  the free iterates between the end of the binding
 period and the following return, \( ii \)) the return, \( iii \)) the binding period
 following the return. 
\begin{equation}\label{sum}
\sum_{\nu_{i}+p_{i}+1}^{\nu_{i+1}+p_{i+1}}D_{j}=
\sum_{\nu_{i}+p_{i}+1}^{\nu_{i+1}-1}D_{j} + D_{\nu_{i+1}} +
\sum_{\nu_{i+1}+1}^{\nu_{i+1}+p_{i+1}}D_{j}.  \end{equation}
We shall
estimate each of the three terms in separate
arguments.  For the first,  notice that since \(
\omega_{\nu_{i+1}}\subseteq\hat I_{\rho_{i}, m}\subset \Delta^{+} \), 
condition (A1)
implies \( |(f^{\nu_{i+1}-j})'c_{j}(a)|\geq C_{1}e^{\lambda
(\nu_{i+1}-j)} \). Therefore by  Lemma~\ref{pardist} we have \(
|\omega_{j}|\leq C_{1}^{-1}\mathcal{D}_2\mathcal{D}_3
e^{-\lambda (\nu_{i+1}-j)}|\omega_{\nu_{i+1}}|
\). Moreover \( d(\omega_{j})\geq \delta^{\iota}-2I_{r_{\delta^+}}
/r_{\delta^+}^2\geq\delta^{\iota}/2 \), therefore, using also the fact
that \( \delta^{\iota}\geq e^{-\rho_{i+1}} \),
 we have
\begin{equation}\label{lemma1eq1}
\sum_{j=\nu_{i}+p_{i}+1}^{\nu_{i+1}-1}D_{j}
\leq 2C_{1}^{-1}\mathcal{D}_2\mathcal{D}_3
\sum_{j=\nu_{i}+p_{i}+1}^{\nu_{i+1}-1} e^{- \lambda(\nu_{i+1}-j)}|
\omega_{\nu_{i+1}}| \delta^{-\iota} \leq \frac{2C_{1}^{-1}\mathcal{D}_2\mathcal{D}_3
e^{-\lambda}}{1-e^{-\lambda}}
|\omega_{\nu_{i+1}}|e^{\rho_{i+1}}.  
\end{equation} 
For the second
term we immediately have \begin{equation}\label{lemma1eq2}
D_{\nu_{i+1}}\leq 2 |\omega_{\nu_{i+1}}| e^{ \rho_{i+1}}.
\end{equation}
since \( D_{\nu_{i+1}} \) is the supremum of \(
|\omega_{\nu_{i+1}}|/|c_{\nu_{i+1}}| \) over all points \(
c_{\nu_{i+1}} \) in \( \omega_{\nu_{i+1}} \) and \(
|c_{\nu_{i+1}}|\geq \frac{1}{2} e^{-\rho_{i+1}}\) by definition.  

The estimation for the third term is the trickiest, or at least the
least intuitive.  
By Lemma~\ref{pardist} we have, for all \( j\in
[\nu_{i+1}+1, \nu_{i+1}+p_{i+1}] \), 
\begin{equation}\label{deriv}
|\omega_{j}|\leq \mathcal{D}_2\mathcal{D}_3|\omega_{\nu_{i+1}}| 
\sup_{a\in\omega}
\{|(f_{a}^{j-\nu_{i+1}})'(c_{\nu_{i+1}})(a)|\}.  
\end{equation}
We therefore need an \emph{upper bound} for 
\[ 
|(f_{a}^{j-\nu_{i+1}})'(c_{\nu_{i+1}})(a)| = 
|(f^{j-\nu_{i+1}-1})'(c_{\nu_{i+1}+1}(a))| \cdot 
|f'_{a}(c_{\nu_{i+1}}(a))|.
\]
For the second term in the product we just have 
\(
|f'_{a}(c_{\nu_{i+1}}(a))|\leq 2L_1^{-1}
e^{-\rho_{i+1}} 
\) by \eqref{quadratic0}. Thus, substituting into \eqref{deriv} gives
\begin{equation}\label{deriv1a}
|\omega_{j}|\leq \mathcal{D}_2\mathcal{D}_3 2L_1^{-1}
\sup_{a\in\omega}
\{|(f^{j-\nu_{i+1}-1})'(c_{\nu_{i+1}+1}(a))| \} |\omega_{\nu_{i+1}}| 
e^{-\rho_{i+1}}.  
\end{equation}
Thus it only remains to get an upper bound for the derivative during
the binding period. 
Fix \( a\in\omega \) and let  
$\gamma_{j-\nu_{i+1}-1}:=|c_{j-\nu_{i+1}-1}(a)-c_j(a)|$. 
Then we have \( |\gamma_{0}|\geq L_1
e^{-2\rho_{i+1}} \) and, by the definition of binding periods, 
\(
|\gamma_{j-\nu_{i+1}-1}|\leq e^{-\alpha_1 (j-\nu_{i+1}-1)} 
\). 
Therefore using 
the Mean Value Theorem  and Lemma~\ref{dist1} which says that all
derivatives are comparable during the binding period, we get
\[
|(f^{j-\nu_{i+1}-1})'(c_{\nu_{i+1}+1}(a))|\leq \mathcal{D}_1
|\gamma_{j-\nu_{i+1}-1}|/|\gamma_{0}|\leq \mathcal{D}_1 
e^{-\alpha_1 (j-\nu_{i+1}-1)}L_1^{-1}
e^{2\rho_{i+1}}.  
\] 
Substituting this into \eqref{deriv1a} then gives 
\begin{equation}\label{deriv1b}
|\omega_{j}|\leq 2 \mathcal{D}_1 \mathcal{D}_2\mathcal{D}_3 L_1^{-2} 
e^{-\alpha_1 (j-\nu_{i+1}-1)}|\omega_{\nu_{i+1}}| e^{\rho_{i+1}}.  
\end{equation}
To bound \( d(\omega_{j}) \) we just observe that the definition of
binding period and the 
slow recurrence condition $(SR)_{j-\nu_{i+1}-1}$ imply, for \( j>
\nu_{i+1}+1 \), 
\begin{equation}\label{deriv1c}
    \begin{aligned}
|c_j(a)| &\geq |c_{j-\nu_{i+1}-1}(a)| - e^{-\alpha_{1} (j-\nu_{i+1}-1)} 
\\ &\geq e^{-\alpha_{0}(j-\nu_{i+1}-1)} - e^{-\alpha_{1} (j-\nu_{i+1}-1)}
= e^{-\alpha_0(j-\nu_{i+1}-1)}(1-e^{-(\alpha_1-\alpha_0)(j-\nu_{i+1}-1)})
\\ &\geq e^{-\alpha_0(j-\nu_{i+1}-1)}(1-e^{-(\alpha_1-\alpha_0)}).
\end{aligned}
\end{equation}
Notice that for \( j=\nu_{i+1}+1 \) we just have \( |c_j(a)|  \geq 1
\geq e^{-\alpha_0(j-\nu_{i+1}-1)}(1-e^{-(\alpha_1-\alpha_0)}) \) so,
formally, the inequality holds in this case also. 
Thus, \eqref{deriv1b} and \eqref{deriv1c} and the fact that \(
\alpha_{1} \geq \alpha_{0} \) give 
\begin{equation}\label{deriv1d}
\frac{|\omega_{j}|}{d(\omega_{j})}\leq 
\frac{2 \mathcal{D}_1 \mathcal{D}_2\mathcal{D}_3 L_1^{-2} 
e^{-\alpha_1 (j-\nu_{i+1}-1)}|\omega_{\nu_{i+1}}| e^{\rho_{i+1}}}
{e^{-\alpha_0(j-\nu_{i+1}-1)}(1-e^{-(\alpha_1-\alpha_0)})}
\leq \frac{2 \mathcal{D}_1 \mathcal{D}_2\mathcal{D}_3 L_1^{-2} 
|\omega_{\nu_{i+1}}| e^{\rho_{i+1}}}
{1-e^{-(\alpha_1-\alpha_0)}}.
\end{equation}
Substituting \eqref{lemma1eq1}, \eqref{lemma1eq2} and \eqref{deriv1d} 
into \eqref{sum} and using the definition of \( \hat D \) in
\eqref{dhat} gives the statement in Sublemma \eqref{lemma1}. 
\end{proof}

Substituting \eqref{distsum1} in \eqref{distsum} we get 
\begin{equation}\label{sum1} 
    \sum_{j=0}^{\nu_{q}+p_{q}} D_{j}=
\sum_{i=0}^{q-1}\sum _{\nu_{i}+p_{i}+1}^{\nu_{i+1}+p_{i+1}}D_{j}\leq
\hat D \sum_{i=1}^{q}|\omega_{\nu_{i}}|e^{\rho_{i}}.
\end{equation} 
Now we subdivide the sum on the right hand side into partial sums
corresponding to return times with the same return depth \( r \)
denoted by \( i: \rho_{i}=r \):
\begin{equation}\label{sum2} 
    \sum_{i=1}^{q}|\omega_{\nu_{i}}|e^{
\rho_{i}} =\sum_{r\geq r_{\delta^{+}}}e^{ r}\sum_{i: \rho_{i}=r}
|\omega_{\nu_{i}}|.  
\end{equation} 
In the next sublemma we estimate the total contribution of of returns 
corresponding to a fixed return depth \( r \). 

\begin{sublemma} \label{sumr} For
any \( r\geq r_{\delta^{+}} \), 
\[
\sum_{i:\rho_{i}=r}|\omega_{\nu_{i}}|\leq \left(2+
e\frac{(\log\delta^{-\iota})^{2}}{(\log\delta^{-\iota}-1)^{2}}\right) 
\frac{e^{- r}}{r^{2}}
\sum_{j=0}^{\infty}
\biggl(\frac{C_{1}^{-1}\delta^{\iota (1-\gamma_1)}}{(\log
\delta^{\iota})^2}\biggr)^j.  
\]
\end{sublemma}

\begin{proof} Let \( \mu_{j}=\nu_{i_{j}}, j=1, \dots, m
\) be the subsequence of returns and escapes
with return depths equal to \( r \). 
Using Lemma~\ref{uniformbind} and condition (A1) we have for all \( a\in\omega
\) and \( j=1, \dots, m-1 \), \[
|(f^{\mu_{i+1}-\mu_{i}}_{a})'(c_{\mu_{i}}(a))|\geq 
C_{1}\mathcal D_2\mathcal D_3e^{(1-\gamma_1) r} r^2\geq
C_{1}\mathcal D_2\mathcal D_3e^{(1-\gamma_{1}) r_{\delta^{+}}}r_{\delta^+}^2.\]
Therefore by Lemma~\ref{pardist}, 
\( |\omega_{\mu_{i}}|/|\omega_{\mu_{i+1}}|\leq
C_{1}^{-1}\delta^{\iota (1-\gamma_1)}/
(\log\delta^{\iota})^2 \) and 
\[
\sum_{j=1}^{m}|\omega_{\mu_{i}}|\leq \sum_{j=0}^{m-1}
\biggl(\frac{C_{1}^{-1}\delta^{\iota (1-\gamma_1)}}{(\log
\delta^{\iota})^2}\biggr)^j |\omega_{\mu_{m}}| \leq
 \sum_{j=0}^{\infty}\biggl(\frac{C_{1}^{-1}\delta^{\iota (1-\gamma_1)}}{(\log
\delta^{\iota})^2}\biggr)^j|\omega_{\mu_{m}}|.  
\]
Recall that $\omega_m$ possibly spreads across three contiguous
elements of $\mathcal{I}^+$. Two of these have length at most \(
e^{-r}/r^{2} \) and the third one has length at most \(
e^{-(r-1)}/(r-1)^{2} \leq (er^{2}/(r-1)^{2})(e^{-r}/r^{2}) \leq  
(er_{\delta^{+}}^{2}/(r_{\delta^{+}}-1)^{2})(e^{-r}/r^{2}) \).
This gives 
\[
|\omega_{\mu_{m}}|\leq \left(2+
e\frac{(\log\delta^{-\iota})^{2}}{(\log\delta^{-\iota}-1)^{2}}\right) 
\frac{e^{-r}}{r^2}. 
\]
\end{proof} 

Notice that the convergence of the sum on the right hand side is
guaranteed by the definition of  \( \gamma_{1} \) in \eqref{min}. 
Indeed some
straightforward rearrangements shows that condition \eqref{min}
implies that
\( {C_{1}^{-1}\delta^{\iota (1-\gamma_1)}}/{(\log
\delta^{\iota})^2} <1\).

Substituting
the estimate of Sublemma~\ref{sumr} into 
\eqref{sum2} now gives 
\[ 
\sum_{i=1}^{q}|\omega_{\nu_{i}}|e^{
\rho_{i}} =\sum_{r\geq r_{\delta^{+}}} e^{ r}
\sum_{i: \rho_{i}=r}
|\omega_{\nu_{i}}| \leq 
\sum_{r\geq r_{\delta^{+}}} 
\left[\left(2+
e\frac{(\log\delta^{-\iota})^{2}}{(\log\delta^{-\iota}-1)^{2}}\right) 
\sum_{j=0}^{\infty}
\biggl(\frac{C_{1}^{-1}\delta^{\iota (1-\gamma_1)}}{(\log
\delta^{\iota})^2}\biggr)^j\right] \frac{1}{r^{2}}.
\]
Summing over \( r \) and substituting into \eqref{sum1} then gives 
\begin{align*}
\sum_{j=0}^{\nu_{q}+p_{q}} D_{j} &
\leq \hat D 
\left[\left(2+
e\frac{(\log\delta^{-\iota})^{2}}{(\log\delta^{-\iota}-1)^{2}}\right) 
\sum_{j=0}^{\infty}
\biggl(\frac{C_{1}^{-1}\delta^{\iota (1-\gamma_1)}}{(\log
\delta^{\iota})^2}\biggr)^j\right] \frac{1}{r_{\delta^{+}}-1} 
\\ &=
\hat D 
\left[\left(2+e\frac{(\log\delta^{-\iota})^{2}}
{(\log\delta^{-\iota}-1)^{2}}\right) 
\frac{(\log\delta^{\iota})^2} {(\log\delta^{\iota})^2 - 
C_{1}^{-1}\delta^{\iota (1-\gamma_1) }}
\right] \frac{1}{\log\delta^{-\iota}-1} = \hat D \hat{\hat D}.
\end{align*}
This completes the proof of  Lemma~\ref{posexp} for \( k \leq \nu_{q}+p_{q}+1
\).  
If \( k>\nu_{q}+p_{q}+1 \) we consider the additional terms \(
D_{j} \) restricting ourselves to some subinterval \(
\bar\omega\subset \omega \) with \( \bar\omega_{k}\subseteq\Delta^{+}
\). Clearly the preceding estimates are unaffected by this
restriction.  Using (A1) and Lemma~\ref{pardist} we have 
\[ 
|\bar\omega_{j}|\leq 
C_{1}^{-1}\mathcal{D}_2\mathcal{D}_3e^{-\lambda (k-j)}|\bar\omega_{k}|\leq
C_{1}^{-1}\mathcal{D}_2\mathcal{D}_3  e^{-\lambda
(k-j)}\delta^{\iota}, 
\] 
and therefore ,using also the fact that \(
|c_{j}(a)|\geq \delta^{\iota} \) since \(
\omega_{j}\cap\Delta^{+}=\emptyset \), we get 
\[
\sum_{j=\nu_{q}+p_{q}+1}^{k-1}D_{j}\leq C_{1}^{-1}\mathcal{D}_2\mathcal{D}_3
\sum_{j=\nu_{q}+p_{q}+1}^{k-1}
e^{-\lambda (k-j)} \leq C_{1}^{-1}\mathcal{D}_2\mathcal{D}_3
\sum_{i=1}^{\infty} e^{-\lambda i}
= \frac{C_{1}^{-1}\mathcal{D}_2\mathcal{D}_3
e^{-\lambda}}{1-e^{-\lambda}}.  
\]

\section{Positive measure in parameter space}
\label{s:posmeas}

In this section we prove Proposition~\ref{parameter}.  We divide the
proof into 4 sections.

\subsection{Large deviations}\label{deviations}
Recall first the
definition of \( \Omega^{(n)} \) in \eqref{exc}:
\[ 
\Omega^{(n)}=\{a\in\Omega^{(n-1)}: \mathcal E^{(n)}(a) \leq \alpha n\}.
\]
We will show here that the \emph{average value} of \( \mathcal E^{(n)} \)
on \( \Omega^{(n-1)} \) is \emph{low} and therefore \emph{most}
parameter in \( \Omega^{(n-1)} \) make it into \( \Omega^{(n)} \).
More precisely, let \( \gamma, \gamma_{2} \) be as in Section 
\ref{constcond}, then we have the following 
\begin{proposition} For every \( n\geq 1 \), 
\label{average}
    \begin{equation*}
    \int_{\Omega^{(n-1)}} \hspace{-.8cm} e^{\gamma_{2}\cal{E}{n}} 
    \leq \left(1+\sum_{R\geq r_{\delta}}e^{(\gamma -1) 
    R}\right)^{n}
    |\Omega|.
    \end{equation*}
\end{proposition}

We will prove this proposition in the next three sections. First we
show how, by a simple large deviation argument, it implies Proposition
\ref{parameter} and therefore our Main Theorem. 
The definition of \( \Omega^{(n)} \) gives 
\[ 
|\Omega^{(n-1)}|-|\Omega^{(n)}|=|\Omega^{(n-1)}\setminus\Omega^{(n)}|=
|\{\omega\in\hcal Pn=\cal Qn: e^{\gamma_{2}\cal En}\geq 
e^{\gamma_{2}\alpha n}\}|. 
\]
Therefore, using Chebychev's inequality (large deviations) 
and Proposition~\ref{average} we have
\[
|\Omega^{(n-1)}|-|\Omega^{(n)}| 
\leq 
e^{-\gamma_{2} \alpha n}\int_{\Omega^{(n-1)}}\hspace{-.5cm}e^{\gamma_{2}\cal En}
\leq 
\left[e^{-\gamma_{2} \alpha}\left(1+  \sum_{R\geq r_{\delta}} 
 e^{(\gamma -1) R}\right)\right]^{n}
|\Omega|,
\]
which implies, recall the definition of \( \tilde\eta \) in Section
\ref{constcond}, 
\begin{equation*} 
   |\Omega^{(n)}|\geq |\Omega^{(n-1)}|-\left[e^{-\gamma_{2} \alpha}
   \left(1+  \sum_{R\geq r_{\delta}} 
 e^{(\gamma-1) R}\right)\right]^{n} |\Omega| = 
 |\Omega^{(n-1)}| - \tilde\eta |\Omega|.
\end{equation*}
Iterating this expression and using the fact that no exclusions are
made before time N, and using also the definition of \( \eta \), we have  
\[ 
|\Omega^{*}|\geq \left(1-\sum_{j=N}^{\infty}
\left[e^{-\gamma_{2} \alpha}\left(1+  \sum_{R\geq r_{\delta}} 
 e^{(\gamma -1) R}\right)\right]^{j} \right) |\Omega| =
 (1-\eta)|\Omega|.
\]

\subsection{Renormalization properties of the combinatorics}
\label{altcomb}

The proof of Proposition~\ref{average}is quite subtle. 
In this section 
we give an alternative combinatorial description of the parameters in \( 
\Omega^{(n-1)} \) and state a key technical Proposition~\ref{main} in 
terms of this combinatorial description. 
This description is crucial to the 
argument and 
highlights some remarkable renormalization properties of the 
construction. 
We shall also show how
Proposition~\ref{main} implies Proposition~\ref{average}. 

Recall that \hcal Pn is the partition of \( 
\Omega^{(n-1)} \) which takes into account the dynamics at time \( n \) 
and which restricts to the partition \( \mathcal P^{(n)} \) of \( 
\Omega^{(n)} \) after the exclusion of a certain elements of \( \hcal 
Pn \).
To each \( \omega\in\hcal Pn \) is 
associated a sequence 
\( 0=\eta_{0}<\eta_{1}<\dots <\eta_{s}\leq n, \ s=s(\omega)\geq 0 
\) of escape times and 
a corresponding sequence of escaping components 
\(
    \omega\subseteq\omega^{(\eta_{s})}\subseteq\dots
    \subseteq\omega^{(\eta_{0})}
\)  \quad 
with 
\(\omega^{(\eta_{i})}\subseteq\Omega^{(\eta_{i})} \)
and 
\( \omega^{(\eta_{i})}\in\mathcal P^{(\eta_{i})}.\)
To simplify the formalism we also define some ``fake'' escapes by 
letting
\(
    \omega^{(\eta_{i})}=\omega
\) 
for all  \(s+1\leq i\leq n\). In this way we have a
well defined parameter 
interval \( \omega^{(\eta_{i})} \) associated to \( \omega\in
\hcal Pn \) for each \( 0\leq i\leq n \). 
Notice that for two intervals \( \omega, \tilde\omega\in \hcal Pn \) and 
any \( 0\leq i\leq n \), the corresponding intervals 
\( \omega^{(\eta_{i})} \) and \( \tilde\omega^{(\eta_{i})} \) 
are either disjoint or coincide.
Then we define
\[ 
Q^{(i)}=\bigcup_{\omega\in\hcal Pn}\omega^{(\eta_{i})}
 \] 
 and let 
 \(
 \cal Qi = \{\omega^{(\eta_{i})}\}
 \)
 denote the natural partition of 
 \( Q^{(i)} \) into intervals of the form \( \omega^{(\eta_{i})} \). 
 Notice that 
 \(
 \Omega^{(n-1)}= Q^{(n)}\subseteq\dots \subseteq
 Q^{(0)}=\Omega^{(0)} 
 \) 
 and 
  \(\cal Qn=\hcal Pn \)
 since
 the number \( s \) of escape times is always strictly 
 less than \( n \) and therefore in particular 
 \( \omega^{(\eta_{n})}=\omega \) 
 for all \( \omega\in\hcal Pn \).
 For a given \( \omega=\omega^{(\eta_{i})}\in \cal Qi, 
 \ 0\leq i \leq n-1 \) we let 
 \[ 
  Q^{(i+1)}(\omega) = \{\omega'=\omega^{(\eta_{i+1})}
 \in \cal Q{i+1}: \omega'\subseteq \omega\} 
 \] 
 denote all the elements of \( \cal Q{i+1} \) which are contained 
 in \( \omega \) 
 and let \( \cal Q{i+1}(\omega) \) denote the corresponding partition.
Then we define a function \( \Delta\cal Ei: Q^{(i+1)}(\omega) \to 
\mathbb N  \) by
\[ 
\Delta\cal Ei (a)=\cal E{\eta_{i+1}}(a)
-\cal E{\eta_{i}}(a). 
\]
 This gives the total sum of all essential return 
depths associated to the itinerary of the element \( \omega'
\in\cal Q{i+1}(\omega) \) containing \( a \), between the escape 
at time \( \eta_{i} \) and the escape at time \( \eta_{i+1} \).  
Clearly \( \Delta\cal Ei (a) \) is constant on elements of 
\( \cal Q{i+1}(\omega) \). 
Finally we let 
\[
 \cal Q{i+1}(\omega, R) = \{\omega'\in \mathcal Q^{(i+1)}: 
 \omega'\subseteq \omega, \Delta\cal Ei (\omega')= R\}.
 \]
 Notice that the entire construction given here depends on \( n \). 
The main motivation for this construction and is the following 

 \begin{proposition}\label{main}
 For all \( i \leq n-1\), \( \omega\in \cal Qi \)  and \( R\geq 0 \) we have 
 \[ 
 \sum_{\tilde\omega\in \cal Q{i+1}(\omega, R)} \!\!\!\!\! |\tilde\omega|
 \leq e^{(\gamma_{1}+\gamma_{0} -1) R} |\omega|.
 \]
 \end{proposition}

 This says that the probability of accumulating a large 
 total return depth between one escape and the next is exponentially 
 small. 
The strategy for proving this result is straightforward. We show that 
the size of each partition element \( \omega'\in \cal Q{i+1}(\omega, R) \) 
is exponentially small in \( R \) and then use a combinatorial 
argument to show that the total number of such elements cannot be too 
large. The proposition follows immediately from Lemmas~\ref{metric}  and \ref{combestlem} in the next two sections. First
however we show how Proposition~\ref{main} implies Proposition
\ref{average}.

\begin{proof}[Proof of Proposition~\ref{average} assuming Proposition~\ref{main}] 
 Notice first of all that    
    \begin{equation}\label{main1}
    \int_{\Omega^{(n-1)}} \hspace{-.8cm} e^{\gamma_{2}\cal{E}{n}} = 
      \sum_{\omega\in\cal Q{n}} \hspace{-.2cm} e^{\gamma_{2} \cal En (\omega)} |\omega|.
    \end{equation}
    Thus it is sufficient to bound the right hand side.
Let  \(  0\leq i \leq n-1\),
and \( \omega\in \cal Qi \) and recall that  
\( \cal En = \Delta\cal E0 + \dots + \Delta\cal E{n-1} \) 
    and \( \Delta\cal Ei \) is constant on elements of \( \cal Q{i} \).
    Then we write
\begin{equation}\label{main2}
\sum_{\omega'\in\cal Q{i+1}(\omega)} \hspace{-.5cm} e^{\gamma_{2}\Delta\cal Ei (\omega')} 
|\omega'|  = \hspace{-.5cm}
\sum_{\omega'\in\cal Q{i+1}(\omega, 0)}\hspace{-.6cm}|\omega'| + 
\sum_{R\geq r_{\delta}} 
e^{\gamma_{2}R} \hspace{-.6cm}\sum_{\omega'\in\cal Q{i+1}(\omega, R)}
\hspace{-.6cm}|\omega'|. 
\end{equation}
For the first term in the sum we just use the trivial bound 
\begin{equation}\label{main3}
\sum_{\omega'\in\cal Q{i+1}(\omega, 0)}\hspace{-.6cm}|\omega'| 
\leq |\omega|. 
\end{equation}
For the second, we use Proposition~\ref{main} and the definition of \( 
\gamma \) to get 
\begin{equation}\label{main4}
\sum_{R\geq r_{\delta}} 
e^{\gamma_{2}R} \hspace{-.6cm}\sum_{\omega'\in\cal Q{i+1}(\omega, R)}
\hspace{-.6cm}|\omega'| \leq 
\sum_{R\geq r_{\delta}} 
 e^{(\gamma_{0}+\gamma_{1}+ \gamma_{2}-1) R}  |\omega| 
 = 
 \sum_{R\geq r_{\delta}} 
  e^{(\gamma-1) R}  |\omega|. 
\end{equation}
Substituting \eqref{main2} and \eqref{main3} into \eqref{main1} we
get
\begin{equation}\label{i}
    \sum_{\omega'\in\cal Q{i+1}(\omega)} \hspace{-.5cm} e^{\gamma_{2}\Delta\cal Ei (\omega')} 
    |\omega'| \leq
\left(1+  \sum_{R\geq r_{\delta}}  
 e^{(\gamma-1) R}\right) |\omega|. 
\end{equation}
Now, to obtain a bound for \eqref{main1} recall that by construction
each \( \omega^{(n)}\in \mathcal Q^{(n)} \) belongs to a nested
sequence of intervals \( \omega^{(n)}\subseteq \omega^{(n-1)}\subseteq
\dots \subseteq \omega^{(0)}=\Omega\) which each \( \omega_{i} \)
belonging to \( \mathcal Q^{(i)}(\omega^{(i-1)} \) for \( n\geq i\geq 1 \). 
Therefore we can write \eqref{main1} as
\[ 
\sum_{\omega^{(1)}\in \cal Q1 (\omega^{(0)})} 
\hspace{-.8cm} e^{\gamma_{2}\Delta\cal E0 (\omega^{(1)})}
\hspace{-.8cm}
\sum_{\omega^{(2)}\in \cal Q{2}(\omega^{(1)})} 
\hspace{-.8cm} e^{\gamma_{2}\Delta\cal E1 (\omega^{(2)})}
\hspace{-.3cm}\dots
\hspace{-1cm} \sum_{\omega^{(n-1)}\in \cal Q{n-1}(\omega^{(n-2)})} 
\hspace{-1.2cm} e^{\gamma_{2}\Delta\cal E{n-1}(\omega^{(n-1)})}
\hspace{-.7cm}
\sum_{\omega=\omega^{(n)}\in \cal Q{n}(\omega^{(n-1)})} 
\hspace{-.7cm}  e^{\gamma_{2}\Delta\cal E{n-1}(\omega^{(n)})}|\omega|.
\]
Notice the \textit{nested} nature of the expression. Applying \eqref{i} 
repeatedly gives
\[ 
\sum_{\omega\in\cal Qn} \hspace{-.2cm} e^{\gamma_{2}\cal En(\omega)} |\omega| 
\leq
\left(1+  \sum_{R\geq r_{\delta}} 
 e^{(\gamma-1) R}\right)^{n}|\Omega|. 
\]
\end{proof}

\subsection{Metric estimates}

\begin{lemma}
\label{metric}
 For all \(  0\leq i \leq n-1\), \( \omega\in \cal Qi \), \( R\geq 0 \) and 
 \(\tilde\omega\in \cal Q{i+1}(\omega, R)\) 
  we have
 \[ 
 |\tilde\omega|\leq e^{(\gamma_{1}-1) R} |\omega|. 
 \]
 \end{lemma}

 \begin{proof}
From the construction of \( \tilde\omega \) 
there is a nested sequence of intervals 
\[
\tilde\omega=\omega^{(\nu_{s+1})} 
\subseteq\omega^{(\nu_{s})}\subseteq
\dots\subseteq\omega^{(\nu_{1})}\subseteq\omega^{(\nu_{0})}
=\omega 
\]
where \( \omega \) has an escape at time \( \nu_{0} \), 
 each \( \omega^{(\nu_{j})}, \ j=1,\ldots, s  \) has an essential return 
at time \( \nu_{j} \) (intuitively 
\( \omega^{(\nu_{j})} \) is created as a consequence of the 
intersection of \( \omega^{(\nu_{j-1})} \) with \( \Delta \) 
at time \( \nu_{j} \)).
Write
\begin{equation}\label{ratios}
\frac{|\tilde\omega|}{|\omega|}=\frac{|\omega^{(\nu_{1})}|}
{|\omega^{(\nu_{0})}|}
\frac{|\omega^{(\nu_{2})}|}
{|\omega^{(\nu_{1})}|}
\dots \frac{|\omega^{(\nu_{s})}|}{|\omega^{(\nu_{s-1})}|} 
\frac{|\tilde\omega|}{|\omega^{(\nu_{s})}|}. 
 \end{equation}
We shall estimate the right hand side of \eqref{ratios} term by term. 
We start by considering the terms 
\( {{\vert}\omega^{(\nu_{j+1})}\vert}/{\vert\omega^{(\nu_j)}\vert} \)
for   \( j=1,\ldots, s-1 \) for which both \(
\nu_{j} \) and \( \nu_{j+1} \) are essential returns. The idea is to
compare the two parameter intervals by comparing their images and
applying the bounded distortion condition. In this case we choose to
compare their images at some intermediate time between the return at
time \( \nu_{j} \) and the return at time \( \nu_{j+1} \), more
specifically at the end of the binding period following the return at 
time \( \nu_{j} \). Thus, using the bounded distortion condition
we have 
$$
\frac{{\vert}\omega^{(\nu_{j+1})}\vert}{\vert\omega^{(\nu_j)}\vert}=
\frac{{\vert}c_{\nu_j+p_j+1}'(a)\vert}{{\vert}c_{\nu_j+p_j+1}'(b)\vert}\cdot
\frac{\vert\omega^{(\nu_{j+1})}_{\nu_j+p_j+1}\vert}{\vert\omega^{(\nu_j)}
_{\nu_j+p_j+1}\vert}
< 
\mathcal D\frac{\vert\omega^{(\nu_{j+1})}_{\nu_j+p_j+1}\vert}
{\vert\omega^{(\nu_j)}_{\nu_j+p_j+1}\vert}.
$$
To bound the numerator we use   Lemma~\ref{uniformbind} 
 to get 
 \[
 \vert\omega^{(\nu_j)}_{\nu_j+p_j+1}\vert \geq 
 \Gamma_{1} e^{-\gamma_{1}r_{j}}.
 \] 
 To bound the denominator notice 
 that there may (or may not) be a sequence of inessential returns 
 between time \( \nu_j+p_j+1 \) and time \( \nu_{j+1} \). In any case,
 the 
 accumulated derivative taken over all free and bound periods is \(
 \geq 1 \). Therefore, using 
 \eqref{pardisteq}, 
 condition (A1) and Lemma~\ref{uniformbind} we get  
\begin{equation*}
|\omega^{(\nu_{j+1})}_{\nu_{j+1}}|\geq C_{1}\mathcal 
 D_{2}^{-1}\mathcal D_{3}^{-1} |\omega^{(\nu_{j+1})}_{\nu_j+p_j+1}|. 
 \end{equation*}
Thus, by the definition of \( \Gamma_{1} \) in \eqref{Gamma1}, 
 \begin{equation}\label{metric0}
 \frac{{\vert}\omega^{(\nu_{j+1})}\vert}{\vert\omega^{(\nu_j)}\vert}
 \leq \mathcal D\frac{\vert\omega^{(\nu_{j+1})}_{\nu_j+p_j+1}\vert}
{\vert\omega^{(\nu_j)}_{\nu_j+p_j+1}\vert}
\leq \frac{\mathcal D\mathcal D_{2}\mathcal D_{3}}{\Gamma_{1}C_{1}} 
 e^{-r_{j+1}+\gamma_{1}r_{j}} \leq  e^{-r_{j+1}+\gamma_{1}r_{j}}.
 \end{equation}
 For the last term in ratios we just use the trivial bound
  \[ 
  |\tilde\omega|\leq |\omega^{(\nu_{s})}|,
  \] 
and so we get 
\begin{equation}\label{metric1}
\frac{|\tilde\omega|}{|\omega|}\leq 
\frac{|\omega^{(\nu_{1})}|}{|\omega^{(\nu_{0})}|} 
e^{\sum_{j=1}^{s-1} -r_{j+1}+\gamma_{1}r_{j}} 
= 
\frac{|\omega^{(\nu_{1})}|}{|\omega^{(\nu_{0})}|}
e^{r_{1}-\gamma_{1}r_{s}+ \sum_{j=1}^{s} (\gamma_{1}-1) r_{j}}.
\end{equation}
Estimates for  \(  {|\omega^{(\nu_{1})}|}/{|\omega^{(\nu_{0})}|}  \).
are different from the other cases, in principle, since \( \nu_{0} \)
is an escape and not a return time. However, 
if \( \nu_{0} \) is an \emph{essential escape} then we can actually apply 
the binding period estimates and, using exactly the same arguments as 
above we get 
\(  {|\omega^{(\nu_{1})}|}/{|\omega^{(\nu_{0})}|} \leq e^{-r_{1} + 
\gamma_{1}r_{0}}\) and substituting this into \eqref{metric1} and 
using the fact that \( r_{0}< r_{s} \) we get the statement in the 
Lemma in this case. 

It remains to consider the case in which \( \nu_{0} 
\) is a \emph{substantial escape}, i.e. \( 
\omega_{\nu_{0}}^{(\nu_{0})}\) lies outside \( \Delta^{+} \) and 
satisfies \( |\omega_{\nu_{0}}^{(\nu_{0})}|\geq 
\delta^{\iota}/(\log\delta^{-\iota})^{2} \). 
We consider two separate cases depending on whether
its image \( \omega^{(\nu_0)}_{\nu_1} \) satisfies 
$\omega^{(\nu_0)}_{\nu_1}\subset\Delta^+$ or not. 
Suppose first that $\omega^{(\nu_0)}_{\nu_1}\subset\Delta^+$. 
Then we can apply the bounded 
distortion as above and for some 
$a\in{\omega}^{(\nu_0)}$ and $b\in\omega^{(\nu_1)}$ 
we get 
\begin{equation}\label{metric2}
\frac{|\omega^{(\nu_1)}|}{|\omega^{(\nu_0)}|}=
\frac{|c_{\nu_1}'(a)|}{|c_{\nu_1}'(b)|}
\cdot\frac{|\omega^{(\nu_1)}_{\nu_1}|}{|\omega^{(\nu_0)}_{\nu_1}|}
<\mathcal D
\frac{|\omega^{(\nu_1)}_{\nu_1}|}{|\omega^{(\nu_0)}_{\nu_1}|}
\leq \frac{\mathcal D e^{-r_{1}}}{|\omega^{(\nu_0)}_{\nu_1}|}.
\end{equation}
To bound the denominator we use once more condition (A1) and
\eqref{pardisteq} to get 
\begin{equation}\label{metric3}
|\omega^{(\nu_0)}_{\nu_1}| \geq \frac{1}{\mathcal D_{2}\mathcal 
D_{3}}\min_{a\in\omega^{(\nu_{0})}}\{|(f^{\nu_{1}-\nu_{0}})'(c_{\nu_{0}}(a))|\}
|\omega^{(\nu_0)}_{\nu_0}| \geq \frac{1}{\mathcal D_{2}\mathcal 
D_{3}}\frac{\delta^{\iota}}{(\log\delta^{-\iota})^{2}}.  
\end{equation}
Substituting \eqref{metric3} into \eqref{metric2} and then \eqref{metric2} 
into  \eqref{metric1}, and using the fact that \( r_{s}\geq r_{\delta} = 
\log \delta^{-1} \), we get 
\begin{equation}\label{metric4}
   \frac{|\tilde\omega|}{|\omega|}\leq 
    {\mathcal D\mathcal D_{2}\mathcal D_{3}}
  \frac{e^{-\gamma_{1}r_{s}}(\log\delta^{-\iota})^{2}}{\delta^{\iota}} 
  e^{(\gamma_{1}-1) R}
    \leq 
    \frac{\mathcal D\mathcal D_{2}\mathcal D_{3}}{C_{1}}
  \frac{\delta^{\gamma_{1}}(\log\delta^{-\iota})^{2}}{\delta^{\iota}} 
  e^{(\gamma_{1}-1) R}.
 \end{equation}   
By straightforward rearrangement and taking logs we have 
\[ 
\frac{\mathcal D\mathcal D_{2}\mathcal D_{3}}{C_{1}}
  \frac{\delta^{\gamma_{1}}(\log\delta^{-\iota})^{2}}{\delta^{\iota}} 
  \leq 1 
  \Longleftrightarrow 
  \gamma_{1} \geq 
  \iota+\frac{
  \log (\mathcal D \mathcal D_{2} \mathcal 
    D_{3} C_{1}^{-1})+ 2\log\log\delta^{-\iota} }{\log\delta^{-1}}.
\]
The inequality on the right is satisfied by \eqref{C3}, and thus we obtain
our result in this case. 
Now, if
$\omega^{(\nu_0)}_{\nu_1}\not\subset\Delta^+$, the bounded distortion 
result applies only to the subinterval 
$\overline{\omega}^{(\nu_0)}\subset \omega^{(\nu_0)}$ 
such that \( \overline{\omega}^{(\nu_0)}_{\nu_{1}}\subset\Delta^{+} \)
and we get 
\[
\frac{|\omega^{(\nu_1)}|}{|\overline\omega^{(\nu_0)}|}=
\frac{|c_{\nu_1}'(a)|}{|c_{\nu_1}'(b)|}
\cdot\frac{|\omega^{(\nu_1)}_{\nu_1}|}{|\overline\omega^{(\nu_0)}_{\nu_1}|}
<\mathcal D
\frac{|\omega^{(\nu_1)}_{\nu_1}|}{|\overline\omega^{(\nu_0)}_{\nu_1}|}
\leq \frac{\mathcal D e^{-r_{1}}}{|\overline\omega^{(\nu_0)}_{\nu_1}|}.
\]
However in this case 
recall that 
\( \omega^{(\nu_0)}_{\nu_1} \) necessarily intersects \( \Delta \) 
since \( \nu_{1} \) is an essential return time for \(  
\omega^{(\nu_1)}   \subset \omega^{(\nu_0)} \). Therefore 
we get immediately 
\( |\overline\omega^{(\nu_0)}_{\nu_1}|\geq \delta^{\iota}/2 \) using 
simply the observation that $\overline{\omega}^{(\nu_0)}_{\nu_1}$ 
intersects both
$\partial\Delta^+$ and $\partial\Delta$. The final estimate follows as
in the previous paragraph. 
\end{proof}

\subsection{The counting argument}\label{combest}

\begin{lemma}\label{combestlem}
    For all \( 0\leq i
\leq n-1\), 
\( \omega\in \mathcal{Q}^{(i)}\) and 
\( R\geq r_{\delta} \), 
we have 
\[ 
\#
\mathcal{Q}^{(i+1)}(\omega, R) \leq e^{\gamma_{0} R}.  
\] 
\end{lemma} 

Before starting the proof, 
 recall that each 
\( \tilde\omega\in \mathcal{Q}^{(i+1)}(\omega, R) \) has a combinatorial
itinerary, associated to the iterates between the \( i \)'th escape
at time \( \eta_{i} \) and the \( i+1 \)'th escape at time \(
\eta_{i+1} \), 
specified by a sequence     
\begin{equation}\label{seq}
(\pm r_{1},
m_{1}), (\pm r_{2}, m_{2}), \dots, ( \pm r_{s}, m_{s}) 
\end{equation}
with 
\begin{equation}\label{seqcond}
s\geq 1, \quad |r_{1}|+\dots+|r_{s}|=R, \quad |r_{j}|\geq r_{\delta}, 
\quad m_{j}\in [1, r_{j}^{2}], \quad j=1, \dots, s.
\end{equation}  
The combinatorics is not unique but the multiplicity of partition
elements with the same combinatorics can be controlled. Indeed this is
essentially the main reason for introducing the sets \( \mathcal
Q^{(i)} \). 
\begin{sublemma}\label{multiplicity}
The cardinality of elements of \( \mathcal Q^{(i+1)}(\omega, R) \) having
the same combinatorial itinerary is at most \( r_{\delta}^{3} \).
 \end{sublemma}   
 \begin{proof}
The first time, after
     time \( \eta_{i} \), that \(\omega \) intersects \( \Delta^{+} \) in
     a chopping time, every subinterval which arises out of the `chopping'
     has either an escape time, in which case the sequence above is empty
  or an essential return time in which case
     a unique pair of integers \( r_{1} \) and \( m_{1} \) are associated
     to it.  Thus no two elements created up to this time can share the
     same sequence. Fixing one of these subintervals which has an essential
     return we consider higher iterates until the next time that it
     intersects \( \Delta^{+} \). At this time it is further subdivided
     into subintervals. Those which have essential returns at this time all
     have another uniquely defined pair of integers \( r_{2} \) and \(
     m_{2} \) associated to them. However multiplicity can occur for those
     which have escape times: all the subintervals which fall in \(
     \Delta^{+}\setminus\Delta \) have an escape at this time, and
     therefore belong to \( \mathcal Q^{(i)} \) and we do not consider
     further iterates, but all share the same first (and only) term of the
     associated sequence of return depths. However the number of such
     subintervals can be estimated by the number of elements of the
     partition \( \mathcal I^{+}|_{\Delta^{+}\setminus\Delta} \) plus at
     most two elements which can escape by falling outside \( \Delta^{+}
     \). The number of such intervals is then \( \leq 2
     (r_{\delta}-r_{\delta^+}) r_{\delta}^{2}+2 \leq r_{\delta}^{3}
     \). In the case of the intervals which have two or more returns we
     repeat the same reasoning to get the result. 
  \end{proof}   
  
     \begin{proof}[Proof of Lemma~\ref{combest}] 
Sublemma~\ref{multiplicity} reduces the proof of Lemma~\ref{combest}
to a purely combinatorial calculation of the cardinality of the set \( 
\mathcal S_{R} \) of all possible
sequences of the form \eqref{seq} satisfying the constraints given in 
\eqref{seqcond}. Let us denote by \( \mathcal S_{R}(s) \) the subset
of \( \mathcal S_{R} \) given by sequences with some fixed number \( s \)
of terms, and by \( \mathcal S_{R}^{+}(s) \) the subsequence of these 
given by considering only positive \( r_{j}'s \). Notice that the
crucial condition \( |r_{j}| \geq r_{\delta} \) implies that \( s\leq 
R/r_{\delta} \), and therefore we have
\begin{equation}\label{combest1}
\#\mathcal S_{R} \leq \sum_{s\leq R/r_{\delta}}  \#\mathcal S_{R}(s) 
\leq \sum_{s\leq R/r_{\delta}} 2^{s} \#\mathcal S_{R}^{+}(s). 
\end{equation}
To obtain a bound for \( \#\mathcal S_{R}^{+}(s)  \) observe first of 
all that for a given sequence \( (r_{1}, \ldots, r_{s}) \) the the
terms \( (m_{1}, \ldots, m_{s}) \) contribute an additional factor of 
\begin{equation}\label{mterms}
\prod_{j=1}^{s}r_{j}^{2} = e^{ 2 \sum_{i=1}^{s}\log r_{j}}=
e^{ 2 \sum_{i=1}^{s}\frac{\log r_{j}}{r_{j}} r_{j}} 
\leq e^{ 2 \sum_{i=1}^{s}\frac{\log r_{\delta}}{r_{\delta}} r_{j}} = 
e^{ \frac{2  \log r_{\delta}}{r_{\delta}} R}.
\end{equation}
It remains therefore only to estimates the number of possible
sequences \(  (r_{1}, \ldots, r_{s}) \) with \( r_{j}\geq r_{\delta} \)
and \( r_{1}+\dots + r_{s}=R \). This number is bounded above by 
\[ 
\begin{pmatrix} R-1 \\ s \end{pmatrix} \leq \begin{pmatrix} R \\ s \end{pmatrix} 
    = \frac{R!}{(R-s) !\ s!}.
\]
Indeed, consider a row of \( R  \) objects. The number of
ways of partitioning such a set into exactly \( s \) non-empty subsets
is equivalent to the number of ways of selecting \( s \) objects out
of the \( R-1 \) objects following the first object in the row.
Indeed, once such \( s \) objects have been chosen we can define the
partition as being formed by the consecutive objects which follow one 
of the choices (including the chosen object itself) until the object
preceding the next chosen object. 

Using Stirling's approximation formula for factorials 
\( k! \in [1, 1+\frac{1}{4k}] \sqrt{2\pi k}k^{k}e^{-k} \) we get 
\begin{equation}\label{rterms1} 
    \begin{pmatrix} R\\ s\end{pmatrix} 
    =
\frac{R !}{(R-s)! s !}  \leq 
\frac{R^{R}}{(R-s)^{R-s}s^{s}} 
=
\left(\frac{R}{R-s}\right)^{R} 
\left(\frac{R-s}{s}\right)^{s}.
\end{equation} 
To estimate the first term we use the fact that \(
s\leq R/r_{\delta}\) to get
\begin{equation}\label{rterms2}
\left(\frac{R}{R-s}\right)^{R}
\leq 
\left(\frac{R}{R-\frac{R}{r_{\delta}}}\right)^{R}
\leq
\left(\frac{R}{R(1-\frac{1}{r_{\delta}})}\right)^{R}
=
\left(1 -\frac{1}{r_{\delta}}\right)^{-R}
\leq 
e^{R \log (1+ \frac{2}{r_{\delta}})} \leq e^{\frac{2}{r_{\delta}} R}.  
\end{equation}
Notice that the last two inequalities follow from the
Taylor series \( (1-x)^{-1} = 1+x+ x^{2} + \dots \leq 1+2x \) with \(
x=1/r_{\delta} \) and using the fact that \( r_{\delta} \gg 2 \), and 
from the fact that   \( \log (1+x) < x \) for all \( x>0 \). 
To estimate the second term we write first of all 
\begin{equation}\label{rterms3}
\left(\frac{R-s}{s}\right)^{s} \leq \left(\frac{R}{s}\right)^{s}
\leq 
\left[\left(\frac{s}{R}\right)^{-\frac{s}{R}}\right]^{R} 
\leq
\left[\left(\frac{1}{r_{\delta}}\right)^{-\frac{1}{r_{\delta}}}\right]^{R} 
= e^{\frac{\log r_{\delta}}{r_{\delta}} R}. 
\end{equation}
The third inequality uses the fact that \( x^{-x} \) is monotonically decreasing 
to \( 1 \) as \( x\to 0 \) and \(
s/R\leq 1/r_{\delta} \). Now, substituting \eqref{rterms2} and
\eqref{rterms3} into \eqref{rterms1} and multiplying by \eqref{mterms}
gives 
\[ 
\#\mathcal S_{R}^{+}(s)  \leq e^{\frac{2+3\log
r_{\delta}}{r_{\delta}}R}.
\]
Substituting this into \eqref{combest1} we then get 
\[ 
\#\mathcal S_{R}  \leq \sum_{s\leq R/r_{\delta}} 2^{s} e^{\frac{2+3\log
r_{\delta}}{r_{\delta}}R} = 
\sum_{s\leq R/r_{\delta}} e^{\frac{2+\log 2+3\log
r_{\delta}}{r_{\delta}}R} = \frac{R}{r_{\delta}}e^{\frac{2+\log 2+3\log
r_{\delta}}{r_{\delta}}R}.
\]
Finally, taking into account the bound on the multiplicity of
intervals with the same combinatorics, given by Sublemma
\ref{multiplicity}, we get 
\begin{equation}\label{combest2}
\# \mathcal{Q}^{(i+1)}_n(\omega, R)  
\leq r_{\delta}^{3} \#\mathcal S_{R}
\leq  R r_{\delta}^{2} e^{\frac{2+\log 2+3\log
r_{\delta}}{r_{\delta}}R} = e^{\log R + 2 \log r_{\delta}} 
e^{\frac{2+\log 2+3\log
r_{\delta}}{r_{\delta}}R}.
\end{equation}
We now use the fact that 
\[ 
\log R = \frac{\log R}{R} R \leq \frac{\log r_{\delta}}{r_{\delta}} R 
\quad \text{and} \quad \log r_{\delta} = \frac{\log r_{\delta}}{R}R
\leq \frac{\log r_{\delta}}{r_{\delta}} R,
\]
and thus, substituting into \eqref{combest2}, we get 
\[ 
\#\mathcal S_{R}  \leq  e^{\frac{2+\log 2+5\log
r_{\delta}}{r_{\delta}}R}.
\]
By the definition of \( \gamma_{0}  \) in \eqref{gamma0} 
we obtain the result. 

\end{proof}


\begin{bibsection}
    \begin{biblist}
	\bib{AlvLuzPindim1}{article}{
	  author={Alves, Jos{\'e} F.},
	  author={Luzzatto, Stefano},
	  author={Pinheiro, Vilton},
	  title={Lyapunov exponents and rates of mixing for one-dimensional maps.},
	  journal={Ergodic Th. \& Dyn. Syst.},
	  volume={24},
	  year={2004},
	}
	\bib{ArbMat04}{article}{
	  author={Arbieto, Alexander},
	  author={Matheus, Carlos},
	  title={Decidability of Chaos for Some Families of Dynamical Systems},
	  journal={Foundations of Computational Mathematics},
	  volume={4},
	  pages={269--275},
	  date={2004},
	}
	\bib{AviMor03b}{article}{
	  author={Avila, Artur},
	  author={Gustavo Moreira, Carlos},
	  title={Bifurcations of unimodal maps},
	  booktitle={Dynamical systems. Part II},
	  series={Pubbl. Cent. Ric. Mat. Ennio Giorgi},
	  pages={1\ndash 22},
	  publisher={Scuola Norm. Sup.},
	  place={Pisa},
	  date={2003},
	}
	\bib{AviMor03}{article}{
	  author={Avila, Artur},
	  author={Moreira, Carlos Gustavo},
	  title={Statistical properties of unimodal maps: smooth families with negative Schwarzian derivative},
	  language={English, with English and French summaries},
	  note={Geometric methods in dynamics. I},
	  journal={Ast\'erisque},
	  number={286},
	  date={2003},
	  pages={xviii, 81\ndash 118},
	}
	\bib{AviLyuMel03}{article}{
	  author={Avila, Artur},
	  author={Lyubich, Mikhail},
	  author={de Melo, Welington},
	  title={Regular or stochastic dynamics in real analytic families of unimodal maps},
	  journal={Invent. Math.},
	  volume={154},
	  date={2003},
	  number={3},
	  pages={451\ndash 550},
	}
	\bib{BenCar85}{article}{
	  author={Benedicks, Michael},
	  author={Carleson, Lennart},
	  title={On iterations of $1-ax^2$ on $(-1,1)$},
	  date={1985},
	  journal={Ann. of Math.},
	  volume={122},
	  pages={1\ndash 25},
	}
	\bib{BenCar91}{article}{
	  author={Benedicks, Michael},
	  author={Carleson, Lennart},
	  title={The dynamics of the He\'non map},
	  date={1991},
	  journal={Ann. of Math.},
	  volume={133},
	  pages={73\ndash 169},
	}
	\bib{BruLuzStr03}{article}{
	  author={Bruin, Henk},
	  author={\href {www.ic.ac.uk/~luzzatto}{Stefano Luzzatto}},
	  author={\href {http://www.maths.warwick.ac.uk/~strien}{Sebastian van Strien}},
	  title={\href {http://dx.doi.org/10.1016/S0012-9593(03)00025-9}{Decay of correlations in one-dimensional dynamics}},
	  journal={Ann. Sci. \'Ec. Norm. Sup.},
	  volume={36},
	  number={4},
	  pages={621--646},
	  year={2003},
	}
	\bib{BruSheStr03}{article}{
	  author={Bruin, Henk},
	  author={Shen, Weixiao},
	  author={van Strien, Sebastian},
	  title={Invariant measures exist without a growth condition},
	  journal={Comm. Math. Phys.},
	  volume={241},
	  date={2003},
	  number={2-3},
	  pages={287\ndash 306},
	}
	\bib{ColEck83}{article}{
	  author={Collet, P.},
	  author={Eckmann, J.-P.},
	  title={Positive Lyapunov exponents and absolute continuity for maps of the interval},
	  journal={Ergodic Theory Dynam. Systems},
	  volume={3},
	  date={1983},
	  number={1},
	  pages={13\ndash 46},
	  issn={0143-3857},
	}
	\bib{ColEckLan80}{article}{
	  author={Collet, P.},
	  author={Eckmann, J.-P.},
	  author={Lanford, O. E., III},
	  title={Universal properties of maps on an interval},
	  journal={Comm. Math. Phys.},
	  volume={76},
	  date={1980},
	  number={3},
	  pages={211\ndash 254},
	}
	\bib{GraSwi97}{article}{
	  author={Graczyk, Jacek},
	  author={{\'S}wiatek, Grzegorz},
	  title={Generic hyperbolicity in the logistic family},
	  journal={Ann. of Math. (2)},
	  volume={146},
	  date={1997},
	  number={1},
	  pages={1\ndash 52},
	}
	\bib{Jak01}{article}{
	  author={Jakobson, M.~V.},
	  title={Piecewise smooth maps with absolutely continuous invariant measures and uniformly scaled Markov partitions},
	  journal={Proceedings of Symposia in Pure Mathematics},
	  volume={69},
	  pages={825--881},
	  date={2001},
	}
	\bib{Jak04}{article}{
	  author={Jakobson, M.~V.},
	  title={Parameter choice for families of maps with many critical points},
	  book={Modern Dynamical Systems and Applications},
	  publisher={Cambridge University Press},
	  date={2004},
	}
	\bib{Jak78}{article}{
	  author={Jakobson, M.~V.},
	  title={Topological and metric properties of one-dimensional endomorphisms},
	  journal={Sov. Math. Dokl.},
	  volume={19},
	  date={1978},
	  pages={1452--1456},
	}
	\bib{Jak81}{article}{
	  author={Jakobson, M.~V.},
	  title={Absolutely continuous invariant measures for one\ndash parameter families of one\ndash dimensional maps},
	  date={1981},
	  journal={Comm. Math. Phys.},
	  volume={81},
	  pages={39\ndash 88},
	}
	\bib{Kel90}{article}{
	  author={Keller, Gerhard},
	  title={Exponents, attractors and Hopf decompositions for interval maps},
	  journal={Ergodic Theory Dynam. Systems},
	  volume={10},
	  date={1990},
	  number={4},
	  pages={717\ndash 744},
	}
	\bib{KelNow92}{article}{
	  author={Keller, Gerhard},
	  author={Nowicki, Tomasz},
	  title={Spectral theory, zeta functions and the distribution of periodic points for Collet-Eckmann maps},
	  journal={Comm. Math. Phys.},
	  volume={149},
	  date={1992},
	  number={1},
	  pages={31\ndash 69},
	}
	\bib{Koz03}{article}{
	  author={Kozlovski, Oleg},
	  title={Axiom A maps are dense in the space of unimodal maps in the $C\sp k$ topology},
	  journal={Ann. of Math. (2)},
	  volume={157},
	  date={2003},
	  number={1},
	  pages={1\ndash 43},
	}
	\bib{KozSheStr04}{article}{
	  author={Kozlovski, Oleg},
	  author={Shen, Weixiao},
	  author={\href {http://www.maths.warwick.ac.uk/~strien}{Sebastian van Strien}},
	  title={Density of hyperbolicity in dimension one},
	  status={Preprint},
	  date={2004},
	  eprint={\url {http://www.maths.warwick.ac.uk/~strien/Publications/rigid4june.ps}},
	}
	\bib{Lan82}{article}{
	  author={Lanford, Oscar E., III},
	  title={A computer-assisted proof of the Feigenbaum conjectures},
	  journal={Bull. Amer. Math. Soc. (N.S.)},
	  volume={6},
	  date={1982},
	  number={3},
	  pages={427\ndash 434},
	}
	\bib{Luz05}{article}{
		  author={\href {http://www.ma.ic.ac.uk/~luzzatto}{Stefano Luzzatto}},
		  title={Stochastic behaviour in non-uniformly expanding maps},
		  journal={Handbook of Dynamical Systems},
		  publisher={Elsevier},
		  date={2005},
		}
		\bib{LuzTuc99}{article}{
		  author={\href {http://www.ma.ic.ac.uk/~luzzatto}{Stefano Luzzatto}},
		  author={\href {http://www.math.uu.se/~warwick/}{Tucker, Warwick}},
		  title={Non-uniformly expanding dynamics in maps with singularities and criticalities},
		  journal={\href {http://www.ihes.fr/IHES/Publications/Publications.html}{Inst. Hautes \'Etudes Sci. Publ. Math.}},
		  number={89},
		  date={1999},
		  pages={179\ndash 226},
		}
		\bib{LuzVia00}{article}{
		  author={\href {http://www.ma.ic.ac.uk/~luzzatto}{Stefano Luzzatto}},
		  author={\href {http://www.impa.br/~viana}{Viana, Marcelo}},
		  title={Positive Lyapunov exponents for Lorenz-like families with criticalities},
		  language={English, with English and French summaries},
		  note={G\'eom\'etrie complexe et syst\`emes dynamiques (Orsay, 1995)},
		  journal={Ast\'erisque},
		  number={261},
		  date={2000},
		  pages={xiii, 201\ndash 237},
		}
	\bib{KacMisMro04}{book}{
	  author={Kaczynski, Tomasz},
	  author={Mischaikow, Konstantin},
	  author={Mrozek, Marian},
	  title={Computational homology},
	  series={Applied Mathematical Sciences},
	  volume={157},
	  publisher={Springer-Verlag},
	  place={New York},
	  date={2004},
	  pages={xviii+480},
	}
	\bib{Mis02}{article}{
	  author={Mischaikow, Konstantin},
	  title={Topological techniques for efficient rigorous computation in dynamics},
	  journal={Acta Numer.},
	  volume={11},
	  date={2002},
	  pages={435\ndash 477},
	}
	\bib{ZglNis01}{article}{
	  author={Zgliczy{\'n}ski, Piotr},
	  author={Mischaikow, Konstantin},
	  title={Rigorous numerics for partial differential equations: the Kuramoto-Sivashinsky equation},
	  journal={Found. Comput. Math.},
	  volume={1},
	  date={2001},
	  number={3},
	  pages={255\ndash 288},
	}
	\bib{Lla04}{article}{
	  author={de la Llave, Rafael},
	  title={A tutorial on KAM theory},
	  status={Preprint},
	  date={2004},
	}
	\bib{Lor63}{article}{
	  title={Deterministic nonperiodic flow},
	  author={Lorenz, E. D.},
	  journal={J. Atmosph. Sci.},
	  volume={20},
	  pages={130\ndash 141},
	  date={1963},
	}
	\bib{Lyu02}{article}{
	  author={Lyubich, Mikhail},
	  title={Almost every real quadratic map is either regular or stochastic},
	  journal={Ann. of Math. (2)},
	  volume={156},
	  date={2002},
	  number={1},
	  pages={1\ndash 78},
	}
	\bib{Lyu97}{article}{
	  author={Lyubich, Mikhail},
	  title={Dynamics of quadratic polynomials. I, II},
	  journal={Acta Math.},
	  volume={178},
	  date={1997},
	  number={2},
	  pages={185\ndash 247, 247\ndash 297},
	  issn={0001-5962},
	}
	\bib{Man85}{article}{
	  author={Ma{\~n}{\'e}, Ricardo},
	  title={Hyperbolicity, sinks and measure in one-dimensional dynamics},
	  journal={Comm. Math. Phys.},
	  volume={100},
	  date={1985},
	  number={4},
	  pages={495\ndash 524},
	}
	\bib{MelStr93}{book}{
	  author={de Melo, Welington},
	  author={\href {http://www.maths.warwick.ac.uk/~strien}{Sebastian van Strien}},
	  title={One-dimensional dynamics},
	  series={Ergebnisse der Mathematik und ihrer Grenzgebiete (3) [Results in Mathematics and Related Areas (3)]},
	  volume={25},
	  publisher={Springer-Verlag},
	  place={Berlin},
	  date={1993},
	  pages={xiv+605},
	}
	\bib{Mis81}{article}{
	  author={Misiurewicz, Michal},
	  title={Absolutely continuous measures for certain maps of an interval},
	  journal={Inst. Hautes \'Etudes Sci. Publ. Math.},
	  number={53},
	  date={1981},
	  pages={17\ndash 51},
	}
	\bib{NowStr91}{article}{
	  author={Nowicki, Tomasz},
	  author={\href {http://www.maths.warwick.ac.uk/~strien}{Sebastian van Strien}},
	  title={Invariant measures exist under a summability condition for unimodal maps},
	  journal={Invent. Math.},
	  volume={105},
	  date={1991},
	  number={1},
	  pages={123\ndash 136},
	}
	\bib{Ryc88}{article}{
	  author={Rychlik, Marek Ryszard},
	  title={Another proof of Jakobson's theorem and related results},
	  journal={Ergodic Theory Dynam. Systems},
	  volume={8},
	  date={1988},
	  number={1},
	  pages={93\ndash 109},
	}
	\bib{SimTat91}{article}{
	  author={Sim{\'o}, C.},
	  author={Tatjer, J. C.},
	  title={Windows of attraction of the logistic map},
	  booktitle={European Conference on Iteration Theory (Batschuns, 1989)},
	  pages={335\ndash 342},
	  publisher={World Sci. Publishing},
	  place={River Edge, NJ},
	  date={1991},
	}
	\bib{SimTat05}{article}{
	  author={Sim{\'o}, C.},
	  author={Tatjer, J. C.},
	  title={personal communication},
	  date={2005},
	}
	\bib{ThiTreYou94}{article}{
	  author={Thieullen, Ph.},
	  author={Tresser, C.},
	  author={Young, L.-S.},
	  title={Positive Lyapunov exponent for generic one-parameter families of unimodal maps},
	  journal={J. Anal. Math.},
	  volume={64},
	  date={1994},
	  pages={121\ndash 172},
	}
	\bib{Thu99}{article}{
	  author={Thunberg, Hans},
	  title={Positive exponent in families with flat critical point},
	  journal={Ergodic Theory Dynam. Systems},
	  volume={19},
	  date={1999},
	  number={3},
	  pages={767\ndash 807},
	}
	\bib{Tsu93}{article}{
	  author={Tsujii, Masato},
	  title={Positive Lyapunov exponents in families of one-dimensional dynamical systems},
	  journal={Invent. Math.},
	  volume={111},
	  date={1993},
	  number={1},
	  pages={113\ndash 137},
	}
	\bib{Tuc02}{article}{
	  author={Tucker, Warwick},
	  title={A rigorous ODE solver and Smale's 14th problem},
	  journal={Found. Comput. Math.},
	  volume={2},
	  date={2002},
	  number={1},
	  pages={53\ndash 117},
	}
	\bib{Tuc99}{article}{
	  author={Tucker, Warwick},
	  title={The Lorenz attractor exists},
	  language={English, with English and French summaries},
	  journal={C. R. Acad. Sci. Paris S\'er. I Math.},
	  volume={328},
	  date={1999},
	  number={12},
	  pages={1197\ndash 1202},
	}
	\bib{UlaNeu47}{article}{
	  author={Ulam, S.},
	  author={von Neumann, J.},
	  title={On combination of stochastic and deterministic processes},
	  date={1947},
	  journal={Bull. AMS},
	  volume={53},
	  pages={1120},
	}
	\bib{Yoc01}{article}{
	  author={Yoccoz, Jean-Christophe},
	  title={Jakobson's Theorem},
	  status={Manuscript},
	  date={2001},
	}
	\bib{You92}{article}{
	  author={Young, L.-S.},
	  title={Decay of correlations for certain quadratic maps},
	  journal={Comm. Math. Phys.},
	  volume={146},
	  date={1992},
	  number={1},
	  pages={123\ndash 138},
	}

\end{biblist}
\end{bibsection}
\end{document}